\documentclass[12pt,centertags,oneside]{amsart}
\usepackage{amsmath,amstext,amsthm,amscd,typearea}
\usepackage{amssymb}
\usepackage{fancyhdr}
\usepackage[mathscr]{eucal}
\usepackage{mathrsfs}
\usepackage{concmath}    
\usepackage{charter} 
\usepackage{dsfont}  
\usepackage{pdfsync}
\usepackage[backref=page,bookmarks,bookmarksnumbered]{hyperref}
\usepackage{color}
\usepackage{comment}
\usepackage[a4paper,width=16.2cm,top=2.8cm,bottom=2.8cm]{geometry}

\numberwithin{equation}{section}
 
\newcommand{\field}[1]{\mathbb{#1}} 

\newcommand{\R}{\field{R}}
\newcommand{\C}{\field{C}}

\DeclareMathOperator{\Ker}{Ker}

\DeclareMathOperator{\Dom}{Dom}

\DeclareMathOperator{\rank}{rk}
\DeclareMathOperator{\Id}{Id}
\DeclareMathOperator{\supp}{supp} 
\DeclareMathOperator{\tr}{Tr}

\newtheorem{thm}{Theorem}[section]
\newtheorem{lemma}[thm]{Lemma}
\newtheorem{prop}[thm]{Proposition}
\newtheorem{cor}[thm]{Corollary}
\theoremstyle{definition}
\newtheorem{rem}[thm]{Remark}
\newtheorem{example}[thm]{Example}
\theoremstyle{definition}
\newtheorem{defn}[thm]{Definition}

\newcommand{\be}{\begin{eqnarray}}
\newcommand{\ee}{\end{eqnarray}}

\begin{document}
\pagestyle{fancy}
\lhead{} 
\chead{}  
\rhead{} 
\lfoot{}  
\cfoot{\thepage} 
\rfoot{}
\renewcommand{\headrulewidth}{0pt}
\renewcommand{\footrulewidth}{0pt} 
    
\title{Nakano-Griffiths inequality, holomorphic Morse inequalities,
and extension theorems for $q$-concave domains}

\author{Bingxiao Liu}
\address{Bingxiao Liu, Great Bay University, School of Sciences,
Dongguan 523000, China}
\address{Universit{\"a}t zu K{\"o}ln, Department Mathematik/Informatik,
Abteilung Mathematik, Weyertal 86--90, 50931 K{\"o}ln, Germany}
\email{liubx@gbu.edu.cn; bingxiao.liu@uni-koeln.de}

\author{George Marinescu}
\address{George Marinescu, Universit{\"a}t zu K{\"o}ln,
Department Mathematik/Informatik, Abteilung Mathematik,
Weyertal 86--90, 50931 K{\"o}ln, Germany}
\email{gmarines@math.uni-koeln.de}

\author{Huan Wang}
\address{Huan Wang, Institute of Mathematics, Henan Academy of Sciences,
Zhengzhou 450046, Henan, China}
\email{huanwang@hnas.ac.cn}

\thanks{B. Liu was partially supported by DFG Priority Program 2265
`Random Geometric Systems' (Project-ID 422743078) from Cologne and the
ANR--DFG project `QuasiDy\,--\,Quantization, Singularities, and Holomorphic
Dynamics' (ANR-21-CE40-0016) from Universit\'{e} de Lille (France).}

\thanks{G. Marinescu is partially supported by the DFG funded projects
SFB TRR 191 `Symplectic Structures in Geometry, Algebra and Dynamics'
(Project-ID 281071066\,--\,TRR 191), DFG Priority Program 2265 `Random
Geometric Systems' (Project-ID 422743078), and the ANR--DFG project
`QuasiDy\,--\,Quantization, Singularities, and Holomorphic Dynamics'
(Project-ID 490843120).}

\thanks{H. Wang is supported by the CSC Visiting Scholarship
(No.\ 202508410099), Natural Science Foundation of Henan Province of China
(Project No.\ 252300420324), Science and Technology Tackling Project of
Henan Province (Project No.\ 262102520039), High-level Talent Research
Start-up Project Funding of Henan Academy of Sciences
(Project No. 241819110), the Scientific and Technological Research Project
of Henan Academy of Sciences (Project No. 20262319005), Guest Program of
MPIM Bonn, and IMU Simons Research Fellowship for Developing Countries
(Identifier: 2290).}

\date{\today}

\begin{abstract}

We establish a general Nakano--Griffiths inequality with boundary conditions
and apply it to derive holomorphic Morse inequalities for domains satisfying
analytic convexity assumptions. The proof of these inequalities relies
on analyzing the spectral spaces of the Laplace operator with
$\overline{\partial }$-Neumann boundary conditions. As an application, we
obtain a criterion for Moishezon 1-concave manifolds.

We further apply the holomorphic Morse inequalities on Levi $q$-concave
domains in conjunction with the Kohn--Rossi extension theorem to obtain
the following result. Let $X$ be a compact complex manifold of dimension
$n $ equipped with a holomorphic line bundle that is semi-positive everywhere
and positive at least at one point, and let $D\subset X$ be a smooth
$q$-concave domain ($1 \le q \le n-1$). We prove that every
$\bar{\partial }_{b}$-closed $(0,\ell )$-form on $bD$ with values in a holomorphic
vector bundle, admits a meromorphic extension to $D$ for all
$q \le \ell \le n-1$.

\bigskip

\noindent{{\sc Mathematics Subject Classification 2020}:
32W10, 32F10, 32D15, 32A25.}

\smallskip

\noindent{{\sc Keywords}: Kohn--Rossi extension, holomorphic Morse inequalities,
$q$-concave manifold, semi-positive line bundle}
\end{abstract}
\maketitle

\tableofcontents

\section{Introduction}
\label{intro}

A foundational tool in complex geometry is the Bochner--Kodaira--Nakano
identity and its variants, which have far-reaching implications, from vanishing
theorems to holomorphic Morse inequalities and extension results. In this
paper, we develop a general version of the Nakano--Griffiths inequality
with boundary conditions, unifying the pseudoconvex/positive and pseudoconcave/negative
settings. This is achieved via a refined Bochner--Kodaira--Nakano formula
that incorporates boundary terms in full generality, and extends the classical
results of Andreotti--Vesentini \cite{AV65}, Griffiths
\cite{Griffiths1966}, and their analytic extensions in \cite{MM} (see {Theorems~\ref{thm:BKN-variant} and \ref{thm:NG-inequality}}).

Our new inequality leads naturally to holomorphic Morse inequalities for
domains and manifolds satisfying analytic convexity conditions. Rather
than relying on scaling arguments, we adopt a geometric approach based
on the spectral theory of the $\overline{\partial}$-Laplacian with
$\overline{\partial}$-Neumann boundary conditions. Our method provides
new index ranges and accommodates weaker curvature assumptions than those
previously required. In particular, we establish holomorphic Morse inequalities
for Levi $q$-concave domains (see {Theorem~\ref{thm:4.10important}}) and
for $q$-concave complex manifolds (see {Theorem~\ref{TheoremB}}), refining
and extending earlier results \cite{RB:05,MM,M:96}. As an application
of this analytic framework, we obtain a geometric criterion for Moishezon
1-concave manifolds.

We conclude with an application to the extension problem for
$\bar\partial _{b}$-closed forms on the boundary of pseudoconcave domains,
where no general vanishing theorem is available. This question has its
roots in classical results such as Hartogs' theorem
\cite[Theorem 2.3.2]{Hor:90}, which ensures that holomorphic functions
on $M \setminus K \subset \mathbb{C}^{n}$, $n \ge 2$, extend to $M$ when
$K \subset M$ is compact and $M \setminus K$ is connected. Bochner
\cite{Bo:43} refined this by allowing the extension of functions satisfying
the tangential Cauchy--Riemann equations on the boundary $bM$. Further
generalizations to complex manifolds were given in
\cite{AH:72,TL:89}.

A key method, introduced by Ehrenpreis \cite{Ehr:56}, reduces such extension
problems to solving the $\bar\partial $-equation with compact support,
interpreting obstructions via $\bar\partial $-cohomology. This has been
developed in \cite{And:73,AH:72,Hor:90,KR:65,TL:89}.

In our setting, we study the extension of $\bar\partial _{b}$-closed forms
with values in a holomorphic vector bundle from the boundary of a pseudoconcave
domain. Unlike the case of holomorphic functions, this requires more refined
tools. We combine holomorphic Morse inequalities for $q$-concave domains
with the Kohn--Rossi extension theorem to obtain precise results. Specifically,
let $X$ be a compact complex manifold of dimension $n$, equipped with a
holomorphic line bundle that is semi-positive everywhere and strictly positive
at one point, and let $D \subset X$ be a smooth $q$-concave domain with
$1 \le q \le n-1$. We prove that every $\bar\partial _{b}$-closed
$(0,\ell )$-form on $bD$ with values in a holomorphic vector bundle, admits
a meromorphic extension to $D$ for all $q \le \ell \le n-1$ (see {Theorem~\ref{thmextension1}}). Our method also yields new vanishing results, including
for domains with Levi-flat boundaries.

\subsection{Nakano-Griffiths inequality with boundary terms}
\label{sec1.1}

Our approach to the holomorphic Morse inequalities is based on a novel
Bochner-Kodaira-Nakano formula with boundary terms. Let
$(X,\Theta )$ be a Hermitian manifold of dimension $n$, and
$M:=\{x\in X:\varrho (x)<0\}$ be a relatively compact domain with a smooth
defining function $\varrho $ such that $|d\varrho |=1$ on $bM$. Let
$\mathscr{L}_{\varrho}$ denote the Levi form for $bM$ (see {\eqref{eq:2.2Levi}}). For a holomorphic vector bundle $F \to X$ with a smooth
Hermitian metric $h^{F}$, let $\nabla ^{F}$ denote the associated Chern
connection and $R^{F}=(\nabla ^{F})^{2}$ the Chern curvature. We refer
to Section~\ref{sec-NG} for the notations used in the following theorem.

\begin{thm}[A Bochner-Kodaira-Nakano formula with boundary terms]
\label{thm:BKN-variant}
Let $M\Subset X$ be given as above. Let
$\tau \in \mathscr{C}^{\infty}(X,[0,1])$ be such that $\tau =1$ near
$bM$ and $\supp \tau $ is contained in an open neighborhood of $bM$ where
$d\varrho $ does not vanish. Fix $\varepsilon \in [0,1]$. Then for
$s\in B^{0,q}(M,F)$, we have
\begin{equation}
\begin{split} &\|\overline \partial ^{F} s\|^{2}_{L^{2}(M)}+\|
\overline \partial ^{F*} s\|^{2}_{L^{2}(M)}
\\
= &\left \langle [\sqrt{-1}R^{F},\Lambda ]_{\varepsilon} s,s\right
\rangle _{L^{2}(M)} +2(1-\varepsilon )\left \langle \tau R^{F}(e^{(1,0)}_{
\mathfrak{n}}, e^{(0,1)}_{\mathfrak{n}})s,s\right \rangle _{L^{2}(M)}\\
&+
\varepsilon \sum _{j=1}^{n}\|\widetilde{\nabla }^{TX}_{
\overline{\omega }_{j}}s\|^{2}_{L^{2}(M)}
\\
&+\int _{bM}\langle [\sqrt{-1}\mathscr{L}_{\varrho},\Lambda ]^{bM}_{
\varepsilon}s,s\rangle _{h} dv_{bM}+ (1-\varepsilon )\int _{bM} e_{
\mathfrak{n}}(|d\varrho |)|s|^{2}_{h} dv_{bM}
\\
&+\big\langle \sum ^{n}_{j,k=1} R^{K^{\ast}_{X}}(\omega _{j},
\overline \omega _{k})\overline{\omega }^{k}\wedge \iota _{
\overline{\omega }_{j}}s,s\big\rangle _{L^{2}(M)}-(1-\varepsilon )
\int _{M} |s|^{2}_{h} \sqrt{-1}\partial \overline{\partial }\Omega _{
\tau}
\\
&+Q^{F}_{T,\varepsilon}(s,s)+(1-\varepsilon )\left (\Psi ^{F}_{\tau}(s,s)+
2\|\sqrt{\tau}\,\widetilde{\nabla}^{TX}_{e^{(0,1)}_{\mathfrak{n}}}s\|^{2}_{L^{2}(M)}
\right )
\\
&+(1-\varepsilon )\left (\sum _{j=1}^{n}\|\widetilde{\nabla }^{TX}_{
\omega _{j}}s\|^{2}_{L^{2}(M)} -2\|\sqrt{\tau}\,\widetilde{\nabla}^{TX}_{e^{(1,0)}_{
\mathfrak{n}}}s\|^{2}_{L^{2}(M)}\right ).
\end{split}
\label{eq:3.42Nov}
\end{equation}
\end{thm}
The proof of {\eqref{eq:3.42Nov}} is based on the Bochner-Kodaira-Nakano
formula with boundary terms \cite[Theorem 1.4.21]{MM} (cf.\, {Theorem~\ref{thm-bkn-mm}}), which uses a Riemannian formalism and the approach of
Griffiths \cite[Theorem 7.4]{Griffiths1966}, based on the existence of
the $(n-1,n-1)$-form $\Omega _{\tau}$, specific to Hermitian manifolds.
For $\varepsilon =1$, formula {\eqref{eq:3.42Nov}} reduces to the Bochner-Kodaira-Nakano
formula with boundary terms {\eqref{eq:3.17Nov}} by Andreotti and Vesentini
\cite[p.\ 113]{AV65} and by Griffiths \cite[(7.14)]{Griffiths1966}, which
in turn follows from \cite[Theorem 1.4.21]{MM}.

An important feature of {\eqref{eq:3.42Nov}} is the introduction of the modified
curvature commutator
$[\sqrt{-1}R^{F},\Lambda ]_{\varepsilon} := [\sqrt{-1}R^{F},\Lambda ]+
\varepsilon \mathrm{Tr}_{\Theta}[\sqrt{-1}R^{F}]$ (cf.\ {\eqref{eq:MLB1}})
and similarly for the Levi form,
$[\sqrt{-1}\mathscr{L}_{\varrho},\Lambda ]^{bM}_{\varepsilon}$ (cf.\ {\eqref{eq:MLB2}}).
The curvature commutator $[\sqrt{-1}R^{F},\Lambda ]$ plays an important
role in the applications of the Bochner-Kodaira-Nakano formula
\cite{Demailly1985,Dem1986,MM}, such as vanishing theorems, holomorphic
Morse inequalities, or Bergman kernel asymptotics. In {\eqref{eq:3.42Nov}}, the main difficulty is to identify an operator of the
same nature on the boundary $bM$, namely
$[\sqrt{-1}\mathscr{L}_{\varrho},\Lambda ]^{bM}_{\varepsilon}$, where we
need to develop new techniques inspired by the remarkable calculations
from \cite[VII. \S 5]{Griffiths1966} used to derive a Nakano-Griffiths
inequality \cite[(7.8)]{Griffiths1966}; see Subsection \ref{ss:3.2variant}.

If $F$ is a line bundle and $R^{F}$ is negative, then
$[\sqrt{-1}R^{F},\Lambda ]_{\varepsilon}$ is a positive operator on
$(0,q)$-forms ($q\leq n-1$) for sufficiently small
$\varepsilon \geq 0$; if $\varepsilon =1$ (or $\varepsilon $ is sufficiently
close to $1$), $[\sqrt{-1}R^{F},\Lambda ]_{\varepsilon}$ becomes negative;
similarly for the term
$[\sqrt{-1}\mathscr{L}_{\varrho},\Lambda ]^{bM}_{\varepsilon}$, see {\eqref{eq:4.12dec24}}. This suggests that a formula like {\eqref{eq:3.42Nov}} is needed to deal with negative bundles on a domain
with a pseudoconcave boundary, parallel to a formula such as {\eqref{eq:3.17Nov}} adapted to the situation of positive bundles on a domain
with pseudoconvex boundary. Our formula {\eqref{eq:3.42Nov}} actually unifies
these two scenarios.

Since we do not assume $\Theta $ to be K\"{a}hler, the torsion tensor
$T$ (see {\eqref{eq:3.5torsion}}) appears in {\eqref{eq:3.42Nov}}, and is encoded
in the term $Q^{F}_{T,\varepsilon}(s,s)$ (see {Definition~\ref{def:3.7QFT}}). If $(X,\Theta )$ is K\"{a}hler, or more generally, if
$d\Theta =0$ in a neighborhood of $M$, then $T=0$ in $\overline{M}$, and
we have $Q^{F}_{T,\varepsilon}(s,s)=0$.

As a consequence of the equality {\eqref{eq:3.42Nov}}, we obtain a generalization
of \cite[(7.8)]{Griffiths1966}.

\begin{thm}[Nakano-Griffiths inequality with boundary terms]%
\label{thm:NG-inequality}
Let $M=\{\varrho <0\}\Subset X$ be a relatively compact domain with a smooth
boundary $bM$. Let $\tau \in \mathscr{C}^{\infty}(X,[0,1])$ be such that
$\tau =1$ near $bM$ and $\supp \tau $ is contained in an open neighborhood
of $bM$ where $d\varrho $ does not vanish. Then there exists a constant
$C=C(X,g^{TX},M,\varrho ,\tau )>0$, which is independent of
$(F,h^{F})$, such that, for any $\varepsilon \in (0,1]$ and any
$s\in B^{0,q}(M,F)$, we have
\begin{equation}
\begin{split} &(1+\varepsilon )\left (\|\overline \partial ^{F} s\|^{2}_{L^{2}(M)}+
\|\overline \partial ^{F*} s\|^{2}_{L^{2}(M)}\right )+C(1+
\frac{1}{\varepsilon})\|s\|^{2}_{L^{2}(M)}
\\
&\quad \geq \,\left \langle [\sqrt{-1}R^{F},\Lambda ]_{\varepsilon} s,s
\right \rangle _{L^{2}(M)}+2(1-\varepsilon )\left \langle \tau R^{F}(e^{(1,0)}_{
\mathfrak{n}}, e^{(0,1)}_{\mathfrak{n}})s,s\right \rangle _{L^{2}(M)}
\\
&\qquad \qquad +\int _{bM}\langle [\sqrt{-1}\mathscr{L}_{\varrho},
\Lambda ]^{bM}_{\varepsilon}s,s\rangle _{h} dv_{bM}+(1-\varepsilon )
\int _{bM} e_{\mathfrak{n}}(|d\varrho |)|s|^{2}_{h} dv_{bM}
\\
&\qquad \qquad +(1-\varepsilon )\left (\|(\widetilde{\nabla}^{TX})^{1,0}s
\|^{2}_{L^{2}(M)}-2\|\sqrt{\tau}\widetilde{\nabla}^{TX}_{e^{(1,0)}_{
\mathfrak{n}}}s\|^{2}_{L^{2}(M)}\right ).
\end{split}
\label{eq:3.44Nov}
\end{equation}
\end{thm}

Formula {\eqref{eq:3.44Nov}} is more involved than
\cite[(7.8)]{Griffiths1966}. Indeed, in \cite{Griffiths1966}, it is assumed
that $|d\varrho |=1$ near $bM$ (not only on $bM$), so that
$e_{\mathfrak{n}}(|d\varrho |)\equiv 0$ on $bM$ and the second boundary
term (the integral on $bM$) in {\eqref{eq:3.44Nov}} do not appear. The
last line in {\eqref{eq:3.44Nov}} will always be non-negative. In
\cite{Griffiths1966}, a modified Hermitian metric on $F$ near $bM$ is used
to derive \cite[(7.8)]{Griffiths1966} (see
\cite[the last part of VII. \S 6]{Griffiths1966}). This is not needed in {\eqref{eq:3.44Nov}}.

When $\varepsilon =1$ in {\eqref{eq:3.44Nov}}, we get
\cite[Corollary 1.4.22]{MM} (up to constants), which can be used to prove
the holomorphic Morse inequalities for $q$-convex manifolds or domains
(see \cite[Section 3.5]{MM}). However, for the $q$-concave case, we need
to apply {\eqref{eq:3.44Nov}} for sufficiently small $\varepsilon >0$, in
particular, to obtain the optimal fundamental estimates in {Proposition~\ref{prop:4.3Jan25}}; see Subsection \ref{ss4.2Jan25}. In Subsection \ref{ss:4.5pqcorona}, we combine these two set-ups to obtain the weak holomorphic
Morse inequalities for $(p,q)$-coronas and $(p,q)$-convex-concave manifolds
(see {Definition~\ref{def:pqcorona}}).

Let $H^{0,\ell}(\overline{M},F)$ denote the Dolbeault cohomology groups
up to the boundary, defined by the Dolbeault complex consisting of differential
forms on $M$ that are smooth up to the boundary; see Subsection \ref{ss:cohomology}. From \eqref{eq:3.44Nov}, we can also recover the classical
Andreotti-Tomassini vanishing theorem for both $q$-concave and $q$-convex
manifolds; see Subsection \ref{ss:6.1AT}. In fact, we obtain a more general
vanishing theorem as follows. For a Hermitian line bundle $L\to X$ and
an integer $k\in \mathbb{N}$, $L^{k}:=L^{\otimes k}$ denotes the $k$-th
tensor power of $L$.

\begin{thm}%
\label{thm_vanishing}
Let $X$ be a complex manifold of dimension $n\geq 2$. Let
$1\leq q\leq n-1$, and let $M\Subset X$ be a smooth domain in $X$. Let
$E$ be a holomorphic vector bundle on $X$.

(a) Assume that the Levi form $\mathscr{L}_{bM}$ has at least $n-q$ negative
eigenvalues on $bM$. Let $(L,h^{L})$ be a Hermitian line bundle on
$X$ such that $c_{1}(L,h^{L})<0$ on $\overline{M}$. Then, there exists
an integer $k_{0}=k_{0}(L,E)$ such that
\begin{equation}
H^{0,\ell}(\overline{M},L^{k}\otimes E)=0 \quad \,\text{ for }\, k
\geq k_{0}\,\text{ and }\,\ell \leq n-q-1.
\label{eq1.3}
\end{equation}

(b) Assume that the Levi form $\mathscr{L}_{bM}$ has at least $n-q$ positive
eigenvalues on $bM$. Let $(L,h^{L})$ be a Hermitian line bundle on
$X$ such that $c_{1}(L,h^{L})>0$ on $\overline{M}$. Then, there exists
an integer $k_{0}=k_{0}(L,E)$ such that
\begin{equation}
H^{0,\ell}(\overline{M},L^{k}\otimes E)=0 \quad \,\text{ for }\, k
\geq k_{0}\,\text{ and }\,\ell \geq q.
\label{eq1.4}
\end{equation}
The cohomology groups $H^{0,\ell}(\overline{M},L^{k}\otimes E)$ in the
above statements can be replaced by the spaces of harmonic forms
$\mathscr{H}^{0,\ell}(M,L^{k}\otimes E)$ with respect to an arbitrary Hermitian
metric on $E$.
\end{thm}

Finally, in Subsection \ref{ss:6.3Leviflat}, employing {Theorem~\ref{thm:NG-inequality}, we derive a vanishing theorem for negative line
bundles on relatively compact domains of complex dimension $n\geq 2$ with
a smooth Levi-flat boundary.

\subsection{Holomorphic Morse inequalities for \texorpdfstring{$q$}{q}-concave domains}
\label{sec1.2}

Let $F$ be a holomorphic vector bundle on the complex manifold $X$. We denote by
$H^{\ell}(X,F)$ the $\ell $-th sheaf cohomology group of $X$ with values
in $F$. Let $M\Subset X$ be a relatively compact domain with a smooth boundary.
Given the smooth metrics on $F$ and $X$, we can consider the Kodaira
Laplacians on $M$ with $\overline \partial $-Neumann boundary conditions.
We denote by $\mathscr{H}^{0,\ell}(M,F)$ the corresponding space of harmonic
forms. For a Hermitian line bundle $(L,h^{L})$, let $R^{L}$ be the Chern
curvature form and let
$c_{1}(L,h^{L}):=\frac{\sqrt{-1}}{2\pi} R^{L}$ be the first Chern form
associated with $h^{L}$. We denote by
$(L^{k}, h^{L^{k}}):=(L^{\otimes k}, (h^{L})^{\otimes k})$ the $k$-th tensor
power.

The formula {\eqref{eq:3.44Nov}}, together with suitable geometric conditions,
implies the optimal fundamental estimates for the Kodaira Laplacians with
$\overline \partial $-Neumann boundary conditions on $\overline{M}$; therefore,
{Theorems~\ref{thm:4.10important} and \ref{thm:4.13new}} follow from the
abstract holomorphic Morse inequalities for $L^{2}$ cohomology; see Subsection \ref{sub:2.3abstract} and \cite[Section 3.2]{MM}. We have the following
holomorphic Morse inequalities for the (Levi) $q$-concave domains.

\begin{thm}%
\label{thm:4.10important}
Let $M\Subset X$ be a smooth domain of a Hermitian manifold
$(X,\Theta )$ of dimension $n$ such that the Levi form
$\mathscr{L}_{bM}$ (see {\eqref{eq:2.2Levi}}) has at least $n-q$ negative
eigenvalues, $1\leq q\leq n-1$. Let $(E, h^{E})$ and $(L,h^{L})$ be holomorphic
Hermitian vector bundles on $X$ with $\rank (L)=1$.

If $c_{1}(L,h^{L})\leq 0$ on $M\setminus K_{L}$ for some compact subset
$K_{L}\subset M$, then for $\ell \leq n-q-1$,
\begin{equation}
\begin{split} \limsup _{k\rightarrow \infty}k^{-n} \dim \mathscr{H}^{0,
\ell}(M,L^{k}\otimes E)&\leq \frac{\rank (E)}{n!}\int _{K_{L}(\ell ,h^{L})}(-1)^{
\ell}c_{1}(L,h^{L})^{n},
\\
\limsup _{k\rightarrow \infty}\sum _{j=0}^{\ell }(-1)^{\ell -j}k^{-n}
\dim \mathscr{H}^{0,j}(M,L^{k}\otimes E)&\leq \frac{\rank (E)}{n!}
\int _{K_{L}(\leq \ell \,, h^{L})}(-1)^{\ell}c_{1}(L,h^{L})^{n},
\end{split}
\label{eq:4.58Jan25-1}
\end{equation}
where $K_{L}(\ell ,h^{L})$ is the subset of $K_{L}$ consisting of all points
$x$ such that $\sqrt{-1}R^{L}_{x}$ is non-degenerate and has exactly
$\ell $ negative eigenvalues, and
$K_{L}(\leq \ell ,h^{L})=\cup _{0\leq j \leq \ell } K_{L}(j,h^{L})$.

Moreover, if $1\leq q\leq n-2$, we also have
\begin{equation}
\liminf _{k\rightarrow \infty}k^{-n} \dim \mathscr{H}^{0,0}(M,L^{k}
\otimes E)\geq \frac{\rank (E)}{n!}\int _{K_{L}(\leq 1\,,h^{L})} c_{1}(L,h^{L})^{n}.
\label{eq:4.58Jan25-2}
\end{equation}
The spaces $\mathscr{H}^{0,\ell}(M,L^{k}\otimes E)$ in {\eqref{eq:4.58Jan25-1}} and {\eqref{eq:4.58Jan25-2}} can be replaced by the
cohomology groups $H^{0,\ell}(\overline{M},L^{k}\otimes E)$.
\end{thm}

When $E=\mathbb{C}$, the Morse inequalities {\eqref{eq:4.58Jan25-1}} for
$0\leq \ell \leq n-q-1$ were also implicitly given by Berman in
\cite[Theorems 6.5 and 6.6]{RB:05}, via the scaling method for Bergman
kernels. Our method develops the Bochner techniques along Griffiths
\cite{Griffiths1966} and adapts the abstract holomorphic Morse inequalities
for $L^{2}$ cohomology as in \cite{MM}, which also leads to a more geometric
approach for the classical Andreotti-Tomassini vanishing theorems
\cite{AT:69}; see {Theorems~\ref{thm_vanishing} and \ref{vanish}}.

The holomorphic Morse
inequalities for $q$-concave manifolds, established in
\cite{M:96,M:2016-q}, are a consequence of {Theorem~\ref{thm:4.10important}}. In this paper, our method is based on the fundamental
estimates for the Kodaira Laplacians with $\overline\partial $-Neumann
boundary conditions.

\begin{thm}[{\cite{M:96}}]%
\label{TheoremB}
Let $X$ be a $q$-concave complex manifold of dimension $n\geq 2$, where
$1\leq q \leq n-1$. Let $E,L$ be holomorphic vector bundles on $X$ with
$\rank (L)=1$. Let $h^{L}$ be a smooth Hermitian metric on $L$ such that
$c_{1}(L,h^{L})\leq 0$ on $X\setminus K_{L}$ for some compact subset
$K_{L}$. Then for $j\leq n-q-1$, we have the weak Morse inequalities: as
$k\rightarrow \infty $,
\begin{equation}
\dim H^{j}(X,L^{k}\otimes E)\leq \rank (E)\frac{k^{n}}{n!}\int _{K_{L}(j,h^{L})}(-1)^{j}c_{1}(L,h^{L})^{n}+o(k^{n}).
\label{eq:1.4Jan25}
\end{equation}
For $j\leq n-q-2$, we have the strong Morse inequalities: as
$k\rightarrow \infty $,
\begin{equation}
\sum _{\ell =0}^{j}(-1)^{j-\ell}\dim H^{\ell}(X,L^{k}\otimes E)\leq
\rank (E)\frac{k^{n}}{n!}\int _{K_{L}(\leq \,j,h^{L})}(-1)^{j}c_{1}(L,h^{L})^{n}+o(k^{n}).
\label{eq:1.3Jan25}
\end{equation}
Moreover, if $1\leq q\leq n-2$, we also have
\begin{equation}
\dim H^{0}(X,L^{k}\otimes E)\geq \rank (E)\frac{k^{n}}{n!}\int _{K_{L}(
\leq 1\,,h^{L})} c_{1}(L,h^{L})^{n}+o(k^{n}).
\label{eq:1.5Jan25}
\end{equation}
\end{thm}
\begin{rem}
\label{rem1.6}
The proof of the above theorem will be given in Subsection \ref{ss:4.3Jan25}. Note that in \cite{M:96}, the inequality {\eqref{eq:1.4Jan25}} was stated only for all $j\leq n-q-2$; here we can
state it for all $j\leq n-q-1$ using our approach. Additionally, we would
like to make the following remarks:

(i) Let $\varrho \in \mathscr{C}^{\infty}(X,\mathbb{R})$ be a sublevel
exhaustion function for $X$ with exceptional set $K$ (see {Definition~\ref{defmfd}}) in {Theorem~\ref{TheoremB}}, and set
$X_{c}:=\{x\in X\;:\;\varrho (x)<c\}$ for $c\in \mathbb{R}$. We can modify
the assumption from
\textit{$c_{1}(L,h^{L})\leq 0$ on $X\setminus K_{L}$} to the more relaxed
condition \textit{$c_{1}(L,h^{L})\leq 0$ on $X_{c}\setminus K_{L}$}, where
$c$ is a regular value of $\varrho $ with
$K\cup K_{L}\Subset X_{c}\Subset X$.

(ii) If $X$ is a compact complex manifold, then we can take
$K_{L}=X$, and all the assumptions in {Theorem~\ref{TheoremB}} hold trivially,
so that the inequalities {\eqref{eq:1.4Jan25}} and {\eqref{eq:1.3Jan25}} are
exactly the holomorphic Morse inequalities proven by Demailly in
\cite[Th\'{e}or\`{e}me 0.1]{Demailly1985}.

(iii) When $X$ is noncompact, we further develop the ideas in
\cite[Chapter 3]{MM} and give a different approach to \cite{M:96}. Holomorphic
Morse inequalities for noncompact manifolds were proved in various contexts;
by Nadel-Tsuji \cite{NT88} for some complete K\"{a}hler manifolds with
Ricci negative curvature; by Bouche \cite{Bouche1989} for $q$-convex manifolds;
by Marinescu \cite{Mar92b,M:96,TCM:01} for $q$-concave manifolds, weakly
$1$-complete manifolds, and covering manifolds; and by Berman
\cite{RB:05} for smooth domains satisfying the $Z(q)$ condition. 
The paper \cite{HMZ26} presents an approach based on the heat kernel
to the Morse inequalities for domains satisfying the $Z(q)$ condition.
In this
paper, parallel to Ma-Marinescu \cite[Chapter 3]{MM}, we present a unified
treatment. The Bochner-Kodaira-Nakano-type formulas with boundary terms
({Theorems~\ref{thm:BKN-variant} and \ref{thm:NG-inequality}}) play a fundamental
role in this treatment. For further details, we refer to Subsection \ref{ss:4.5pqcorona}.
\end{rem}

\subsection{Characterization of
Moishezon manifolds and compactification}
\label{sec1.3}

We further apply the Bochner-Kodaira-Nakano formula with a boundary term
to obtain asymptotic lower bounds for the dimension of the spaces of holomorphic
sections of a semi-positive line bundle over a $1$-concave manifold. They
generalize the well-known lower bounds obtained by Siu \cite{Siu84} and
Demailly \cite{Demailly1985} in order to solve the Grauert-Riemenschneider
conjecture, which provides a characterization of Moishezon manifolds in
terms of semi-positive line bundles. The corresponding problem for
$1$-concave manifolds was considered in
\cite{RB:05,HLM22,M:96,M:2016-q}.

Recall that an $n$-dimensional compact complex manifold $M$ or a $q$-concave
complex manifold $M$ ($1\leq q\leq n-2$) is called
\textit{Moishezon} if $\alpha (M)=n$, where $\alpha (M)$ denotes the transcendence
degree of the field of meromorphic functions on $M$ (see
\cite[Definition 2.2.12 and Subsection 3.4]{MM} and \cite[\S 5]{M:96}).
This definition also works for a compact irreducible complex space (see
\cite[Example 3.4.2(ii)]{MM}). According to a classical theorem of Moishezon,
a compact Moishezon manifold admits a proper modification
$\widehat{X}\to X$, given by a finite number of blow-ups with smooth centers,
such that $\widehat{X}$ is a projective algebraic manifold (see
\cite[Theorem 2.2.16]{MM}).

The following criterion for Moishezon $q$-concave manifolds is essentially
derived from \cite[Theorem 1.1]{M:96} (or equivalently,
{\eqref{eq:1.5Jan25}}) in combination with the above definition and
\cite[Theorem 3.4.7]{MM}.
\begin{thm}
\label{thm1.7}
Let $X$ be a $q$-concave complex manifold of dimension $n\geq 3$, where
$1\leq q \leq n-2$. Let $L$ be a holomorphic line bundle on $X$ with a
smooth Hermitian metric $h^{L}$ such that $c_{1}(L,h^{L})\leq 0$ on
$X\setminus K_L$ for some compact subset $K_L$. If
\begin{equation}
\int _{K_{L}(\leq 1\,,h^{L})} c_{1}(L,h^{L})^{n} >0,
\label{eq:1.5Jan25Moishezon}
\end{equation}
then $X$ is Moishezon.
\end{thm}
A holomorphic line bundle $L\to X$ on a $q$-concave manifold $X$ is called
big if
\begin{equation*}
\limsup _{k\rightarrow \infty}k^{-n}\dim H^{0}(X, L^{k})>0.
\end{equation*}
By {\eqref{eq:1.5Jan25}}, the criterion {\eqref{eq:1.5Jan25Moishezon}} also
gives the bigness of $L$. In this case, for $k$ sufficiently large, the
sections of $H^{0}(X, L^{k})$ provide local coordinates on a Zariski open
set of $X$.

Now let us focus on the $1$-concave case, also known as the strictly pseudoconcave
case. As we mentioned in the beginning of this part, we would like to study
the semi-positive line bundle $(L,h^{L})$.

\begin{thm}%
\label{thm:4.13new}
Let $X$ be a connected complex manifold of dimension $n\geq 3$. Let
$M$ be a relatively compact domain in $X$ with a smooth boundary
$bM$, and assume that $\varrho $ is a defining function for $bM$ with
$M=\{x\in X:\varrho (x)<0\}$. Assume that the Levi form
$\mathscr{L}_{\varrho}$ is strictly negative on $bM$. Let $E$ be a holomorphic
vector bundle on $X$ of rank $\rank (E)\geq 1$.

Let $(L,h^{L})$ be a Hermitian line bundle on $X$ such that
$c_{1}(L,h^{L})\geq 0$ on an open neighborhood of $bM$ in
$\overline{M}$.

(I) Assume that there exists an open neighborhood $U$ of $bM$ and a plurisubharmonic
function $\phi \in \mathscr{C}^{\infty}(U,\mathbb{R})$ such that
$c_{1}(L,h^{L})=\frac{\sqrt{-1}}{2\pi}\partial \overline{\partial }
\phi $ on $U$. Then we have%
\begin{equation}
\begin{split} & \liminf _{k\rightarrow +\infty} k^{-n}\dim H^{0}(M,L^{k}
\otimes E)
\\
&\quad \geq \frac{\rank (E)}{n!}\int _{M(\,\leq 1,h^{L})} c_{1}(L,h^{L})^{n}+
\frac{\rank (E)}{n!}\int _{bM}c_{1}(L,h^{L})^{n-1}\wedge
\frac{\sqrt{-1}}{2\pi}\partial \phi ,
\end{split}
\label{eq:4.61part1}
\end{equation}
where $bM$ has the induced orientation as the boundary of $M$.

(II) Assume $c_{1}(L,h^{L})\geq 0$ on the whole of $\overline{M}$. If we can take
$\phi =-\varrho $ in the above situation, that is,
$c_{1}(L,h^{L})=-\frac{\sqrt{-1}}{2\pi}\partial \overline{\partial }
\varrho $ on $U$, then
\begin{align}
&\lim _{k\rightarrow +\infty} k^{-n}\dim H^{0}(M,L^{k}\otimes E)\notag\\
&\quad =
\frac{\rank (E)}{n!}\int _{M} c_{1}(L,h^{L})^{n}-\frac{\rank (E)}{n!}
\int _{bM}c_{1}(L,h^{L})^{n-1}\wedge \frac{\sqrt{-1}}{2\pi}\partial
\varrho .
\label{eq:4.62part2}
\end{align}
\end{thm}
When $E=\mathbb{C}$, {\eqref{eq:4.62part2}} was proven by Berman
\cite[Corollary 6.7]{RB:05} and \cite{RB:10}. In his results, it is only
assumed that
$c_{1}(L,h^{L})=-\frac{\sqrt{-1}}{2\pi}\mathscr{L}_{\varrho}$ on
$bM$. A version of {\eqref{eq:4.61part1}} was derived in
\cite[(2.1)]{M:2016-q}, under the assumption that $(L,h_{L})$ is positive
on $\overline{M}$ and $\mathscr{L}_{\varrho}$ has at least two negative
eigenvalues on $bM$.

By {Theorem~\ref{thm:4.13new}} and \cite[Theorem 5.2]{M:96} (cf.\ also
\cite{HLM22,M:2016-q}), we obtain the following result for the compactification
of strictly pseudoconcave domains.

\begin{cor}%
\label{cor:1.7new}
Let $X$ be a connected complex manifold of dimension $n\geq 3$, and let
$M$ be a connected, strictly pseudoconcave, relatively compact domain in
$X$ with a smooth boundary $bM$. Let $(L,h^{L})$ be a Hermitian line bundle
on $X$ such that $c_{1}(L,h^{L})\geq 0$ holds in an open neighborhood of
$bM$ in $\overline{M}$. We assume that there exists an open neighborhood
$U$ of $bM$ and a plurisubharmonic function
$\phi \in \mathscr{C}^{\infty}(U,\mathbb{R})$ such that
$c_{1}(L,h^{L})=\frac{\sqrt{-1}}{2\pi}\partial \overline{\partial }
\phi $ on $U$, and
\begin{equation}
\int _{M(\,\leq 1 , h^{L})} c_{1}(L,h^{L})^{n} > -\int _{bM}c_{1}(L,h^{L})^{n-1}
\wedge \frac{\sqrt{-1}}{2\pi}\partial \phi .
\label{eq:4.62partI}
\end{equation}
Then $L$ is big, and $M$ can be compactified to a Moishezon manifold; that
is, there exists a compact Moishezon manifold $\widetilde{M}$ and an open
set $M_{0}\subset \widetilde{M}$ such that $M$ is biholomorphic to
$M_{0}$.
\end{cor}

We illustrate the applications of the above corollary with the following
example.
\begin{example}[Moishezon $1$-concave domains in a compact irreducible complex space with isolated singularities]
\label{exmp1.10}
Let $X$ be a compact irreducible complex space of dimension
$n\geq 3$ with at most isolated singularities, denoted by
$X_{\mathrm{sing}}$ (see \cite[\S 3.4]{MM}). Set
$X_{\mathrm{reg}}:=X\setminus X_{\mathrm{sing}}$ and let
$D\subset X$ be a finite set that contains $X_{\mathrm{sing}}$.

Let $L\to X$ be a holomorphic line bundle equipped with a smooth Hermitian
metric $h^{L}$ on $X\setminus D$ such that, in a local chart near each
point $x\in D$, $h^{L}$ as well as all its derivatives (actually, it is
enough to consider all the derivatives up to order $2$) are uniformly
bounded. This is the case when $h^{L}$ is smooth on the whole of $X$ in
the sense that $h^{L}$ is smooth on $X_{\mathrm{reg}}$ and for each
$x_{\alpha}\in X_{\mathrm{sing}}$, there exists an open chart
$U_{\alpha}$ of $X$ and a holomorphic embedding
$\iota _{\alpha}: U_{\alpha}\hookrightarrow \mathbb{C}^{N_{\alpha}}$ such
that $x_{\alpha}\in U_{\alpha}$,
$\iota _{\alpha}(x_{\alpha})=0\in \mathbb{C}^{N_{\alpha}}$, and
$(L|_{U_{\alpha}},h^{L}|_{U_{\alpha}})=\iota_\alpha^*(\mathbb {C}^{N_\alpha}\times\mathbb {C},h_\alpha)$ is the pullback of the trivial line
bundle on $\mathbb{C}^{N_{\alpha}}$ together with a smooth Hermitian metric.

We assume $c_{1}(L,h_{L})\geq 0$ on an open deleted neighborhood of
$D$. For each point $x\in D$, we fix a small open domain $U_{x}$ of
$x$ such that
\begin{itemize}
\item each $U_{x}$ is a relatively compact, strictly pseudoconvex domain
with a smooth boundary;
\item on an open neighborhood $V_{x}$ of $\overline{U}_{x}$, the line bundle
$L$ can be holomorphically trivialized so that, in this trivialization, $h^{L}$ has a local weight
$\phi _{x}\in \mathscr{C}^{\infty}(V_{x}\setminus \{x\},\mathbb{R})$. We
also take $U_{x}$ and $V_{x}$ such that
$c_{1}(L,h_{L})|_{V_{x}\setminus \{x\}}\geq 0$. By our assumption, every
derivative of $\phi _{x}$ in the deleted neighborhood
$V_{x}\setminus \{x\}$ is uniformly bounded.
\item the open subsets $V_{x}$, $x\in D$, are disjoint from each other.
\end{itemize}
Observe that one can always choose a family
$\{x\in U_{x}\Subset V_{x}\Subset X\}_{x\in D}$ satisfying the above conditions.

Now we set $M:=X\setminus \cup _{x\in D} \overline{U}_{x}$. Then
$bM$ is the disjoint union of $bU_{x}$, but with the opposite orientation.
Then $M$ becomes a strictly pseudoconcave domain with a smooth boundary.
Our assumption on $(L,h^{L})$ implies that, for $x\in D$, the smooth form
$c_{1}(L,h^{L})^{n-1}\wedge \frac{\sqrt{-1}}{2\pi}\partial \phi _{x}$ is
uniformly bounded on a small deleted neighbourhood of $x$, then by
Stokes' theorem, we have
\begin{equation}
\int _{U_{x}\setminus \{x\}} c_{1}(L,h^{L})^{n} =-\int _{bU_{x}} c_{1}(L,h^{L})^{n-1}
\wedge \frac{\sqrt{-1}}{2\pi}\partial \phi _{x}\,.
\label{eq1.14}
\end{equation}
As a consequence, we have
\begin{equation}
\int _{X_{\mathrm{reg}}(\,\leq 1,h^{L})} c_{1}(L,h^{L})^{n}=\int _{M(
\,\leq 1, h^{L})}c_{1}(L,h^{L})^{n} + \int _{bM} c_{1}(L,h^{L})^{n-1}
\wedge \frac{\sqrt{-1}}{2\pi}\partial \phi ,
\label{eq:1.10new2025}
\end{equation}
where $\phi =\{\phi _{x}\}_{x\in D}$. By \cite[Theorem 3.4.10]{MM}, if
the integral in {\eqref{eq:1.10new2025}} is strictly positive, then
$X$ is Moishezon, which implies that $M$ is Moishezon as a consequence
of {Corollary~\ref{cor:1.7new}}.

The subsequent stability result for Moishezon manifolds is established
by combining the idea in \cite{M:2016-q} with {Corollary~\ref{cor:1.7new}} and {\eqref{eq:1.10new2025}}. The proof is provided at the
end of Section~\ref{ss:4.4semipositive}.
\end{example}
\begin{thm}%
\label{thm:stability}
Under the same setup as in the above example, assume in addition that
$c_{1}(L,h^{L})>0$ on $X\setminus D$. Then, for any sufficiently small
deformation $J'$ of the complex structure $J$ on $M$ and any sufficiently
small deformation of the complex structure on $L$ that makes the original smooth
line bundle $L|_{M} \to M$ holomorphic with respect to the new complex
structure, the manifold $M$ is Moishezon 1-concave and can be compactified
to a Moishezon manifold.
\end{thm}

\subsection{Kohn-Rossi extension theorem for
\texorpdfstring{$q$}{q}-concave domains}
\label{sec1.4}

Let $\mu $ be the projection from the space of smooth forms with values
in holomorphic bundles on the boundary to its subspace consisting of complex
tangential forms with values in bundles, and let
$\overline \partial _{b}$ be the Kohn-Rossi operator on the space of 
smooth forms with values in holomorphic vector bundles on the boundary;
see Subsection \ref{sec-dbarb}.

Under the same conditions on $M\Subset X$ as {Theorem~\ref{thm:4.10important}}, if $L$ is semi-positive throughout $M$, then {\eqref{eq:4.58Jan25-1}} indicates that
\begin{equation*}
\dim \mathscr{H}^{0,n-\ell -1}({M},L^{-k}\otimes E^{*})=o(k^{n}),
\quad \text{ as $k\rightarrow \infty $ when $q \leq \ell $}\,;
\end{equation*}
see {Corollary~\ref{cor_qconcave}}. This estimate, along with the Kohn-Rossi
criterion ({Theorem~\ref{extensioncr}}), allows for the substitution of cohomology
vanishing theorems typically used in the extension results for holomorphic
forms, leading to the following theorem.

\begin{thm}%
\label{thmextension1}
Let $X$ be a connected, compact complex manifold of dimension
$n\geq 2$. Let $(E,h^{E})$ and $(L,h^{L})$ be holomorphic Hermitian vector
bundles over $X$ and $\rank (L)=1$. Assume that $L$ is semi-positive on
$X$ and positive at least at one point. Let $1\leq q\leq n-1$ and assume
that $M$ is a relatively compact domain in $X$ with a smooth boundary
$bM$ on which the Levi form restricted to the analytic tangent space has
at least $n-q$ negative eigenvalues. Then, for $q\leq \ell \leq n-1$, there
exists a non-zero holomorphic section $s\in H^{0}(X,L^{k_{0}})$ for some
$k_{0}\in \mathbb{N}$, such that for every $\overline \partial _{b}$-closed
form $\sigma \in \Omega ^{0,\ell}(bM, E)$, there exists a
$\overline \partial $-closed extension $S$ of
$s\otimes \sigma \in \Omega ^{0,\ell}(bM, L^{k_{0}}\otimes E)$, i.e.,
\begin{equation}
\nonumber
S\in \Omega ^{0,\ell}(\overline{M}, L^{k_{0}}\otimes E)\,, \quad
\text{$\overline \partial S=0$ on $M$}\,, \quad
\text{$\mu (S|_{bM})=\mu (s\otimes \sigma )$ on $bM$.}
\end{equation}
\end{thm}

Note that the hypothesis on the line bundle $L$ implies that $X$ is a Moishezon
manifold by Siu's criterion \cite[Theorem 2.2.7]{MM}. The hypotheses of
the theorem can be weakened as follows. It is enough to assume that the
holomorphic line bundle $L$ is big and semipositive in a neighborhood of
the domain $M$. As a corollary, we obtain that under these conditions,
if the Levi form on $bM$ has one negative eigenvalue everywhere, then every
smooth $(0, n-1)$-form on $bM$ can be extended meromorphically over
$M$; see {Theorem~\ref{thmextension2}}. We will apply vanishing theorems
to prove the $\overline \partial _{b}$-extension results in {Theorems~\ref{extpos} and \ref{extneg}}. A preliminary version of {Theorem~\ref{thmextension1}} can be found in \cite{M97}.

For $k\in \mathbb{N}$, the space
$H^{0}(\mathbb{CP}^{n},\mathscr{O}(k))$ is canonically isomorphic to the
space of homogeneous polynomials of degree $k$. By choosing
$X=\mathbb{CP}^{n}$, $L=\mathscr{O}(1)$, $E$ trivial,
$\alpha =\sigma $, and $A=S\otimes s^{-1}$ in {Theorem~\ref{thmextension1}}, we obtain the following interesting case.

\begin{cor}
\label{cor1.13}
Let $M\subset \mathbb{CP}^{n}$ be an open domain with a smooth boundary
$bM$. Let $1\leq q\leq n-1$, and assume that the Levi form of $bM$ restricted to the
analytic tangent space $T^{(1,0)}bM$ has at least $n-q$ negative eigenvalues.
Then, for each $q\leq j\leq n-1$, there exists a homogeneous polynomial
$f_{j}\in \mathbb{C}[z_{0},z_{1},\cdots ,z_{n}]$ of degree
$k_{j}\in \mathbb{N}$ such that for every $\overline \partial _{b}$-closed
form $\alpha \in \Omega ^{0,j}(bM)$, there exists a
$\overline \partial $-closed extension $A$ of $\alpha $ on $M$ outside
the zero set $Z(f_{j}):=\{ [z]\in \mathbb{CP}^{n}: f_{j}(z)=0 \}$, that
is,
\begin{equation}
\nonumber
A\in \Omega ^{0,j}\left (\overline{M}\setminus Z(f_{j})\right )\,,
\quad \text{$\overline \partial A=0$ on $M\setminus Z(f_{j})$}\,,
\quad
\text{$\mu \left (A|_{bM\setminus Z(f_{j})}\right )= \mu \left (
\alpha |_{bM\setminus Z(f_{j})}\right )$}.
\end{equation}
\end{cor}

\subsection{Organization of the paper}
\label{sec1.5}

This paper is structured as follows. Section~\ref{preconvex} introduces
the notation and reviews known results on cohomology groups for complex
manifolds satisfying certain convexity or concavity conditions. In Section~\ref{sec-NG}, we derive our version of the Bochner-Kodaira-Nakano formula
and establish the Nakano-Griffiths inequality including boundary contributions.
Section~\ref{sec_cvhmi} is devoted to a systematic investigation of Morse
inequalities for $q$-concave manifolds and domains, based on our Nakano-Griffiths
inequality. In Section~\ref{proofddbar1}, we prove the Kohn-Rossi extension
theorem for $q$-concave domains; see {Theorem~\ref{thmextension1}}. Finally,
in Section~\ref{sec-vanish}, we rederive the Andreotti-Tomassini vanishing
theorem and apply it to the $\overline \partial _{b}$ extension problem;
in addition, we obtain a vanishing theorem for domains with Levi-flat boundaries.

\section{Preliminaries}
\label{preconvex}

Let $(X, \Theta )$ be a connected Hermitian manifold of dimension
$n$, and let $(F, h^{F})$ be a holomorphic Hermitian vector bundle over
$X$. Let $p,q\in \mathbb{N}$. We denote by $\Omega ^{p,q}(X,F)$ the space
of smooth $(p,q)$-forms with values in $F$, and by
$\langle \cdot ,\cdot \rangle _{h^{F},\Theta}$ the pointwise Hermitian
metric on $\Omega ^{p,q}(X,F)$ induced by $\Theta $ and $h^{F}$. Let
$M\Subset X$ be a relatively compact domain with a smooth boundary
$bM$. We denote by $\overline{M}=M\bigcup bM$ the closure of $M$ in
$X$. Throughout this section, we use these notations without further notice.

\subsection{Complex convexity}
\label{ss:2.1}

The theory of $q$-convexity/concavity was introduced by Andreotti and Grauert
in \cite{AG:62} and is one of the basic tools in the study of the geometry
of noncompact complex spaces. We recall here some important facts of the
theory.

\begin{defn}%
\label{defmfd}
Let $M$ be a connected complex manifold of dimension $n$. Let
$1\leqslant q\leqslant n$.

(i) $M$ is called \emph{$q$-convex} if there exists a smooth function
$\varphi :M\longrightarrow [a,b)$, where $a\in \mathbb{R}$,
$b\in \mathbb{R}\cup \{+\infty \}$, such that the \emph{sublevel set}
$M_{c}=\{x\in M: \varphi (x)<c\}$ is relatively compact in $M$ for all
$c\in [a,b)$ and $\sqrt{-1}\partial \overline \partial \varphi $ has at
least $n-q+1$ positive eigenvalues outside a compact set $K$ (called the
exceptional set).

(ii) $M$ is called \emph{$q$-concave} if there exists a smooth function
$\varphi :M\longrightarrow (a,b]$, where
$a\in \mathbb{R}\cup \{-\infty \}$ and $b\in \mathbb{R}$, such that the
\emph{superlevel set} $M_{c}=\{x\in M: \varphi (x)>c\}$ is relatively compact
in $M$ for all $c\in (a,b]$ and
$\sqrt{-1}\partial \overline \partial \varphi $ has at least $n-q+1$ positive
eigenvalues outside a compact set $K\subset M$, called
\emph{an exceptional compact set}. In the sequel, we also use
$\varrho :=-\varphi $ as the sublevel exhaustion function and then
$M_{c}=\{x\in M: \varrho (x) < c\} \Subset M$ for all
$c \in [-b, -a)$.

In all these cases, we call $\varphi $ or $\varrho $ an
\emph{exhaustion function}.
\end{defn}

\begin{example}
We provide here some examples of $q$-concave manifolds (cf.\ \cite{M:96}).

(1) Let $X$ be a compact K\"ahler space of pure dimension $n$ and
$Y$ an analytic subset of pure dimension $q$ containing the singular locus
of $X$. Then $X\setminus Y$ is a $(q+1)$-concave manifold. In particular,
the regular locus of a projective algebraic variety with isolated singularities
is $1$-concave.

(2) Let $X$ be a compact complex manifold and $Y$ an analytic subset of
pure dimension $q$. Then $X\setminus Y$ is a $(q+1)$-concave manifold.

(3) If $\pi :X\rightarrow Y$ is a proper holomorphic map between a complex
manifold $X$ and a $q$-concave manifold $Y$ such that the dimension of
its fibers does not exceed $r$, then $X$ is $(q+r)$-concave.
\end{example}

\begin{example}
\label{em-hyperconcave}
A hyperconcave (or hyper $1$-concave) manifold is a $1$-concave manifold
with $a=-\infty $; cf.\ \cite[Definition 3.4.1]{MM}. A geometric example
due to Siu--Yau is a complete K\"ahler manifold of finite volume and bounded
negative sectional curvature; see more examples in
\cite[Example 3.4.2]{MM}.

We also need the following terminology concerning the Levi convexity of
domains $M$ with smooth boundaries in a complex manifold $X$ of dimension
$n$, cf.\ \cite{FK:72}. We assume that there exists a real smooth function
$\varrho $ on $X$ such that $M=\{x\in X :\varrho (x)<0\}$,
$bM=\{x\in X: \varrho (x)=0\}$, and $d \varrho (x)\neq 0$ for any
$x\in bM$. In addition, we assume that $|d\varrho |=1$ on $bM$. Then, we
call $\varrho $ a defining function of $M$. Let
$TX\otimes _{\mathbb{R}}\mathbb{C}=T^{(1,0)}X\oplus T^{(0,1)}X$ be the
splitting of the complexified tangential bundle. The analytic tangent space
to $bM$ at $x\in bM$ is given by
\begin{equation*}
T_{x}^{(1,0)}bM:=\{ v\in T_{x}^{(1,0)}X : \partial \varrho (v)=0 \}.
\end{equation*}

The Levi form of $\varrho $ is the $2$-form
$ \mathscr{L}_{\varrho}\in \mathscr{C}^{\infty }(bM, T^{(1,0)*}bM
\otimes T^{(0,1)*}bM )$ given by
\begin{equation}
\mathscr{L}_{\varrho}(U,\overline{V}):=(\partial \overline \partial
\varrho )(U,\overline{V}),\quad
\text{for $U, V\in T_{x}^{(1,0)}bM$, $x\in bM$}.
\label{eq:2.2Levi}
\end{equation}
The number of positive and negative eigenvalues of the Levi form is independent
of the choice of the defining function. We will denote the Levi form by
$\mathscr{L}_{bM}$ when the defining function $\varrho $ is not emphasized.
\end{example}

\begin{defn}%
\label{def:Zq}
We say that $M$ satisfies the condition $Z(q)$ if the Levi form
$\mathscr{L}_{\varrho}$ has at least $n-q$ positive or at least
$q+1$ negative eigenvalues at each point of $bM$.
\end{defn}

Immediately, we have the following result (see
\cite[Theorem 6.7]{HL88}).

\begin{prop}%
\label{prop:2.3dec24}
Let $M$ be a $q$-concave manifold with superlevel exhaustion function
$\varphi $ and exceptional compact set $K$. Then the following claims hold
for any $c<\inf _{K} \varphi $:

(a) The manifold $M_{c}=\{x\in M\;:\; \varphi >c\}$ is $q$-concave.

(b) If in addition $c$ is a regular value of $\varphi $, then
$M_{c}$ satisfies the $Z(j)$ condition for $j\leq n-q-1$.
\end{prop}

Let $M\Subset X$ be a smooth domain such that its Levi form
$\mathscr{L}_{bM}$ has at least $n-q$ negative eigenvalues on $bM$ with
$1\leq q\leq n-1$. Let $\varrho $ be its defining function. For
$\varepsilon \in \mathbb{R}$ with $|\varepsilon |$ sufficiently small, the Levi form of
$M_{\varepsilon}:=\{x\in X: \varrho (x)<\varepsilon \}$ still has at least
$n-q$ negative eigenvalues on $bM_{\varepsilon}$. Following
\cite[Theorem 6.7]{HL88}, if we set
$\varrho _{1}=-\exp (-C\varrho )+1$ for a sufficiently large constant
$C>0$, then we can rewrite $M=\{x\in X\;:\; \varrho _{1}<0\}$, and
$\sqrt{-1}\partial \overline{\partial }\varrho _{1}$ has at least
$n-q+1$ negative eigenvalues in a small neighborhood of $bM$.  So
$M$ is a $q$-concave manifold.

\begin{defn}%
\label{def:2.6domain}
A smooth domain $M\Subset X$ is called (Levi) $q$-concave ($1\leq q
\leq n-1$) if the Levi form $\mathscr{L}_{bM}$ has at least $n-q$ negative
eigenvalues on $bM$.
\end{defn}

The following important cohomology finiteness theorem is due to Andreotti-Grauert
\cite{AG:62}; see also \cite{AV65,FK:72,Hor:65,Oh:82}.
\begin{thm}[Andreotti-Grauert]
\label{thm2.7}
Let $M$ be a $q$-concave manifold of dimension $n\geq 2$, and let
$F$ be a holomorphic vector bundle on $M$. Then the sheaf cohomology groups
of $M$ with values in $F$ satisfy
\begin{equation}
\dim H^{j}(M,F)<\infty ,\quad j\leq n-q-1.
\label{eq2.2}
\end{equation}
\end{thm}

\subsection{\texorpdfstring{$\overline \partial $}{overline{∂}}-Neumann problem
and Hodge theory}
\label{ss:cohomology}

Let $p,q\in \mathbb{N}$ and $0\leq p,q\leq n$. Let
$\Omega ^{p,q}(M, F)$ be the space of smooth $(p,q)$-forms with values
in $F$. Let
$\Omega ^{p,q}(\overline{M}, F)\subset \Omega ^{p,q}(M, F)$ be the subspace
of those forms that are smooth up to and including the boundary $bM$. The
$L^{2}$-scalar product on $\Omega ^{p,q}(M, F)$ is given by
\begin{eqnarray}
\langle s_{1},s_{2} \rangle _{L^{2}(M,F)}:=\int _{M} \langle s_{1}(x),
s_{2}(x) \rangle _{h^{F},\,\Theta}\, dv_{M}(x)
\label{eq2.3}
\end{eqnarray}
where $dv_{M}=\Theta ^{n}/ n!$ is the volume form of $M$. We will write
simply $\langle s_{1},s_{2} \rangle _{L^{2}}$ or
$\langle s_{1},s_{2} \rangle _{L^{2}(M)}$ for
$\langle s_{1},s_{2} \rangle _{L^{2}(M,F)}$ if there is no confusion, and
we denote by $\|\cdot \|_{L^{2}}$ the corresponding $L^{2}$-norm and by
$L^{2}_{p,q}({M}, F)$ the $L^{2}$ completion of
$\Omega ^{p,q}(\overline{M},F)$. Let
\begin{equation*}
\overline \partial ^{F}: \Dom (\overline \partial ^{F})\cap L^{2}_{p,q}(M,F)
\rightarrow L^{2}_{p,q+1}({M}, F)
\end{equation*}
be the maximal extension of the Cauchy-Riemann operator. Sometimes we will
use $\overline \partial $ instead of $\overline \partial ^{F}$ to simplify
the notation. Let $\overline \partial ^{F*}$ be the Hilbert space adjoint
of $\overline \partial ^{F}$.

Let $e_{\mathfrak{n}}$ be the inward pointing unit normal vector at
$bM$. We decompose
$e_{\mathfrak{n}}=e^{(1,0)}_{\mathfrak{n}}+e^{(0,1)}_{\mathfrak{n}}
\in T^{(1,0)}X\oplus T^{(0,1)}X$. Set
\begin{equation}
\label{B_pq}
B^{p,q}(M,F)=\{s\in \Omega ^{p,q}(\overline{M},F)\;:\; \iota _{e^{(0,1)}_{
\mathfrak{n}}} s=0 \text{ on } bM\}.
\end{equation}
Then by \cite[Proposition 1.4.19]{MM}, we have
\begin{equation}
B^{p,q}(M,F)=\Dom (\overline \partial ^{F\ast})\cap \Omega ^{p,q}(
\overline{M},F).
\label{eq2.5}
\end{equation}

Furthermore, we define the Gaffney extension of the Kodaira Laplacian by
\begin{equation}
\begin{split} \Dom (\square ^{F})&=\{ s\in \Dom (\overline \partial ^{F})
\cap \Dom (\overline \partial ^{F*}) : \overline \partial ^{F} s \in
\Dom (\overline \partial ^{F*}),~ \overline \partial ^{F*}s \in \Dom (
\overline \partial ^{F}) \},
\\
\square ^{F}&:=\overline \partial ^{F} \overline \partial ^{F*}+
\overline \partial ^{F*}\overline \partial ^{F}.
\end{split}
\label{eq2.6}
\end{equation}
We denote the space of harmonic $(p,q)$-forms with values in $F$ by
\begin{equation}
\mathscr{H}^{p,q}({M}, F):=\Ker (\square ^{F})\cap L^{2}_{p,q}(M,F)=
\{ s\in \Dom (\square ^{F})\cap L^2_{p,q}(M,F): \square ^{F} s=0 \}
\label{eq2.7}
\end{equation}
and the orthogonal projection by
\begin{eqnarray}
H: L^{2}_{p,q}({M}, F)\rightarrow \mathscr{H}^{p,q}({M}, F).
\label{eq2.8}
\end{eqnarray}

We recall the solution of the $\overline \partial $-Neumann problem.
\begin{thm}[{\cite[(3.1.11), (3.1.14)]{FK:72}};{\cite[Theorem 5.11]{KR:65}}]
\label{neumannop}
If $M$ satisfies the $Z(q)$ condition, then there exists a compact operator
(the Neumann operator)
\begin{equation*}
\mathcal{N}: L^{2}_{p,q}({M}, F)\rightarrow L^{2}_{p,q}({M}, F)
\end{equation*}
such that
\begin{itemize}
\item[(a)]
$\mathcal{N}L^{2}_{p,q}({M}, F) \subset \Dom (\square ^{F})$ and
$ L^{2}_{p,q}({M}, F)=\square ^{F} \mathcal{N}L^{2}_{p,q}({M}, F)+
\mathscr{H}^{p,q}({M}, F) $;
\item[(b)] $\mathcal{N}$ commutes with
$\square ^{F},~\overline \partial ^{F},~\overline \partial ^{F*},~H$;
\item[(c)]
$\mathcal{N}(\Omega ^{p,q}(\overline{M}, F))\subset \Omega ^{p,q}(
\overline{M}, F)$ and
$H(\Omega ^{p,q}(\overline{M}, F))\subset \Omega ^{p,q}(\overline{M}, F)$;
\item[(d)]
$ \mathscr{H}^{p,q} ({M}, F)\subset \Omega ^{p,q}(\overline M, F)$ is finite
dimensional.
\end{itemize}
\end{thm}

The reduced $L^{2}$-Dolbeault cohomology is defined by
\begin{equation}%
\overline{H}^{p,q}_{(2)}(M,F):=
\dfrac{\Ker (\overline \partial ^{F})\cap L^{2}_{p,q}(M,F) }{[ {\mathrm{Im}}( \overline \partial ^{F}) \cap L^{2}_{p,q}(M,F)]},
\label{eq2.9}
\end{equation}
where $[V]$ denotes the closure of the space $V$. By the weak Hodge decomposition
(see \cite[(3.1.21) and (3.1.22)]{MM}),
\begin{equation}%
\overline{H}^{p,q}_{(2)}(M,F)\cong \mathscr{H}^{p,q}(M,F).
\label{eq2.10}
\end{equation}
The (non-reduced) $L^{2}$-Dolbeault cohomology is defined by
\begin{equation}
H^{p,q}_{(2)}(M,F):=
\dfrac{\Ker (\overline \partial ^{F})
\cap L^{2}_{p,q}(M,F) }{ {\mathrm{Im}}( \overline \partial ^{F}) \cap L^{2}_{p,q}(M,F)}.
\label{eq2.11}
\end{equation}
We say the fundamental estimate holds in bidegree $(p,q)$ for forms with
values in $F$ if there exists a compact subset $K\subset M$ and
$C>0$ such that
\begin{equation}
\|s\|^{2}_{L^{2}(M)}\leq C(\|\overline \partial ^{F} s\|^{2}_{L^{2}(M)}+
\|\overline \partial ^{F*}s\|^{2}_{L^{2}(M)}+\int _{K}|s|^{2}dv_{M}),
\label{eq:2.12FE}
\end{equation}
for
$s\in \Dom (\overline \partial ^{F})\cap \Dom (\overline \partial ^{F,*})
\cap L^{2}_{p,q}(M,F)$. If the fundamental estimate holds, the reduced
and non-reduced $L^{2}$-Dolbeault cohomology coincide; see
\cite[Theorem 3.1.8]{MM}.

The sheaf cohomology $H^{q}(M,F)$ for holomorphic sections of $F$ is isomorphic
to the Dolbeault cohomology $H^{0,q}(M,F)$. Besides, the Dolbeault cohomology
up to the boundary is defined by
\begin{equation}
H^{p,q}(\overline M,F)=
\dfrac{\Ker (\overline \partial ^{F})\cap \Omega ^{p,q}(\overline M,F)}{ {\mathrm{Im}}( \overline \partial ^{F}) \cap \Omega ^{p,q}(\overline M,F)}.
\label{eq2.13}
\end{equation}

Regarding the relationships among the aforementioned cohomology groups,
we present the following fundamental results.

\begin{thm} [{\cite[(3.1.14), (4.3.1)]{FK:72}},
{\cite[Theorem 3.4.9]{Hor:65}}]
\label{zqzq1}
If $M$ satisfies $Z(q)$, then
\begin{equation}
H^{p,q}(\overline M,F)\simeq H^{p,q}_{(2)}(M,F)\simeq \mathscr{H}^{p,q}(M,F).
\label{eq:2.14iso25}
\end{equation}
If $M$ satisfies $Z(q)$ and $Z(q+1)$, then
\begin{equation*}
H^{p,q}(M, F)\simeq \mathscr{H}^{p,q}({M}, F).
\end{equation*}
The spaces $H^{p,q}(\overline M,F)$ and $H^{p,q}(M, F)$ are independent
of the metrics used on the domain $M$ and $F$, whereas the other spaces
mentioned rely on metric data.
\end{thm}

\begin{thm}[{\cite[Proposition 18]{AG:62}}]%
\label{thm:2.9AG}
If $X$ is a connected $q$-concave manifold of dimension $n\geq 2$ equipped
with an exhaustion function $\varrho $ as in {Definition~\ref{defmfd}} and
$1\leq q\leq n-1$, then for a regular value $c\in \mathbb{R}$ of
$\varrho $ such that the sublevel set $X_{c}$ contains the exceptional
set $K$ of $\varrho $, we have
\begin{equation}
H^{j}(X,F)\simeq H^{j}(X_{c},F)\simeq H^{0,j}(X_{c},F), \quad
\text{for $j\leq n-q-1$}.
\label{eq2.15}
\end{equation}
\end{thm}

\begin{thm}[{\cite[Theorem 3.4.9 and the subsequent remark]{Hor:65}}]%
\label{thm:2.12H65}
Under the same hypothesis as {Theorem~\ref{thm:2.9AG}}, we have for
$0\leq j\leq n-q-1$,
\begin{equation}
\dim H^{j}(X,F)\leq \dim H^{0,j}_{(2)}(X_{c},F) =\dim H^{0,j}(
\overline X_{c},F).
\label{eq2.16}
\end{equation}
\end{thm}

\subsection{Fundamental estimates and abstract holomorphic Morse inequalities}
\label{sub:2.3abstract}

We recall the results presented in \cite[Section 3.2]{MM} on the abstract
holomorphic Morse inequalities for $L^{2}$-cohomology. The main idea is
to show that the spectral spaces of the Laplacian, associated with the
small eigenvalues, inject into the spectral spaces of the Laplacian with
Dirichlet boundary conditions on a smooth, relatively compact domain. The
asymptotics of the latter operator can be calculated explicitly; for example,
see \cite[Theorem 3.2.9]{MM}. We also refer to \cite{Bouche1989,Demailly1985,M:96} for these results.

Let $(X,\Theta )$ be a Hermitian manifold of dimension $n$. Let $L$,
$E$ be holomorphic vector bundles on $X$ equipped with smooth Hermitian
metrics $h^{L}$, $h^{E}$. Throughout, we assume $\mathrm{rank}(L)=1$. Let
$\overline \partial ^{E}_{k}: \Omega _{0}^{0,\bullet}(X,L^{k}\otimes E)
\rightarrow \Omega _{0}^{0,\bullet +1}(X,L^{k}\otimes E)$ be the associated
$\overline{\partial }$-operator, and we use the same notation for its maximal
extension with respect to the $L^{2}$-inner product.

We say that \textit{the fundamental estimate} holds in bidegree $(0,q)$ for
forms with values in $L^{k}\otimes E$, if there exists a compact subset
$K\subset X$ and $C_{0}>0$ such that for sufficiently large $k$, we have
for all
$s\in \Dom (\overline \partial ^{E}_{k})\cap \Dom (
\overline \partial ^{E,*}_{k})\cap L^{2}_{(0,q)}(X,L^{k}\otimes E)$,
\begin{equation}
\|s\|^{2}_{L^{2}(X)}\leq \frac{C_{0}}{k}\left ( \|\overline \partial ^{E}_{k}
s\|^{2}_{L^{2}(X)}+\|\overline \partial ^{E,*}_{k} s\|^{2}_{L^{2}(X)}
\right )+C_{0} \int _{K} |s|^{2} dv_{X}.
\label{eq:2.14FEoptimal}
\end{equation}
Note that {\eqref{eq:2.14FEoptimal}} implies the fundamental estimate {\eqref{eq:2.12FE}} for sufficiently large $k$.

\begin{thm}[{\cite[Theorem 3.2.13]{MM}}]%
\label{thm-l2hmi}
Assume that there exists $m\in \mathbb{N}$ such that the fundamental estimate {\eqref{eq:2.14FEoptimal}} holds for $(0,j)$-forms with values in
$L^{k}\otimes E$ and any $j\leq m$. Let $U\Subset X$ with a smooth boundary
be such that $K\Subset U$. Then as $k\rightarrow \infty $, the following
strong and weak Morse inequalities hold for any $j\leq m$,
\begin{equation}
\sum _{\ell =0}^{j}(-1)^{j-\ell}\dim H_{(2)}^{0,\ell}(X,L^{k}\otimes E)
\leq \rank (E)\frac{k^{n}}{n!}\int _{U(\leq \,j,h^{L})}(-1)^{j}c_{1}(L,h^{L})^{n}+o(k^{n}),
\label{eq2.18}
\end{equation}
and
\begin{equation}
\dim H_{(2)}^{0,j}(X,L^{k}\otimes E)\leq \rank (E)\frac{k^{n}}{n!}
\int _{U(j,h^{L})}(-1)^{j}c_{1}(L,h^{L})^{n}+o(k^{n}).
\label{eq:2.20weak}
\end{equation}
\end{thm}
As a consequence of
\cite[Theorem 3.2.9, Proposition 3.2.11, and (3.2.63)]{MM}, we have the
following result (implicit in the proof of \cite[Theorem 3.2.13]{MM}).
\begin{cor}%
\label{cor:2.13new}
Assume that the fundamental estimate {\eqref{eq:2.14FEoptimal}} holds for
$(0,q)$-forms with values in $L^{k}\otimes E$ for some
$q\in \{0,1,\ldots , n\}$. Let $U\Subset X$ with a smooth boundary be such
that $K\Subset U$. Then as $k\rightarrow \infty $, the weak Morse inequalities {\eqref{eq:2.20weak}} hold at degree $(0,q)$.
\end{cor}

For the space of $L^{2}$-holomorphic sections
$H^{0}_{(2)}(X,L^{k}\otimes E):= H^{0,0}_{(2)}(X,L^{k}\otimes E)$, there
is also a lower bound provided the fundamental estimate holds.

\begin{thm}[{\cite[Theorem 3.2.16]{MM}}]%
\label{thm-l2hmi-2}
Assume that the fundamental estimate {\eqref{eq:2.14FEoptimal}} holds for
$(0,1)$-forms with values in $L^{k}\otimes E$. Let $U\Subset X$ with a
smooth boundary such that $K\Subset U$. Then as
$k\rightarrow \infty $, we have
\begin{equation}
\dim H_{(2)}^{0}(X,L^{k}\otimes E)\geq \rank (E)\frac{k^{n}}{n!}\int _{U(
\leq 1\,,h^{L})}c_{1}(L,h^{L})^{n}+o(k^{n}).
\label{eq2.20}
\end{equation}
\end{thm}

\subsection{The \texorpdfstring{$\overline \partial _{b}$}{dbarb} operator on the boundary}
\label{sec-dbarb}

Let $\Omega ^{p,q}(bM,F)$ be the $(p,q)$-forms with values in $F$ over
$bM$, that is, $\Omega ^{p,q}(bM,F):=\Omega ^{p,q}(X,F)|_{bM}$. We recall
some notions in \cite{KR:65}. A form
$\alpha \in \Omega ^{p,q}(bM,F)$ is called complex normal, if there exists
$\psi \in \Omega ^{p,q-1}(bM,F)$ such that
$\alpha =\psi \wedge (\overline \partial \varrho )|_{bM}$, where
$\varrho $ is a defining function of $M$. We denote by
$\mathcal{C}^{p,q}(bM,F)$ the subspace of $\Omega ^{p,q}(bM,F)$ consisting
of complex normal forms. A form $\beta \in \Omega ^{p,q}(bM,F)$ is called
complex tangential, if
$\langle \alpha (x), \beta (x) \rangle _{h^{F},\Theta}=0$ for every
$\alpha \in \mathcal{C}^{p,q}(bM,F)$ and every $x \in bM$. We denote by
$\mathcal{D}^{p,q}(bM,F)$ the subspace of $\Omega ^{p,q}(bM,F)$ consisting
of complex tangential forms. With respect to the pointwise Hermitian product
$\langle \cdot , \cdot \rangle _{h^{F},\Theta}$, we have
\begin{equation}
\label{eqdecomp}
\Omega ^{p,q}(bM,F)=\mathcal{C}^{p,q}(bM,F)\oplus \mathcal{D}^{p,q}(bM,F).
\end{equation}
Moreover, we denote the projections by
$\mu :\Omega ^{p,q}(bM,F)\rightarrow \mathcal{D}^{p,q}(bM,F)$,
$p,q\in \mathbb{N}$.

A form $\Sigma \in \Omega ^{p,q}(\overline{M},E)$ is a
$\overline \partial $-closed extension of
$\sigma \in \Omega ^{p,q}(bM,E)$ if $\overline \partial \Sigma =0$ on
$M$ and $\mu (\Sigma |_{bM})=\mu (\sigma )$, i.e., $\Sigma $ is
$\overline \partial $-closed and $\Sigma =\sigma $ in the complex tangential
direction on $bM$. We define a map
\begin{equation}
\label{eqdbarb}
\overline \partial _{b}: \Omega ^{p,q}(bM,E) \rightarrow \Omega ^{p,q+1}(bM,E),
\quad \overline \partial _{b} \sigma :=\mu ((\overline \partial
\sigma ')|_{bM})
\end{equation}
where $\sigma '\in \Omega ^{p,q}(\overline{M},E)$ and
$\sigma '|_{bM}=\sigma $. The definition of
$\bar\partial_b$ is independent of the choice of $\sigma'$. We say $\sigma $ is
$\overline \partial _{b}$-closed if
$\overline \partial _{b} \sigma =0$. Kohn-Rossi provided a general criterion
for the extension of $\overline \partial _{b}$-closed forms as follows.
\begin{thm} [{\cite[Proposition 5.13]{KR:65}}{\cite[Theorem 5.3.1]{FK:72}}]
\label{extensioncr}
Suppose $bM$ has property $Z(n-q-1)$. Let $F$ be a holomorphic vector
bundle over $X$. If $\sigma \in \Omega ^{p,q}(bM, F)$, there exists a
$\overline \partial $-closed extension $\Sigma $ of $\sigma $ if and only
if
\begin{equation}
\overline \partial _{b} \sigma =0\quad \text{and} \quad \int _{bM}
\theta \wedge \sigma =0\quad \,\text{ for }\,\theta \in \mathscr{H}^{n-p,n-q-1}({M},
F^{*}).
\label{eq2.23}
\end{equation}
\end{thm}

\section{Bochner-Kodaira-Nakano formula with boundary term}
\label{sec-NG}

Let $(X,J,\Theta )$ be a Hermitian manifold of dimension $n\geq 1$ where
$J$ denotes the complex structure, and define the associated $J$-invariant Riemannian
metric on $TX$ by
\begin{equation}
g^{TX}(\cdot ,\cdot ):=\Theta (\cdot , J\cdot ).
\label{eq3.1}
\end{equation}
The corresponding volume form is $dv_{X}=\Theta ^{n}/n!$. Let
$\nabla ^{TX}$ denote the Levi-Civita connection of $(TX,g^{TX})$, and
let $R^{TX}$ denote the Riemannian curvature tensor. Let
$h^{T^{(1,0)}X}$ (resp.\ $h^{T^{(0,1)}X}$) denote the Hermitian metrics
on $T^{(1,0)}X$ (resp. $T^{(0,1)}X$) induced by $g^{TX}$.

Let $F$ be a holomorphic vector bundle on $X$ equipped with a smooth Hermitian
metric $h^{F}$. In this section, we will use the notation
$\langle \cdot ,\cdot \rangle $ or
$\langle \cdot ,\cdot \rangle _{g^{TX}}$ to denote the bilinear form on
$TX$ or $TX\otimes _{\mathbb{R}}\mathbb{C}$, and we will denote
$\langle \cdot ,\cdot \rangle _{h^{F},\Theta}$ or simply
$\langle \cdot ,\cdot \rangle _{h}$ the Hermitian inner products.

Let $\nabla ^{F}$ denote the Chern connection of $(F,h^{F})$, and let
$R^{F}$ denote the corresponding Chern curvature form. In particular, let
$\nabla ^{T^{(1,0)}X}$ be the Chern connection of the holomorphic Hermitian
vector bundle $(T^{(1,0)}X,h^{T^{(1,0)}X})$. For
$s\in \mathscr{C}^{\infty}(X,T^{(0,1)}X)$,
$\overline s\in \mathscr{C}^{\infty}(X,T^{(1,0)}X)$, we have
\begin{equation*}
\nabla ^{T^{(1,0)}X}\overline s\in \Omega ^{1}(X,T^{(1,0)}X)=
\mathscr{C}^{\infty}(X,\mathbb{C}T^{*}X \otimes T^{(1,0)}X).
\end{equation*}
We define
$\nabla ^{T^{(0,1)}X}s=\overline{\nabla ^{T^{(1,0)}X}\overline s}\in
\mathscr{C}^{\infty}(X,\mathbb{C}T^{*}X \otimes T^{(0,1)}X)$. Set
\begin{equation}
\widetilde{\nabla }^{TX}:=\nabla ^{T^{(1,0)}X}\oplus \nabla ^{T^{(0,1)}X}
\label{eq3.2}
\end{equation}
Then $\widetilde{\nabla }^{TX}$ defines a connection on
$\mathbb{C}\otimes _{\mathbb{R}}TX$ which preserves the real tangent bundle
$TX$ and the Riemannian metric $g^{TX}$.

We define the covariant derivative
\begin{equation}
\widetilde{\nabla }^{TX}_{Y} (\alpha \wedge \beta )=
\widetilde{\nabla }^{TX}_{Y} (\alpha )\wedge \beta +\alpha \wedge
\widetilde{\nabla }^{TX}_{Y} (\beta )
\label{eq3.3}
\end{equation}
for differential forms $\alpha ,\beta $. For $\alpha \in T^{*}X$,
$U,V\in TX$,
\begin{equation}
U(\alpha (V))=(\widetilde{\nabla }^{TX}_{U}\alpha )V+\alpha (
\widetilde{\nabla }^{TX}_{U} V).
\label{eq3.4}
\end{equation}
When we consider the differential forms valued in $F$, we also use the
same notation $\widetilde{\nabla }^{TX}$ for the connection or covariant
derivative given by
$\widetilde{\nabla }^{TX}\otimes \Id _{F} +\Id _{\Lambda (T^{\ast }X)}
\otimes \nabla ^{F}$.

Let $T\in \wedge ^{2} T^{*}X\otimes TX$ be the torsion of the connection
$\widetilde{\nabla }^{TX}$ given by
\begin{equation}
T(U,V)=\widetilde{\nabla }^{TX}_{U} V-\widetilde{\nabla }^{TX}_{V} U-[U,V]
\label{eq:3.5torsion}
\end{equation}
for $U,V\in TX$. By \cite[(1.2.37)]{MM}, $T$ maps
$T^{(1,0)}X\otimes T^{(1,0)}X$ (resp. $T^{(0,1)}X\otimes T^{(0,1)}X$) into
$T^{(1,0)}X$ (resp. $T^{(0,1)}X$), and vanishes on
$T^{(1,0)}X\otimes T^{(0,1)}X$. We will also write
$T=T^{(1,0)}+T^{(0,1)}$. Note that $\Theta $ is K\"{a}hler if and only
if $\nabla ^{TX}$ preserves the splitting
$T^{(1,0)}X\oplus T^{(0,1)}X$ and in this case we have
$\nabla ^{TX}=\widetilde{\nabla }^{TX}$, $T=0$.

In this section, we consider a smooth, relatively compact domain $M$ in
$X$. As in Subsection \ref{ss:2.1}, we will always set
$M=\{x\in X \;:\; \varrho (x) < 0\}$, where
$\varrho \in \mathscr{C}^{\infty}(X,\mathbb{R})$ is a defining function
of $M$ such that $|d\varrho |=1$ on $bM$. The Hermitian metric
$h^{T^{(1,0)}bM}$ on $T^{(1,0)}bM$ is induced by $h^{T^{(1,0)}X}$. In this
setup, given any orthonormal basis $\{\omega _{j}\}_{j=1}^{n}$ of
$(T^{(1,0)}X, h^{T^{(1,0)}X})$ at a point $y\in bM$, we have
\begin{equation}
e^{(1,0)}_{\mathfrak{n}}=-\sum _{j=1}^{n} \bar{\omega}_{j}(\varrho )
\omega _{j},\; e^{(0,1)}_{\mathfrak{n}}=-\sum _{j=1}^{n} \omega _{j}(
\varrho )\bar{\omega}_{j}.
\label{eq3.6}
\end{equation}
Let $dv_{bM}$ be the volume form on $bM$ defined by
\begin{equation}
dv_{X}=d\varrho \wedge dv_{bM} \; \text{ on } bM.
\label{eq3.7}
\end{equation}

\subsection{Bochner-Kodaira-Nakano formula with boundary terms}
\label{ss3.1Nakano}

We recall the Bochner-Kodaira-Nakano formulas with boundary terms given
by Andreotti-Vesentini \cite[p. 113]{AV65}, Griffiths
\cite[(7.14)]{Griffiths1966}, and Ma-Marinescu
\cite[Theorem 1.4.21]{MM}.

First, we introduce the necessary notation. Set the Lefschetz operator
$\mathrm{Lef}_{\Theta}=(\Theta \wedge )\otimes \Id _{F}$ on
$\Lambda ^{\bullet ,\bullet}(T^{\ast }X)\otimes F$, and let
$\Lambda =\iota _{\Theta}$ be its adjoint with respect to the Hermitian
product $\langle \cdot ,\cdot \rangle _{h}$.

\begin{defn}
\label{defn3.1}
The Hermitian torsion operator $\mathcal{T}$ acting on
$\Lambda ^{\bullet ,\bullet}(T^{\ast }X)$ is given as follows
\begin{equation}
\mathcal{T}=\frac{1}{2}\sum _{j,k,m}\langle T(\omega _{j},\omega _{k}),
\bar{\omega}_{m}\rangle \left [2\omega ^{k}\wedge \bar{\omega}^{m}
\wedge \iota _{\omega _{j}}-2\delta _{jm}\omega ^{k}-\omega ^{j}
\wedge \omega ^{k}\wedge \iota _{\omega _{m}}\right ],
\label{eq3.8}
\end{equation}
where $\{\omega _{j}\}_{j=1}^{n}$ is an orthonormal basis of
$(T^{(1,0)}X, h^{T^{(1,0)}X})$ with the dual basis
$\{\omega ^{j}\}_{j=1}^{n}$ of $T^{(1,0)\ast}X$. We will use the same notation
$\mathcal{T}$ for the operator $\mathcal{T}\otimes \Id _{F}$ on
$\Lambda ^{\bullet ,\bullet}(T^{\ast }X)\otimes F$.
\end{defn}

We write $\nabla ^{F}=(\nabla ^{F})^{1,0}+(\nabla ^{F})^{0,1}$ as operators
acting on differential forms valued in $F$. Let
$(\nabla ^{F})^{1,0\ast}$, $(\nabla ^{F})^{0,1\ast}$,
$\mathcal{T}^{\ast}$, $\overline{\mathcal{T}}^{\ast}$ denote the (formal)
adjoints of $(\nabla ^{F})^{1,0}$,
$(\nabla ^{F})^{0,1}=\overline \partial ^{F}$, $\mathcal{T}$,
$\overline{\mathcal{T}}$ with respect to the $L^{2}$-metric or the Hermitian
metric $\langle \cdot ,\cdot \rangle _{h}$. Analogous to
$\square ^{F}$, the holomorphic Kodaira Laplacian is defined as
\begin{equation}
\overline{\square }^{F}=[(\nabla ^{F})^{1,0},(\nabla ^{F})^{1,0\ast}].
\label{eq3.9}
\end{equation}
We recall the Bochner-Kodaira-Nakano formula of Demailly
\cite[Proposition 2.1]{Dem1986}:
\begin{equation}
\square ^{F}=\overline{\square }^{F}+[\sqrt{-1}R^{F},\Lambda ]+ [(
\nabla ^{F})^{1,0},\mathcal{T}^{\ast}]-[\overline \partial ^{F},
\overline{\mathcal{T}}^{\ast}].
\label{eq:3.9Nov}
\end{equation}
Note that the curvature term $[\sqrt{-1}R^{F},\Lambda ]$ is crucial for
various applications of the above formula. By the definition of
$\Lambda $, in an orthonormal basis $\{\omega _{j}\}_{j=1}^{n}$ of
$(T^{(1,0)}X, h^{T^{(1,0)}X})$, we have
\begin{equation}
[\sqrt{-1}R^{F}, \Lambda ]=\sum _{j,k}R^{F} (\omega _{j},
\overline \omega _{k})\omega ^{j}\wedge \iota _{\omega _{k}}-\sum _{j,k}R^{F}(
\omega _{j},\overline \omega _{k})\iota _{\overline \omega _{j}}
\overline \omega ^{k}\wedge .
\label{eq3.11}
\end{equation}
Write
\begin{equation}
\mathrm{Tr}_{\Theta}[\sqrt{-1}R^{F}]=\sum _{j=1}^{n} R^{F} (\omega _{j},
\overline \omega _{j})\in \mathrm{End}(F).
\label{eq:3.11Nov}
\end{equation}
For $s\in \Omega ^{0,q}(X,F)$, $0\leq q \leq n$, we have
\begin{equation}
[\sqrt{-1}R^{F}, \Lambda ]s=\sum _{j,k}R^{F}(\omega _{j},
\overline \omega _{k})\overline \omega ^{k}\wedge \iota _{
\overline \omega _{j}}s- \mathrm{Tr}_{\Theta}[\sqrt{-1}R^{F}]s.
\label{eq3.13}
\end{equation}

Set $\widetilde{F}=F\otimes K^{\ast}_{X}$ with
$K^{\ast}_{X}:=\Lambda ^{n} (T^{(1,0)}X)=\det T^{(1,0)}X$. Since
$K_{X}\otimes K^{\ast}_{X}\simeq \mathbb{C}$, there exists a natural isometry
\begin{equation}
\label{isopsi}
\Psi : \Lambda ^{0,q}(T^{\ast }X)\otimes F \rightarrow \Lambda ^{n,q}(T^{
\ast }X)\otimes \widetilde{F}, \quad \Psi s:=(\omega ^{1}\wedge
\ldots \wedge \omega ^{n}\wedge s) \otimes (\omega _{1}\wedge \ldots
\wedge \omega _{n}).
\end{equation}
Now we consider the smooth, relatively compact domain
$M=\{x\in X: \varrho <0\}\Subset X$ with a defining function
$\varrho \in \mathscr{C}^{\infty}(X)$ satisfying $|d\varrho |=1$ on
$bM$, and the sections in $B^{0,q}(M,F)$. The Bochner-Kodaira-Nakano formula
with a boundary term is as follows:
\begin{thm}[{Ma-Marinescu, \cite[Theorem 1.4.21]{MM}}]%
\label{thm-bkn-mm}
For $q\in \{0,1,\ldots , n\}$ and for $s\in B^{0,q}(M,F)$, we have
\begin{equation}
\begin{split} &\|\overline \partial ^{F} s\|^{2}_{L^{2}(M)}+ \|
\overline \partial ^{F\ast} s\|^{2}_{L^{2}(M)}
\\
=&\|(\nabla ^{\widetilde{F}})^{1,0\ast} \Psi s\|^{2}_{L^{2}(M)}+
\langle \sum _{j,k=1}^{n} R^{F\otimes K_{X}^{\ast}}(\omega _{j},
\overline{\omega }_{k})\overline{\omega }^{k} \wedge \iota _{
\overline{\omega }_{j} }s,s\rangle _{L^{2}(M)}
\\
&-\langle \overline \partial ^{F} s, \Psi ^{-1}\overline{\mathcal{T}}
\Psi s\rangle _{L^{2}(M)}-\langle \Psi ^{-1}\overline{\mathcal{T}}^{
\ast}\Psi s, \overline \partial ^{F\ast} s\rangle _{L^{2}(M)}
\\
&+\langle \overline{\mathcal{T}}^{\ast}\Psi s, (\nabla ^{
\widetilde{F}})^{1,0\ast} \Psi s\rangle _{L^{2}(M)}
\\
&+\int _{bM}\langle \sum _{j,k=1}^{n} \mathscr{L}_{\varrho}(\omega _{j},
\overline{\omega }_{k})\overline{\omega }^{k} \wedge \iota _{
\overline{\omega }_{j} }s,s\rangle _{h} dv_{bM},
\end{split}
\label{eq:3.13Nov}
\end{equation}
where $\{\omega _{j}\}_{j=1}^{n}$ is an orthonormal basis of
$(T^{(1,0)}X, h^{T^{(1,0)}X})$ and $\mathscr{L}_{\varrho}$ is the Levi
form defined in {\eqref{eq:2.2Levi}}. In the boundary term in the last line,
we can replace the sum over the basis $\{\omega _{j}\}_{j=1}^{n}$ by the sum over an orthonormal
basis of $(T^{(1,0)}bM, h^{T^{(1,0)}bM})$ since
$\iota _{e^{(0,1)}_{\mathfrak{n}}}s|_{bM}=0$.
\end{thm}
By substituting the explicit form {\eqref{isopsi}} for the isometry
$\Psi $ into {\eqref{eq:3.13Nov}}, we obtain the following form of the Bochner-Kodaira-Nakano
formula with boundary term, which is due to Andreotti-Vesentini
\cite[p.\,113]{AV65} and Griffiths \cite[(7.14)]{Griffiths1966}.

\begin{thm}
\label{thm3.3}
For $s\in B^{0,q}(M,F)$, we have
\begin{equation}
\begin{split} &\|\overline \partial ^{F} s\|^{2}_{L^{2}(M)}+\|
\overline \partial ^{F*} s\|^{2}_{L^{2}(M)}
\\
=&\sum _{j=1}^{n}\|\widetilde{\nabla }^{TX}_{\overline \omega _{j}}s
\|^{2}_{L^{2}(M)}+\langle \sum _{j,k=1}^{n} R^{F\otimes K^{\ast}_{X}}(
\omega _{j},\overline \omega _{k})\overline \omega ^{k}\wedge \iota _{
\overline{\omega }_{j}}s,s\rangle _{L^{2}(M)}
\\
&+\int _{bM}\langle \sum _{j,k} \mathscr{L}_{\varrho}(\omega _{j},
\overline \omega _{k})\overline \omega ^{k}\wedge \iota _{
\overline \omega _{j}}s,s\rangle _{h} dv_{bM}
\\
&+\langle \overline \partial ^{F} s, \iota _{T^{(0,1)}}s\rangle _{L^{2}(M)}+
\langle \iota ^{*}_{T^{(0,1)}}s,\overline \partial ^{F*}s \rangle _{L^{2}(M)}
\\
&+\sum _{j,k,m=1}^{n}\langle \langle T(\overline \omega _{k},
\overline \omega _{j} ), \omega _{m} \rangle \overline \omega ^{j}
\wedge \iota _{\overline \omega _{m}}s,\widetilde{\nabla }^{TX}_{
\overline \omega _{k}}s \rangle _{L^{2}(M)},
\end{split}
\label{eq:3.17Nov}
\end{equation}
where
$\iota _{T^{(0,1)}}:=\frac{1}{2}\sum _{j,k,m=1}^{n}\langle T(
\overline \omega _{j},\overline \omega _{k}),\omega _{m} \rangle
\overline \omega ^{j}\wedge \overline \omega ^{k}\wedge \iota _{
\overline \omega _{m}}$, and $\iota ^{*}_{T^{(0,1)}}$ denotes the dual
of $\iota _{T^{(0,1)}}$ with respect to the Hermitian metric.
\end{thm}

\subsection{A variant of Bochner-Kodaira-Nakano formula;
Proof of {Theorem~\ref{thm:BKN-variant}}}
\label{ss:3.2variant}

Recall that we always assume $|d\varrho |_{bM}=1$, so that
$|d\varrho |$ is a positive smooth function near $bM$, which might not
be constant. By our definition of $e_{\mathfrak{n}}$, it is clear that
$\sqrt{-1}(e^{(1,0)}_{\mathfrak{n}}-e^{(0,1)}_{\mathfrak{n}})$ is a real
vector field tangent to $bM$; hence, the the derivatives of
$|d\varrho |$ in the normal direction on $bM$ satisfy
\begin{equation}
\frac{1}{2}e_{\mathfrak{n}}(|d\varrho |)|_{bM}=e^{(1,0)}_{
\mathfrak{n}}(|d\varrho |)|_{bM}=e^{(0,1)}_{\mathfrak{n}}(|d\varrho |)|_{bM}
\in \mathscr{C}^{\infty}(bM, \mathbb{R}).
\label{eq3.17}
\end{equation}

Note that the Levi form $\mathscr{L}_{\varrho}$ is defined on
$\Lambda ^{2} T(bM)$. In a local orthonormal frame
$\{\omega _{j}\}_{j=1}^{n-1}$ of $(T^{1,0}bM, h^{T^{1,0}bM})$, we have
\begin{equation}
\mathrm{Tr}_{bM}[\sqrt{-1}\mathscr{L}_{\varrho}]=\sum _{j=1}^{n-1}
\mathscr{L}_{\varrho }(\omega _{j},\overline \omega _{j})\in
\mathscr{C}^{\infty}(bM, \mathbb{R}).
\label{eq3.18}
\end{equation}

\begin{lemma}%
\label{lm:7.3}
Let $\tau \in \mathscr{C}^{\infty}(X,[0,1])$ such that $\tau =1$ near
$bM$ and $\supp \tau $ is included in an open neighborhood of $bM$ where
$d\varrho $ does not vanish. Let
$\Omega _{\tau}\in \Omega ^{n-1,n-1}(X,\mathbb{R})$ be given by
\begin{equation*}
\Omega _{\tau}:=\frac{\Theta ^{n-1}}{(n-1)!}- 2\sqrt{-1}\tau \ast (
\partial \varrho \wedge \bar{\partial}\varrho ),
\end{equation*}
where $\ast $ denotes the Hodge star operator associated to $g^{TX}$. Then
we have
\begin{equation}
\frac{\sqrt{-1}}{2}(\partial -\bar{\partial}) \Omega _{\tau }|_{bM}=
\left ( \mathrm{Tr}_{bM}[\sqrt{-1}\mathscr{L}_{\varrho}]-e_{
\mathfrak{n}}(|d\varrho |)\right )dv_{bM},
\label{eq:7.11B}
\end{equation}
and
\begin{equation}
\Omega _{\tau}|_{bM}=0.
\label{eq:3.21Nov}
\end{equation}
Moreover, if $n\geq 2$, we can choose the function $\tau $ as above such
that
\begin{equation}
\Theta \wedge \Omega _{\tau }> 0 \,,\,\text{ near } M.
\label{eq:7.13B}
\end{equation}
If $n=1$, for any $\varepsilon >0$, we can choose the function
$\tau $ as above such that
\begin{equation}
\Theta \wedge \Omega _{\tau }+\varepsilon dv_{X}> 0 \,,\,
\text{ near } M.
\label{eq:7.14B}
\end{equation}
\end{lemma}
\begin{proof}
At first, we consider a special local frame of $T^{(1,0)}X$ near
$bM$. Let $\{\omega _{j}\}_{j=1}^{n}$ be an orthonormal frame of
$(T^{(1,0)}X, h^{T^{(1,0)}X})$; in particular, near $bM$, we take
\begin{equation}
\omega _{n}=\frac{\sqrt{2}}{|d\varrho |} e^{(1,0)}_{\mathfrak{n}},
\label{eq:3.24Nov}
\end{equation}
where $1/|d\varrho |>0$ is a smooth function near $bM$ with
$(1/|d\varrho |)|_{bM}\equiv 1$. Let $\{\omega ^{j}\}_{j=1}^{n}$ be the
dual frame of $T^{(1,0)\ast}X$; in particular, near $bM$, we have
\begin{equation*}
\omega ^{n}=-\frac{\sqrt{2}}{|d\varrho |}\partial \varrho .
\end{equation*}
As a consequence, the following $(1,1)$-form acting on
$T^{1,0}bM\otimes T^{0,1}bM$ satisfies
\begin{equation}
\omega ^{n}\wedge \overline{\omega }^{n}|_{bM}=0.
\label{eq3.24}
\end{equation}
Here the restriction is the pull-backs under the inclusion
$i_{bM}:bM\to X$.

Then for $k\in \{1,\ldots , n-1\}$ and near $bM$, we have
\begin{equation}
\omega _{k}(\varrho )=\overline{\omega }_{k}(\varrho )=0,
\label{eq3.25}
\end{equation}
and
\begin{equation*}
\omega _{n}(\varrho )=\overline{\omega }_{n}(\varrho )= -
\frac{|d\varrho |}{\sqrt{2}}\cdot
\end{equation*}
Therefore, at $bM$, $\{\omega _{j}\}_{j=1}^{n-1}$ forms a local orthonormal
frame of $(T^{(1,0)}bM, h^{T^{(1,0)}bM})$. Moreover, near $bM$, we have
\begin{equation}
\partial \bar{\partial}\varrho (\omega _{k},\overline{\omega }_{k})=
\begin{cases}
& -( \bar{\partial}\varrho , [\omega _{k},\overline{\omega }_{k}])=
\langle e^{(1,0)}_{\mathfrak{n}},\widetilde{\nabla}^{TX}_{\omega _{k}}
\overline{\omega }_{k}\rangle , \text{ for } k<n;
\\[4pt]
&-\dfrac{1}{|d\varrho |}e_{\mathfrak{n}}(|d\varrho |)+
\dfrac{2}{|d\varrho |^{2}}\langle e^{(1,0)}_{\mathfrak{n}},
\widetilde{\nabla}^{TX}_{e^{(1,0)}_{\mathfrak{n}}}e^{(0,1)}_{
\mathfrak{n}}\rangle , \text{ for } k=n.
\end{cases}
\label{eq3.26}
\end{equation}
Note that the Riemannian volume form associated with $g^{TX}$ is given
as
\begin{equation}
dv_{X}=\frac{\Theta ^{n}}{n!}=\Lambda _{j=1}^{n} (\sqrt{-1}\omega ^{j}
\wedge \overline{\omega }^{j}).
\label{eq:3.28Nov-1}
\end{equation}
For $k\in \{1,\ldots , n\}$, we denote the local real $(n-1,n-1)$ form
\begin{equation}
\Psi _{k\bar{k}}=\sqrt{-1}\iota _{\omega _{k}}\iota _{
\overline{\omega }_{k}}dv_{X}.
\label{eq3.28}
\end{equation}
Then we have
\begin{equation}
\Omega _{\tau}=\sum _{k=1}^{n-1}\Psi _{k\bar{k}}+(1-\tau |d\varrho |^{2})
\Psi _{n\bar{n}}.
\label{eq:3.28Nov-2}
\end{equation}
It follows directly that $\Omega _{\tau}|_{bM}=0$. We have
\begin{equation}
\Theta \wedge \Omega _{\tau }=(n-\tau |d\varrho |^{2})dv_{X}.
\label{eq3.30}
\end{equation}
If $n\geq 2$, we can choose $\tau $ such that
$n-\tau |d\varrho |^{2} > 0$ near the domain $M$. Therefore we obtain {\eqref{eq:7.13B}}. When $n=1$, for any $\varepsilon >0$, we can choose
$\tau $ appropriately such that
$1-\tau |d\varrho |^{2}>-\varepsilon $ near the domain $M$, thereby satisfying {\eqref{eq:7.14B}}.

It remains to prove equation {\eqref{eq:7.11B}}. Initially, it is clear that
\begin{equation}
\iota _{\omega _{n}}dv_{X}|_{bM}=-\frac{1}{\sqrt{2}}dv_{bM}, \,
\iota _{\overline{\omega }_{n}}dv_{X}|_{bM}=-\frac{1}{\sqrt{2}}dv_{bM}.
\label{eq:3.32Novs}
\end{equation}
For $1\leq k \leq n-1$, we compute
\begin{equation}
\begin{split} \partial \Psi _{k\bar{k}}|_{bM} &=-\sqrt{-1}\langle
\partial \overline{\omega }^{n}, \omega _{k}\wedge \overline{\omega }_{k}
\rangle \iota _{\overline{\omega }_{n}}dv_{X}|_{bM}
\\
&=\frac{\sqrt{-1}}{\sqrt{2}}\langle \partial \overline{\omega }^{n},
\omega _{k}\wedge \overline{\omega }_{k}\rangle dv_{bM}
\\
&= -\sqrt{-1}\partial \bar{\partial}\varrho(\omega _{k},\overline{\omega }_{k})
dv_{bM}.
\end{split}
\label{eq3.32}
\end{equation}
Similarly, we have
\begin{equation}
\bar{\partial}\Psi _{k\bar{k}}|_{bM} = \sqrt{-1} \partial
\bar{\partial}\varrho(\omega _{k},\overline{\omega }_{k}) dv_{bM}.
\label{eq3.33}
\end{equation}
Then we have
\begin{equation}
\frac{\sqrt{-1}}{2}(\partial -\bar{\partial})\Psi _{k\bar{k}}|_{bM}=
\partial \bar{\partial}\varrho(\omega _{k},\overline{\omega }_{k}) dv_{bM}.
\label{eq:3.33Nov}
\end{equation}
Now we deal with the last term
$(1-\tau |d\varrho |^{2})\Psi _{n\bar{n}}$. We have
\begin{equation}
\begin{split} (\partial -\bar{\partial})(1-\tau |d\varrho |^{2})\Psi _{n
\bar{n}}|_{bM}&=2(-\omega _{n}(|d\varrho |)\sqrt{-1}\iota _{
\overline{\omega }_{n}}dv_{X}-\overline{\omega }_{n}(|d\varrho |)
\sqrt{-1}\iota _{\omega _{n}}dv_{X})|_{bM}
\\
&=2\sqrt{-1}e_{\mathfrak{n}}(|d\varrho |)dv_{bM}.
\end{split}
\label{eq:3.34Nov}
\end{equation}
Then {\eqref{eq:7.11B}} follows from {\eqref{eq:3.28Nov-2}},
{\eqref{eq:3.33Nov}} and {\eqref{eq:3.34Nov}}. This finishes the proof.
\end{proof}

As a consequence, we obtain a slight extension of
\cite[Theorem 7.3]{Griffiths1966}.
\begin{cor}%
\label{cor:3.5new}
Let $\tau \in \mathscr{C}^{\infty}(X,[0,1])$ and $\Omega _{\tau}$ be as
in {Lemma~\ref{lm:7.3}}. Then for
$g\in \mathscr{C}^{\infty}(\overline{M})$, we have
\begin{equation}
\int _{M} \sqrt{-1}\partial \bar{\partial}g\wedge \Omega _{\tau }-
\int _{M} g\sqrt{-1}\partial \bar{\partial}\Omega _{\tau }=\int _{bM} g
\left ( \mathrm{Tr}_{bM}[\sqrt{-1}\mathscr{L}_{\varrho}]-e_{
\mathfrak{n}}(|d\varrho |)\right ) dv_{bM}.
\label{eq:3.35Nov}
\end{equation}
\end{cor}
\begin{proof}
By {Lemma~\ref{lm:7.3}} and Stokes' formula, we get
\begin{equation}
\begin{split} \int _{bM} g\left ( \mathrm{Tr}_{bM}[\sqrt{-1}
\mathscr{L}_{\varrho}]-e_{\mathfrak{n}}(|d\varrho |)\right )dv_{bM}&=
\frac{1}{2}\int _{bM}g\sqrt{-1}(\partial -\bar{\partial}) \Omega _{
\tau}
\\
&=\frac{1}{2}\int _{M}\sqrt{-1}d\left (g(\partial -\bar{\partial})
\Omega _{\tau}\right )
\\
&=\int _{M} \sqrt{-1}\partial \overline{\partial }g\wedge \Omega _{
\tau}-\int _{M} g\sqrt{-1}\partial \bar{\partial}\Omega _{\tau }
\\
&\qquad +\frac{\sqrt{-1}}{2}\int _{bM} (\partial -
\overline{\partial })g\wedge \Omega _{\tau}.
\end{split}
\label{eq3.37}
\end{equation}
By {\eqref{eq:3.21Nov}}, the last integral on $bM$ vanishes, so that we obtain {\eqref{eq:3.35Nov}}.
\end{proof}

Prior to presenting a more complex variant of the Bochner-Kodaira-Nakano
formula with boundary terms, we introduce some notation to streamline
our expressions.

\begin{defn}[Modified curvature commutators]%
\label{defn:bracket}
For $\varepsilon \in [0,1]$, we define the following operator on
$\Omega ^{p,q}(X,F)$,
\begin{equation}
\label{eq:MLB1}
[\sqrt{-1}R^{F},\Lambda ]_{\varepsilon}:= [\sqrt{-1}R^{F},\Lambda ]+
\varepsilon \mathrm{Tr}_{\Theta}[\sqrt{-1}R^{F}],
\end{equation}
where $\mathrm{Tr}_{\Theta}[\sqrt{-1}R^{F}]$ is defined in {\eqref{eq:3.11Nov}}. Similarly, we define the operator acting on
$\mathcal{D}^{p,q}(bM,F)$ (see {\eqref{eqdecomp}}) by
\begin{equation}
\label{eq:MLB2}
[\sqrt{-1}\mathscr{L}_{\varrho},\Lambda ]^{bM}_{\varepsilon}:= [\sqrt{-1}
\mathscr{L}_{\varrho},\Lambda |_{bM}]+ \varepsilon \mathrm{Tr}_{bM}[
\sqrt{-1}\mathscr{L}_{\varrho}].
\end{equation}
Note that, for $s\in B^{0,q}(M,F)$, we have
$s|_{bM}\in \mathcal{D}^{0,q}(bM,F)$.
\end{defn}

We will also need the following notation which appeared in {Theorem~\ref{thm:BKN-variant}}. Recall that for a $1$-form $\beta $ on $X$,
$\mathrm{Tr}[\nabla \beta ](x)$ is defined with respect to
$(TX, g^{TX})$. If $\{e_{j}\}_{j=1}^{2n}$ is an orthonormal
basis of the real tangent space $T_xX$, then
\begin{equation*}
\mathrm{Tr}[\nabla \beta ](x)=\sum _{j} e_{j} \beta (e_{j})(x)-
\beta (\sum _{j} \nabla ^{TX}_{e_{j}}e_{j})(x).
\end{equation*}
Let $\widetilde{R}^{\Lambda ^{(0,\bullet )}T^{\ast }X}$ denote the curvature
tensor operator associated with the Hermitian connection
$\widetilde{\nabla}^{TX}$ on $\Lambda ^{(0,\bullet )}T^{\ast }X$.
\begin{defn}%
\label{def:3.7QFT}
Given $\tau \in \mathscr{C}^{\infty}(X,[0,1])$, $\Omega _{\tau}$ as
in {Lemma~\ref{lm:7.3}} and $\varepsilon \in [0,1]$; a local orthonormal frame of $(T^{1,0}X, h^{T^{1,0}X})$ denoted by $\{\omega _{j}\}_{j=1}^{n}$. For
$s\in B^{0,q}(M,F)$, we define the following bilinear forms
\begin{equation}
\begin{split} Q^{F}_{T,\varepsilon}(s,s):=&\langle \overline \partial ^{F}
s, \iota _{T^{(0,1)}}s\rangle _{L^{2}(M)}+\langle \iota ^{*}_{T^{(0,1)}}s,
\overline \partial ^{F*}s \rangle _{L^{2}(M)}
\\
&+\sum _{j,k,m=1}^{n}\left \langle \langle T(\overline \omega _{k},
\overline \omega _{j} ), \omega _{m} \rangle \overline \omega ^{j}
\wedge \iota _{\overline \omega _{m}}s,\widetilde{\nabla }^{TX}_{
\overline \omega _{k}}s \right \rangle _{L^{2}(M)}
\\
&-2(1-\varepsilon ){\mathrm{Re}}\int _{M} \sum _{j,k=1}^{n} \langle T(
\omega _{j},\omega _{k}),\overline{\omega }_{k}\rangle \langle
\widetilde{\nabla}^{TX}_{\overline{\omega }_{j}}s,s\rangle _{h} dv_{X},
\end{split}
\label{eq:3.40Nov}
\end{equation}
and
\begin{equation}
\begin{split} \Psi ^{F}_{\tau}(s,s):=&-\left \langle \mathrm{Tr}_{
\Theta}[\sqrt{-1}\widetilde{R}^{\Lambda ^{(0,\bullet )}T^{\ast }X}] s,s
\right \rangle _{L^{2}(M)}
\\
&+\left \langle \tau |d\varrho |^{2}\left (\widetilde{R}^{\Lambda ^{(0,
\bullet )}T^{\ast }X}\otimes \mathrm{Id}_{F}\right )(\omega _{n},
\overline{\omega }_{n})s,s\right \rangle _{L^{2}(M)}
\\
& +2{\mathrm{Re}}\int _{M} \left (\omega _{n}(\tau |d\varrho |^{2}) +\tau |d
\varrho |^{2}\mathrm{Tr}[\nabla \overline{\omega }^{n}]\right )
\langle \widetilde{\nabla}^{TX}_{\overline{\omega }_{n}}s,s\rangle _{h}
dv_{X}
\\
&+2{\mathrm{Re}}\int _{M} \tau |d\varrho |^{2}\langle \widetilde{\nabla}^{TX}_{
\widetilde{\nabla}^{TX}_{\omega _{n}}\overline{\omega }_{n}}s,s
\rangle _{h} dv_{X},
\end{split}
\label{eq3.41}
\end{equation}
where $\omega _{n}$ near $bM$ (on the support of $\tau $) is taken as in {\eqref{eq:3.24Nov}}, that is,
$\omega _{n}=\frac{\sqrt{2}}{|d\varrho |} e^{(1,0)}_{\mathfrak{n}}$.
\end{defn}

The following proposition, whose proof will be given in Subsection \ref{ss:3.4Nov}, is the key ingredient for the proof of {Theorem~\ref{thm:BKN-variant}}.
\begin{prop}%
\label{prop:3.8}
Given $\tau \in \mathscr{C}^{\infty}(X,[0,1])$ and $\Omega _{\tau}$ as
in {Lemma~\ref{lm:7.3}}, then for
$s\in \Omega ^{0,q}(\overline{M},F)$, we have
\begin{align}
 \int _{M}\sqrt{-1}\partial \bar{\partial}|s|^{2}_{h}
\wedge \Omega _{\tau }=&\sum _{j=1}^{n} \|\widetilde{\nabla}^{TX}_{
\omega _{j}}s\|^{2}_{L^{2}(M)}-\sum _{j=1}^{n} \|\widetilde{\nabla}^{TX}_{
\overline{\omega }_{j}}s\|^{2}_{L^{2}(M)}-\left \langle \mathrm{Tr}_{
\Theta}[\sqrt{-1}R^{F}]s,s\right \rangle _{L^{2}(M)}
\notag \\
&-2\|\sqrt{\tau}\widetilde{\nabla}^{TX}_{e^{(1,0)}_{\mathfrak{n}}}s\|^{2}_{L^{2}(M)}+2
\|\sqrt{\tau}\widetilde{\nabla}^{TX}_{e^{(0,1)}_{\mathfrak{n}}}s\|^{2}_{L^{2}(M)}
\notag \\
& -2{\mathrm{Re}}\int _{M} \sum _{j,k=1}^{n} \langle T(\omega _{j},\omega _{k}),
\overline{\omega }_{k}\rangle \langle \widetilde{\nabla}^{TX}_{
\overline{\omega }_{j}}s,s\rangle _{h} dv_{X}
\notag \\
&+2\langle \tau R^{F}(e^{(1,0)}_{\mathfrak{n}}, e^{(0,1)}_{
\mathfrak{n}})s,s\rangle _{L^{2}(M)}+ \Psi ^{F}_{\tau}(s,s).
\label{eq3.42}
\end{align}
\end{prop}

Now we are ready to prove {Theorem~\ref{thm:BKN-variant}}.
\begin{proof}[Proof of {Theorem~\ref{thm:BKN-variant}}]
In this proof, we always use the local orthonormal frame
$\{\omega _{j}\}_{j=1}^{n}$ of $(T^{1,0}X, h^{T^{1,0}X})$ constructed
at the beginning of the proof of {Lemma~\ref{lm:7.3}}, where we have
$\omega _{n}=\frac{\sqrt{2}}{|d\varrho |} e^{(1,0)}_{\mathfrak{n}}$ near
$bM$ (on the support of $\tau $). Applying {\eqref{eq:3.35Nov}} with
$g=|s|^{2}_{h}$ and then combining it with {Proposition~\ref{prop:3.8}},
we get
\begin{equation}
\begin{split} \sum _{j=1}^{n} \|\widetilde{\nabla}^{TX}_{
\overline{\omega }_{j}}s\|^{2}_{L^{2}(M)} =&-\int _{M}|s|^{2}_{h}
\sqrt{-1}\partial \bar{\partial}\Omega _{\tau }+\sum _{j=1}^{n} \|
\widetilde{\nabla}^{TX}_{\omega _{j}}s\|^{2}_{L^{2}(M)}
\\
&-\|\sqrt{\tau}|d\varrho |\widetilde{\nabla}^{TX}_{\omega _{n}}s\|^{2}_{L^{2}(M)}+
\|\sqrt{\tau}|d\varrho |\widetilde{\nabla}^{TX}_{\overline{\omega }_{n}}s
\|^{2}_{L^{2}(M)}
\\
& -2{\mathrm{Re}}\int _{M} \sum _{j,k=1}^{n} \langle T(\omega _{j},\omega _{k}),
\overline{\omega }_{k}\rangle \langle \widetilde{\nabla}^{TX}_{
\overline{\omega }_{j}}s,s\rangle _{h} dv_{X}
\\
&-\left \langle \sum _{j=1}^{n} R^{F}(\omega _{j},\overline{\omega }_{j})s,s
\right \rangle _{L^{2}(M)}+\left \langle \tau |d\varrho |^{2} R^{F}(
\omega _{n},\overline{\omega }_{n})s,s\right \rangle _{L^{2}(M)}
\\
&+\Psi ^{F}_{\tau}(s,s)-\int _{bM}|s|^{2}_{h} \left (\mathrm{Tr}_{bM}[\sqrt{-1}
\mathscr{L}_{\varrho}]-e_{\mathfrak{n}}(|d\varrho |)\right ) dv_{bM}.
\end{split}
\label{eq:3.44proof}
\end{equation}
Now we take $s\in B^{0,q}(M,F)$. We split
\begin{equation*}
\sum _{j=1}^{n} \|\widetilde{\nabla}^{TX}_{\overline{\omega }_{j}}s\|^{2}_{L^{2}(M)}=
\varepsilon \sum _{j=1}^{n} \|\widetilde{\nabla}^{TX}_{
\overline{\omega }_{j}}s\|^{2}_{L^{2}(M)} + (1-\varepsilon ) \sum _{j=1}^{n}
\|\widetilde{\nabla}^{TX}_{\overline{\omega }_{j}}s\|^{2}_{L^{2}(M)},
\end{equation*}
and apply {\eqref{eq:3.44proof}} to the term
$(1-\varepsilon ) \sum _{j=1}^{n} \|\widetilde{\nabla}^{TX}_{
\overline{\omega }_{j}}s\|^{2}_{L^{2}(M)}$ and insert it back to {\eqref{eq:3.17Nov}} to obtain {\eqref{eq:3.42Nov}}.
\end{proof}

\subsection{Nakano-Griffiths inequalities; Proof of {Theorem~\ref{thm:NG-inequality}}}
\label{sec3.3}

In this subsection, we prove {Theorem~\ref{thm:NG-inequality}}. We will use
the following elementary inequality: for any $a,b\in \mathbb{C}$,
$\varepsilon >0$,
\begin{equation}
2|ab|\leq \varepsilon |a|^{2} +\frac{1}{\varepsilon}|b|^{2}.
\label{eq3.44}
\end{equation}

\begin{proof}[Proof of {Theorem~\ref{thm:NG-inequality}}]
First, fix the function $\tau $ as in the theorem. Then there exists a constant
$C_{1}=C_{1}(X,g^{TX},\varrho ,\tau )>0$ such that for all
$\varepsilon \in [0,1]$ and $s\in \Omega ^{0,q}(\overline{M},F)$,
\begin{equation}
\left | \big\langle \sum _{j,k} R^{K^{\ast}_{X}}(\omega _{j},
\overline \omega _{k})\overline{\omega }^{k}\wedge \iota _{
\overline{\omega }_{j}}s,s\big\rangle _{L^{2}(M)}-(1-\varepsilon )
\int _{M} |s|^{2}_{h} \sqrt{-1}\partial \overline{\partial }\Omega _{
\tau}\right |\leq C_{1}\|s\|^{2}_{L^{2}(M)}.
\label{eq:3.46Nov}
\end{equation}

Estimating separately the terms in $Q^{F}_{T,\varepsilon}(s,s)$ involving $\overline \partial ^{F} s$, $\overline \partial ^{F\ast} s$, and the covariant derivatives
$\{\widetilde{\nabla }^{TX}_{\overline{\omega }_{j} } s\}_{j=1}^{n}$, we obtain that there exists a constant
$C_{2}=C_{2}(X,g^{TX},\varrho ,\tau )>0$ such that for all
$\varepsilon \in (0,1]$ and $s\in \Omega ^{0,q}(\overline{M},F)$,
\begin{equation}
\left | Q^{F}_{T,\varepsilon}(s,s)\right |\leq
\frac{C_{2}}{\varepsilon}\|s\|^{2}_{L^{2}(M)}+\varepsilon \left (\|
\overline \partial ^{F} s\|^{2}_{L^{2}(M)}+\|\overline \partial ^{F
\ast} s\|^{2}_{L^{2}(M)}\right ) +\frac{\varepsilon}{4}\sum _{j=1}^{n}
\|\widetilde{\nabla}^{TX}_{\overline{\omega }_{j}}s\|^{2}_{L^{2}(M)} .
\label{eq:3.47Nov}
\end{equation}
Note that $C_{2}$ is determined only by the torsion tensor $T$.

Finally, note that inside $\Psi ^{F}_{\tau}(s,s)$, we can write
\begin{equation}
\widetilde{\nabla }^{TX}_{\widetilde{\nabla }^{TX}_{\omega _{n}}
\overline{\omega }_{n}}s=\sum _{j=1}^{n}\langle \widetilde{\nabla }^{TX}_{
\omega _{n}}\overline{\omega }_{n},\omega _{j}\rangle
\widetilde{\nabla }^{TX}_{\overline{\omega }_{j}}s.
\label{eq3.47}
\end{equation}
Since our choice of $\tau $ ensures that
$\omega _{n}$ is well defined via $e_{\mathfrak{n}}$, we conclude that there exists
a constant $C_{3}=C_{3}(X,g^{TX},\varrho ,\tau )>0$ such that for all
$\varepsilon \in (0,1]$ and $s\in \Omega ^{0,q}(\overline{M},F)$,
\begin{equation}
\left | \Psi ^{F}_{\tau}(s,s)\right |\leq \frac{C_{3}}{\varepsilon}\|s
\|^{2}_{L^{2}(M)}+\frac{\varepsilon}{4}\sum _{j=1}^{n} \|
\widetilde{\nabla}^{TX}_{\overline{\omega }_{j}}s\|^{2}_{L^{2}(M)}.
\label{eq:3.49Nov}
\end{equation}

It is clear that the constants $C_{1}$, $C_{2}$, $C_{3}$ depend only on
the geometry of $(X,g^{TX})$ and the functions $\tau $ and
$\varrho $ (hence also depend on $M$). Finally, by combining {\eqref{eq:3.46Nov}}--{\eqref{eq:3.49Nov}}, we obtain {\eqref{eq:3.44Nov}} from {\eqref{eq:3.42Nov}}.
\end{proof}
\begin{rem}
\label{rem3.9}
From the above proof, {\eqref{eq:3.44Nov}} also holds if we add two more
positive terms to the right-hand side
\begin{equation*}
\frac{\varepsilon}{2}\sum _{j=1}^{n} \|\widetilde{\nabla}^{TX}_{
\overline{\omega }_{j}}s\|^{2}_{L^{2}(M)} + 2(1-\varepsilon )\|\sqrt{
\tau}\widetilde{\nabla}^{TX}_{e^{(0,1)}_{\mathfrak{n}}}s\|^{2}_{L^{2}(M)}.
\end{equation*}
\end{rem}

\subsection{Proof of {Proposition~\ref{prop:3.8}}}
\label{ss:3.4Nov}

Now we focus on the proof of {Proposition~\ref{prop:3.8}}. Let
$\{\omega _{j}\}_{j=1}^{n}$ denote a local orthonormal frame of
$(T^{1,0}X, h^{T^{1,0}X})$. Given
$s\in \Omega ^{0,q}(\overline{M},F)$, we have
\begin{equation}
\sqrt{-1}\partial \bar{\partial}|s|^{2}_{h}\wedge \Omega _{\tau }=
\sum _{j=1}^{n}(\partial \bar{\partial}|s|^{2}_{h})(\omega _{j},
\overline{\omega }_{j}) dv_{X}-2\tau (\partial \bar{\partial}|s|^{2}_{h})(
e^{(1,0)}_{\mathfrak{n}}, e^{(0,1)}_{\mathfrak{n}})dv_{X}.
\label{eq:3.51Novs}
\end{equation}
Note that by $T(\omega _{j},\overline{\omega }_{j})=0$, we have
\begin{equation*}
[\omega _{j},\overline{\omega }_{j}]=\widetilde{\nabla }^{TX}_{
\omega _{j}}\overline{\omega }_{j}-\widetilde{\nabla }^{TX}_{
\overline{\omega }_{j}}\omega _{j}.
\end{equation*}
For $j=1,\ldots , n$, we have
\begin{equation}
\begin{split} (\partial \bar{\partial}|s|^{2}_{h})(\omega _{j},
\overline{\omega }_{j})=&\overline{\omega }_{j} (\omega _{j}(|s|^{2}_{h}))-(
\partial |s|^{2}_{h}, \widetilde{\nabla}^{TX}_{\overline{\omega }_{j}}
\omega _{j})
\\
=&|\widetilde{\nabla}^{TX}_{\omega _{j}}s|^{2}_{h}+|
\widetilde{\nabla}^{TX}_{\overline{\omega }_{j}}s|^{2}_{h}
\\
&-\left \langle \left (\widetilde{R}^{\Lambda ^{(0,\bullet )}T^{\ast }X}
\otimes \mathrm{Id}_{F}+\mathrm{Id}_{\Lambda ^{(0,\bullet )}T^{\ast }X}
\otimes R^{F}\right )(\omega _{j},\overline{\omega }_{j})s,s\right
\rangle _{h}
\\
&+2{\mathrm{Re}}\left \langle \left (\widetilde{\nabla}^{TX}_{\omega _{j}}
\widetilde{\nabla}^{TX}_{\overline{\omega }_{j}}-\widetilde{\nabla}^{TX}_{
\widetilde{\nabla}^{TX}_{\omega _{j}}\overline{\omega }_{j}}\right )s,s
\right \rangle _{h}.
\end{split}
\label{eq:3.52Novs}
\end{equation}
We need the following technical lemmas.
\begin{lemma}%
\label{lm:3.14Nov}
In the above setup, we have
\begin{equation}
\begin{split} &\int _{M} \sum _{j=1}^{n}\left \langle \left (
\widetilde{\nabla}^{TX}_{\omega _{j}}\widetilde{\nabla}^{TX}_{
\overline{\omega }_{j}}-\widetilde{\nabla}^{TX}_{\widetilde{\nabla}^{TX}_{
\omega _{j}}\overline{\omega }_{j}}\right )s,s\right \rangle _{h} dv_{X}
\\
=&-\sum _{j=1}^{n} \|\widetilde{\nabla}^{TX}_{\overline{\omega }_{j}}s
\|^{2}_{L^{2}(M)} - \int _{bM} \langle \widetilde{\nabla}^{TX}_{e^{(0,1)}_{
\mathfrak{n}}}s,s\rangle _{h} dv_{bM} -\int _{M} \sum _{j,k=1}^{n}
\langle T(\omega _{j},\omega _{k}),\overline{\omega }_{k}\rangle
\langle \widetilde{\nabla}^{TX}_{\overline{\omega }_{j}}s,s\rangle _{h}
dv_{X}.
\end{split}
\label{eq:3.51Nov}
\end{equation}
\end{lemma}
\begin{proof}
Note that we cannot choose the local orthonormal frame
$\{\omega _{j}\}_{j=1}^{n}$ globally on $\overline{M}$, but the integrands
in {\eqref{eq:3.51Nov}} are clearly well-defined global functions on
$\overline{M}$.

First, we have
\begin{equation}
L_{\omega _{j}}(\langle \widetilde{\nabla }^{TX}_{\overline{\omega }_{j}}s,s
\rangle _{h} dv_{X})=\left \langle \widetilde{\nabla}^{TX}_{\omega _{j}}
\widetilde{\nabla}^{TX}_{\overline{\omega }_{j}}s,s\right \rangle _{h}
dv_{X}+|\widetilde{\nabla }^{TX}_{\overline{\omega }_{j}}s|^{2}_{h} dv_{X}
+ \langle \widetilde{\nabla }^{TX}_{\overline{\omega }_{j}}s,s
\rangle _{h} L_{\omega _{j}}dv_{X},
\label{eq3.52}
\end{equation}
where $L_{\omega _{j}}$ denotes the Lie derivative on differential forms.

Using the formula {\eqref{eq:3.28Nov-1}} for $dv_{X}$, we obtain
\begin{equation}
L_{\omega _{j}} dv_{X}=\sum _{k=1}^{n}(\langle T(\omega _{j},\omega _{k}),
\overline{\omega }_{k}\rangle -\langle \omega _{j},
\widetilde{\nabla }^{TX}_{\omega _{k}}\overline{\omega }_{k}\rangle ) dv_{X}.
\label{eq3.53}
\end{equation}

Then we have
\begin{equation}
\begin{split} \sum _{j=1}^{n} \left \langle \left (\widetilde{\nabla}^{TX}_{
\omega _{j}}\widetilde{\nabla}^{TX}_{\overline{\omega }_{j}}-
\widetilde{\nabla}^{TX}_{\widetilde{\nabla}^{TX}_{\omega _{j}}
\overline{\omega }_{j}}\right )s,s\right \rangle _{h} dv_{X}=& \sum _{j=1}^{n}
L_{\omega _{j}}(\langle \widetilde{\nabla }^{TX}_{\overline{\omega }_{j}}s,s
\rangle _{h} dv_{X})-\sum _{j=1}^{n} |\widetilde{\nabla }^{TX}_{
\overline{\omega }_{j}}s|^{2}_{h} dv_{X}
\\
& -\sum _{j,k=1}^{n} \langle T(\omega _{j},\omega _{k}),
\overline{\omega }_{k}\rangle \langle \widetilde{\nabla }^{TX}_{
\overline{\omega }_{j}}s,s\rangle _{h} dv_{X}.
\end{split}
\label{eq3.54}
\end{equation}
It remains to prove that
\begin{equation}
\int _{M} \sum _{j=1}^{n} L_{\omega _{j}} (\langle
\widetilde{\nabla }^{TX}_{\overline{\omega }_{j}}s,s\rangle _{h} dv_{X})=
- \int _{bM} \langle \widetilde{\nabla}^{TX}_{e^{(0,1)}_{\mathfrak{n}}}s,s
\rangle _{h} dv_{bM}
\label{eq:3.56Nov}
\end{equation}
We now use the local orthonormal frame $\{\omega _{j}\}_{j=1}^{n}$ of
$(T^{1,0}X, h^{T^{1,0}X})$ constructed at the beginning of the proof
of {Lemma~\ref{lm:7.3}}, where
$\omega _{n}=\frac{\sqrt{2}}{|d\varrho |} e^{(1,0)}_{\mathfrak{n}}$ near
$bM$ (on the support of $\tau $). Using the fact that $\omega _{j}$,
$j=1,\ldots , n-1$, are tangential to $bM$ and Stokes' theorem, we get
\begin{equation}
\int _{M} \sum _{j=1}^{n} L_{\omega _{j}}(\langle \widetilde{\nabla }^{TX}_{
\overline{\omega }_{j}}s,s\rangle _{h} dv_{X})=\int _{M} d\iota _{
\omega _{n}}(\langle \widetilde{\nabla }^{TX}_{\overline{\omega }_{n}}s,s
\rangle _{h} dv_{X})=\int _{bM} \iota _{\omega _{n}}(\langle
\widetilde{\nabla }^{TX}_{\overline{\omega }_{n}}s,s\rangle _{h} dv_{X}),
\label{eq3.56}
\end{equation}
then {\eqref{eq:3.56Nov}} follows from {\eqref{eq:3.24Nov}} and {\eqref{eq:3.32Novs}}, and this finishes the proof.
\end{proof}

Combining the above lemma with {\eqref{eq:3.52Novs}}, we obtain the following
result.
\begin{lemma}%
\label{lm:3.15}
For $s\in \Omega ^{0,q}(\overline{M},F)$, we have
\begin{align}
 &\int _{M} \left (\sum _{j=1}^{n}(\partial
\bar{\partial}|s|^{2}_{h})(\omega _{j},\overline{\omega }_{j})\right )dv_{X}
\notag \\
=&\sum _{j=1}^{n} \|\widetilde{\nabla}^{TX}_{\omega _{j}}s\|^{2}_{L^{2}(M)}-
\sum _{j=1}^{n} \|\widetilde{\nabla}^{TX}_{\overline{\omega }_{j}}s\|^{2}_{L^{2}(M)}-
2{\mathrm{Re}}\int _{bM} \langle \widetilde{\nabla}^{TX}_{e^{(0,1)}_{
\mathfrak{n}}}s,s\rangle _{h} dv_{bM}
\notag \\
& -2{\mathrm{Re}}\int _{M} \sum _{j,k=1}^{n} \langle T(\omega _{j},\omega _{k}),
\overline{\omega }_{k}\rangle \langle \widetilde{\nabla}^{TX}_{
\overline{\omega }_{j}}s,s\rangle _{h} dv_{X}
\notag \\
&-\left \langle \sum _{j=1}^{n}\left (\widetilde{R}^{\Lambda ^{(0,
\bullet )}T^{\ast }X}\otimes \mathrm{Id}_{F}+\mathrm{Id}_{\Lambda ^{(0,
\bullet )}T^{\ast }X}\otimes R^{F}\right )(\omega _{j},
\overline{\omega }_{j})s,s\right \rangle _{L^{2}(M)}.
\label{eq3.57}
\end{align}
\end{lemma}

\begin{lemma}%
\label{lm:3.16}
Let $\tau \in \mathscr{C}^{\infty}(X,[0,1])$ and $\Omega _{\tau}$ be as
in {Lemma~\ref{lm:7.3}} and let $\omega _{n}$ near $bM$ (on the support of
$\tau $) be as in {\eqref{eq:3.24Nov}}, that is,
$\omega _{n}=\frac{\sqrt{2}}{|d\varrho |} e^{(1,0)}_{\mathfrak{n}}$. Then
\begin{equation}
\begin{split} &\int _{M} \tau |d\varrho |^{2}(\partial \bar{\partial}|s|^{2}_{h})(
\omega _{n},\overline{\omega }_{n})dv_{X}
\\
=& \|\sqrt{\tau}|d\varrho |\widetilde{\nabla}^{TX}_{\omega _{n}}s\|^{2}_{L^{2}(M)}-
\|\sqrt{\tau}|d\varrho |\widetilde{\nabla}^{TX}_{\overline{\omega }_{n}}s
\|^{2}_{L^{2}(M)}- 2{\mathrm{Re}}\int _{bM} \langle \widetilde{\nabla}^{TX}_{e^{(0,1)}_{
\mathfrak{n}}}s,s\rangle _{h} dv_{bM}
\\
&-\left \langle \tau |d\varrho |^{2}\left (\widetilde{R}^{\Lambda ^{(0,
\bullet )}T^{\ast }X}\otimes \mathrm{Id}_{F}+\mathrm{Id}_{\Lambda ^{(0,
\bullet )}T^{\ast }X}\otimes R^{F}\right )(\omega _{n},
\overline{\omega }_{n})s,s\right \rangle _{L^{2}(M)}
\\
& -2{\mathrm{Re}}\int _{M} \left (\omega _{n}(\tau |d\varrho |^{2}) +\tau |d
\varrho |^{2}\mathrm{Tr}[\nabla \overline{\omega }^{n}]\right )
\langle \widetilde{\nabla}^{TX}_{\overline{\omega }_{n}}s,s\rangle _{h}
dv_{X}
\\
&-2{\mathrm{Re}}\int _{M} \tau |d\varrho |^{2}\langle \widetilde{\nabla}^{TX}_{
\widetilde{\nabla}^{TX}_{\omega _{n}}\overline{\omega }_{n}}s,s
\rangle _{h} dv_{X}.
\end{split}
\label{eq:3.60Nov}
\end{equation}
\end{lemma}
\begin{proof}
The formula follows from computations analogous to those in the proof of
{Lemma~\ref{lm:3.14Nov}}. Note that near $bM$, we have
\begin{equation}
L_{\omega _{n}} dv_{X}=\mathrm{Tr}[\nabla \overline{\omega }^{n}]dv_{X}.
\label{eq3.59}
\end{equation}
By replacing $dv_{X}$ by $\tau |d\varrho |^{2} dv_{X}$ and proceeding as
in the calculations of {Lemma~\ref{lm:3.14Nov}} we get {\eqref{eq:3.60Nov}}.
\end{proof}

\begin{proof}[Proof of {Proposition~\ref{prop:3.8}}]
This proposition follows directly from {\eqref{eq:3.51Novs}}, {Lemmas~\ref{lm:3.15} and \ref{lm:3.16}}. In particular, the two boundary terms
$-2{\mathrm{Re}}\int _{bM} \langle \widetilde{\nabla}^{TX}_{e^{(0,1)}_{
\mathfrak{n}}}s,s\rangle _{h} dv_{bM}$ cancel out.
\end{proof}

\section{Holomorphic Morse inequalities for \texorpdfstring{$q$}{q}-concave domains}
\label{sec_cvhmi}

In this section, we establish the holomorphic Morse inequalities for Levi
$q$-concave domains stated in {Theorem~\ref{thm:4.10important}}. Our method is
based on the Nakano-Griffiths inequality given in {Theorem~\ref{thm:NG-inequality}}. Moreover, the same strategy also applies to {Theorem~\ref{thm:4.13new}}. In parallel, {Theorem~\ref{thm:4.10important}} yields a
new proof of holomorphic Morse inequalities for $q$-concave manifolds,
as stated in {Theorem~\ref{TheoremB}}. Finally, in Subsection \ref{ss:4.5pqcorona}, we combine our method with the one used in
\cite[Section 3.5]{MM} to study the domains or manifolds with mixed convexity/concavity
such as $(p,q)$-coronas and $(p,q)$-convex-concave complex manifolds.

\subsection{Holomorphic Morse inequalities with modified metrics}
\label{sec4.1}

Let $(E,h^{E})$ and $(L,h^{L})$ be holomorphic Hermitian vector bundles
on a connected Hermitian manifold $(X,\Theta )$ with $\rank (L)=1$. Let
$\overline \partial ^{E}_{k}: \Omega _{0}^{0,\bullet}(X,L^{k}\otimes E)
\rightarrow \Omega _{0}^{0,\bullet +1}(X,L^{k}\otimes E)$ be the associated
$\overline{\partial }$-operator. After taking into account the
$L^{2}$ inner product associated with $\Theta $ and $h^{L}$, $h^{E}$, we
can define the maximal extension of $\overline \partial ^{E}_{k}$, its
Hilbert-space adjoint  $\overline \partial ^{E*}_{k}$, and then the Gaffney
extension $\square ^{E}_{k}$ on
$L^{2}_{0,\bullet}(X,L^{k}\otimes E)$.

Recall that $R^{L}$ denotes the Chern curvature form of $(L, h^{L})$. For
any subset $M\subset X$, we set
\begin{equation}
\begin{split} M(j,h^{L})=\{x\in M\;:\; \sqrt{-1}R^{L}_{x}
\text{ is non-degenerate} \qquad \qquad &
\\
\text{ and has exactly } j \text{ negative eigenvalues}\},&
\end{split}
\label{eq:KjhL}
\end{equation}
and $M(\leq j,h^{L})=\cup_{0\leq \ell \leq j} M(\ell ,h^{L})$.

First, we prove the following result as an intermediate step towards {Theorem~\ref{thm:4.10important}}.

\begin{thm}%
\label{TheoremA}
Let $M\Subset X$ be a smooth domain of a Hermitian manifold
$(X,\Theta )$ of dimension $n$ such that the Levi form
$\mathscr{L}_{bM}$ has at least $n-q$ negative eigenvalues ($1\leq q
\leq n-1$), that is, $M$ is a $q$-concave domain of $X$. Let
$(E, h^{E})$ and $(L,h^{L})$ be holomorphic Hermitian vector bundles on
$X$ with $\rank (L)=1$. Then there exists a smooth Hermitian metric
$h^{L}_{\chi}$ on $L$ (which is obtained by modifying $h^{L}$ near
$bM$) such that, given a sufficiently small neighborhood $U_{1}$ of
$bM$ with $M'=M\setminus \overline{U}_{1}$ having a smooth boundary, we
have the strong Morse inequalities for $j\leq n-q-1$ as
$k\rightarrow \infty $,
\begin{equation}
\sum _{\ell =0}^{j}(-1)^{j-\ell}\dim H^{0,\ell}_{(2)}(M,L^{k}\otimes E)
\leq \frac{\rank (E)k^{n}}{n!}\int _{M'(\leq \,j,h^{L}_{\chi})}(-1)^{j}c_{1}(L,h^{L}_{
\chi})^{n}+o(k^{n}).
\label{eq:4.3Jan25}
\end{equation}
For $j\leq n-q-1$, we have the weak Morse inequalities,
\begin{equation}
\dim H^{0,j}_{(2)}(M,L^{k}\otimes E)\leq \rank (E)\frac{k^{n}}{n!}
\int _{M'(j,h^{L}_{\chi})}(-1)^{j}c_{1}(L,h^{L}_{\chi})^{n}+o(k^{n}).
\label{eq:4.4Jan25}
\end{equation}
Moreover, if $q\leq n-2$, we also have
\begin{equation}
\dim H^{0}_{(2)}(M,L^{k}\otimes E)\geq \rank (E)\frac{k^{n}}{n!}\int _{M'(
\leq 1\,,h^{L}_{\chi})} c_{1}(L,h^{L})^{n}+o(k^{n}).
\label{eq:4.3Jan25-new}
\end{equation}
\end{thm}

A key step in the proof of {Theorem~\ref{TheoremA}} is the following proposition, whose
proof is deferred to Subsection \ref{ss4.2Jan25}.

\begin{prop}%
\label{prop:4.3Jan25}
Let $M\Subset X$ be a $q$-concave domain in a complex manifold $X$ of dimension
$n$, and let $(E, h^{E})$ and $(L,h^{L})$ be holomorphic Hermitian vector
bundles on $X$ with $\rank (L)=1$. Then there exists a defining function
$\varrho \in \mathscr{C}^{\infty}(X,\mathbb{R})$ of $M$ with
$M=X_{0}=\{x\in X\;:\; \varrho (x)<0\}$, a Hermitian metric
$\widetilde{\Theta}$ on $X$ and a modified Hermitian metric
$h^{L}_{\chi}$ on $L$ such that the following property holds: there exist
$C>0$, $C_{0}>0, C'>0$, and $\delta ,\, \varepsilon _{1}\in (0,1)$ such
that given any $c\in [-\delta ,\delta ]$ and any
$\varepsilon \in (0,\varepsilon _{1})$, for any
$s\in B^{0,j}(X_{c},L^{k}\otimes E)$ with
$\supp (s)\subset X_{\delta}\setminus \overline X_{-\delta}$ and
$j\leq n-q-1$, we have
\begin{equation}
(1+\varepsilon )\left (\|\overline \partial _{k}^{E} s\|^{2}_{L^{2}(X_{c})}+
\|\overline \partial _{k}^{E*} s\|^{2}_{L^{2}(X_{c})}\right )+ C(1+
\frac{1}{\varepsilon})\|s\|_{L^{2}(X_{c})}^{2}
\geq (C_{0} k-C')\|s\|_{L^{2}(X_{c})}^{2},
\label{eq4.5}
\end{equation}
where the $L^{2}$ inner product is induced by $\widetilde{\Theta}$ and
$h^{L}_{\chi}$, $h^{E}$, and
$X_{c}:=\{x\in X\;:\; \varrho (x)<c\}$.
\end{prop}

Using a cut-off function near the boundary $bM$ and the fact that
$B^{0,j}(M,L^{k}\otimes E)$ is dense in
$\Dom (\overline \partial ^{E}_{k})\cap \Dom (\overline \partial ^{E*}_{k})
\cap L^{2}_{0,j}(M,L^{k}\otimes E)$ with respect to the graph norm of
$\overline \partial ^{E}_{k}+\overline \partial _{k}^{E*}$, we obtain the
following \emph{optimal} fundamental estimates for $M$ (see also
\cite[(3.5.19)]{MM} and \cite[Proposition 3.4]{LSW23}).
\begin{prop}%
\label{prop:4.4new25}
Let $M\Subset X$ be a $q$-concave domain in a complex manifold $X$ of dimension
$n$, and let $(E, h^{E})$ and $(L,h^{L})$ be holomorphic Hermitian vector
bundles on $X$ with $\rank (L)=1$. Let $\varrho $,
$\widetilde{\Theta}$, $h^{L}_{\chi}$, $\delta $ be as in {Proposition~\ref{prop:4.3Jan25}}.

Then there exists $\widetilde{C}>0$ such that for any given
$\delta _{1} \in \,(0,\delta ]$ and for $k\geq 1$, we have
\begin{equation}
\label{eq-opfe}
\left (1-\frac{\widetilde{C}}{k}\right )\|s\|_{L^{2}(M)}^{2}\leq
\frac{\widetilde{C}}{k}\left (\|\overline \partial ^{E}_{k} s\|_{L^{2}(M)}^{2}+
\|\overline \partial ^{E*}_{k} s\|_{L^{2}(M)}^{2} \right )+\int _{
\overline X_{-\delta _{1}}} |s|^{2} dv_{X}
\end{equation}
for
$s\in \Dom (\overline \partial ^{E}_{k})\cap \Dom (
\overline \partial ^{E*}_{k})\cap L^{2}_{0,j}(M,L^{k}\otimes E)$ with
$j\leq n-q-1$, where the $L^{2}$ norm is induced by the metrics
$\widetilde{\Theta}$, $h^{L}_{\chi}$, $h^{E}$.
\end{prop}

Now we are ready to prove {Theorem~\ref{TheoremA}}.

\begin{proof}[Proof of {Theorem~\ref{TheoremA}}]
The main step is to apply {Theorem~\ref{thm-l2hmi}} (see
\cite[Theorem 3.2.13]{MM}). The optimal fundamental estimate {\eqref{eq-opfe}} holds for $(0,j)$-forms with any $j\leq n-q-1$. This allows
us to directly apply {Theorem~\ref{thm-l2hmi}} to the triplet
$\overline{X}_{-\delta _{1}}\subset M'\Subset M$ in our setting, where
$M'$ is a smooth domain. Explicitly, after working with the modified Hermitian
metrics $\widetilde{\Theta}$ on $X$ and $h^{L}_{\chi}$ on $L$ given in
{Proposition~\ref{prop:4.3Jan25}}, we have, for $j\leq n-q-1$, as
$k\rightarrow \infty $,
\begin{equation}
\begin{split} &\sum _{\ell =0}^{j}(-1)^{j-\ell}\dim H_{(2)}^{0,\ell}(M,L^{k}
\otimes E)\leq \rank (E)\frac{k^{n}}{n!}\int _{M'(\leq \,j,h^{L}_{
\chi})}(-1)^{j}c_{1}(L,h^{L}_{\chi})^{n}+o(k^{n}),
\\
&\dim H_{(2)}^{0,j}(M,L^{k}\otimes E)\leq \rank (E)\frac{k^{n}}{n!}
\int _{M'(j,h^{L}_{\chi})}(-1)^{j}c_{1}(L,h^{L}_{\chi})^{n}+o(k^{n}),
\end{split}
\label{eq:4.8Jan25}
\end{equation}
where each space $H_{(2)}^{0,\ell}(M,L^{k}\otimes E)$ is defined with respect
to $\widetilde{\Theta}$ and $h^{L}_{\chi}$, as in {Lemmas~\ref{lem_metric} and \ref{lem-question2}}, with a function $\chi $ satisfying
condition {\eqref{eq:4.36Jan}}.

By {Theorem~\ref{zqzq1}} and the fact that $M$ satisfies $Z(j)$ condition for
$j\leq n-q-1$, we conclude that {\eqref{eq:4.8Jan25}} holds for the groups
$H_{(2)}^{0,\bullet}(M,L^{k}\otimes E)$ that are defined with respect to
the original metrics $\Theta $ and $h^{L}$. Similarly,
{\eqref{eq:4.3Jan25-new}} follows from {Theorem~\ref{thm-l2hmi-2}}. In this
way, we complete the proof.
\end{proof}

\begin{rem}
\label{rem4.4}
In {Lemma~\ref{lem-question2}} and Subsection \ref{ss:4.3Jan25} we will explicitly
construct a metric $h^{L}_{\chi}$ by using a well-chosen function
$\chi $, in order to prove {Theorem~\ref{TheoremA}}.
\end{rem}

\subsection{Two technical lemmas and proof of {Proposition~\ref{prop:4.3Jan25}}}
\label{ss4.2Jan25}

We need two technical lemmas to prove {Proposition~\ref{prop:4.3Jan25}}.
In this subsection, we fix
$\varrho \in \mathscr{C}^{\infty}(X,\mathbb{R})$ such that
\begin{itemize}
\item There exists an interval $[u,v]$ ($u<v$) of $\mathbb{R}$ which consists
of regular values of $\varrho $ such that for $c\in [u,v]$,
$\varnothing \neq X_{c}:=\{x\in X\; :\; \varrho (x)<c\}\Subset X$.
\item For $c\in [u,v]$, at $x\in bX_{c}$, the Levi form
$\mathscr{L}_{\varrho}:=\partial \overline \partial \varrho |_{bX_{c}}$
has at least $n-q$ negative eigenvalues.
\end{itemize}

\begin{lemma}%
\label{lem_metric}
For any $C_{1}>0$, there exists a Hermitian metric $\Theta $ on $X$, with associated Riemannian metric $g^{TX}$, and
$\varepsilon _{1}\in (0,1)$ such that, for every holomorphic Hermitian
vector bundle $(F,h^{F})$ over $X$ and every $c\in (u,v)$,
the following statements holds.

(i) At each point of $X_{v}\setminus \overline X_{u}$,
\begin{equation}
|d\varrho |_{g^{TX}}=1.
\label{eq:4.8april2025}
\end{equation}
Thus the function $\varrho _{c}:=\varrho -c$ is a defining function
for $X_{c}$; let $e_{\mathfrak{n}}$ denote the metric dual of
$-d\varrho $.

(ii) For $j\in \{0,1,\ldots , n-q-1\}$ and
$0<\varepsilon \leq \varepsilon _{1}$, if
$s\in \Omega _{0}^{0,j}(X_{v}\setminus \overline X_{u},F)$ with
\begin{equation*}
\iota _{e_{\mathfrak{n}}}s|_{bX_{c}}\equiv 0,
\end{equation*}
then at each point $x\in bX_{c}$,
\begin{equation}
\langle [\sqrt{-1}\mathscr{L}_{\varrho _{c}},\Lambda ]^{bX_{c}}_{
\varepsilon }s,s \rangle _{h}(x)\geq C_{1}|s(x)|_{h}^{2}.
\label{eq:4.5dec24}
\end{equation}
\end{lemma}
\begin{proof}
Since any $c\in [u,v]$ is a regular value of $\varrho $, the smooth
$1$-form $d\varrho $ does not vanish on a neighborhood of
$X_{v}\setminus \overline X_{u}$. Let $HX:=\ker d\varrho $ be a (real)
subbundle of $TX$ in the neighborhood of
$X_{v}\setminus \overline X_{u}$. Put
\begin{equation}
H^{(1,0)}X:= T^{(1,0)}X\cap \mathbb{C}HX.
\label{eq:4.11HBX}
\end{equation}
Note that we have $H^{(1,0)}X|_{bX_{c}}=T^{(1,0)}bX_{c}$, the analytic
tangent space to $bX_{c}$ defined in Subsection \ref{ss:2.1}. Set
$H_{b}X:=\{w+\overline{w}\in HX\;:\; w\in H^{(1,0)}X\}$; this is a real subbundle
of $HX$ with corank $1$.

Let $e_{\mathfrak{n}}$ be a smooth vector field in the neighborhood of
$X_{v}\setminus \overline X_{u}$ such that
$e_{\mathfrak{n}}(\varrho )\equiv -1$. Then we have the splitting of the
tangent bundle in the neighborhood of
$X_{v}\setminus \overline X_{u}$
\begin{equation}
TX=H X \oplus \mathbb{R}e_{\mathfrak{n}}=(H_{b} X\oplus \mathbb{R}Je_{
\mathfrak{n}}) \oplus \mathbb{R}e_{\mathfrak{n}}.
\label{eq:4.7dec24}
\end{equation}
Note that the complex structure $J$ preserves $H_{b} X$.

Now we fix any Hermitian metric $\Theta _{0}$ (and the corresponding
$J$-invariant Riemannian metric $g^{TX}_{0}$) such that the splitting {\eqref{eq:4.7dec24}} is orthogonal and $e_{\mathfrak{n}}$ has norm one on
the neighborhood of $X_{v}\setminus \overline X_{u}$. Moreover, we can
write
\begin{equation}
g^{TX}_{0}=g^{HX}_{0}\oplus (d\varrho \otimes d\varrho ),
\label{eq4.12}
\end{equation}
where $g^{HX}_{0}$ is a smooth Euclidean metric on $HX$. Hence
$d\varrho _{c}=d\varrho $ has norm $1$ with respect to $g^{TX}_{0}$. In
particular, let $\Theta ^{b}_{0}$ be the induced Hermitian metric on
$H_{b}X$.

With respect to the above $g^{TX}_{0}$ and proceeding as in the proof of
{Lemma~\ref{lm:7.3}}, we fix a local orthonormal frame
$\{\omega '_{j}\}_{j=1}^{n}$ of $(T^{(1,0)}X,h^{T^{(1,0)}X})$ (and the
dual frame $\{\omega ^{\prime ,j}\}_{j=1}^{n}$) on any small domain
$U$ inside the neighborhood of $X_{v}\setminus \overline X_{u}$, such
that
\begin{equation}
\omega '_{n}=\sqrt{2} e^{(1,0)}_{\mathfrak{n}},\quad \omega ^{\prime ,n}=-
\sqrt{2}\partial \varrho ,
\label{eq4.13}
\end{equation}
and moreover $\{\omega '_{j}\}_{j=1}^{n-1}$ forms an orthonormal frame
of $(H^{(1,0)}X, g^{T^{(1,0)}X}_{0})$. Note that for $c\in (u,v)$,
$H^{(1,0)}X|_{bX_{c}}$ is exactly the analytic tangent space to
$bX_{c}$.

In the sequel, we will work on $U\cap M$ with $M=X_{c}$ for some
$c\in (u,v)$, and set
\begin{equation*}
\varrho _{c}:=\varrho -c.
\end{equation*}
The Levi form on $bM$ is given by
\begin{equation*}
\mathscr{L}_{\varrho _{c}}=\partial \overline \partial \varrho |_{bM}=
\sum _{k=1}^{n-1}\partial \overline \partial \varrho (\omega '_{k},
\overline {\omega }'_{k})\omega ^{\prime ,k}\wedge
\overline {\omega }^{\prime ,k},
\end{equation*}
and the trace with respect to $g_{0}^{TX}$ is given by
\begin{equation*}
\tr _{bM}[\sqrt{-1}\mathscr{L}_{\varrho _{c}}]=\sum _{k=1}^{n-1}\partial
\overline \partial \varrho (\omega '_{k},\overline {\omega }'_{k}).
\end{equation*}
Note that
$\Lambda |_{bM}=-\sqrt{-1}\sum _{k=1}^{n-1} \iota _{
\overline {\omega }'_{k}}\iota _{\omega '_{k}}$.

By {Proposition~\ref{prop:2.3dec24}}, $\mathscr{L}_{\varrho _{c}}$ has at
least $n-q$ negative eigenvalues on $bM$, that is, after a unitary change
on the frame $\{\omega '_{j}\}_{j=1}^{n-1}$ of
$T^{(1,0)}X\cap \mathbb{C}HX$, we can write
\begin{equation}
\mathscr{L}_{\varrho _{c}}=\sum _{j=1}^{n-1} \lambda '_{j} \omega ^{
\prime ,j}\wedge \overline {\omega }^{\prime ,j}, \quad
\text{ where } \lambda '_{1}\leq \cdots \leq \lambda '_{n-q}<0, \quad
\lambda '_{1}\leq \cdots \leq \lambda '_{n-1}.
\label{eq:4.10dec24}
\end{equation}

Fix an integer $j\in \{0,1,\ldots , n-q-1$\}. Take
$s\in \Omega _{0}^{0,j}(X_{v}\setminus \overline X_{u},F)$ with
$\iota _{e_{\mathfrak{n}}}s\equiv 0$. Then there is no
$\overline{\omega }^{\prime ,n}$-component in $s$ with respect to the local
frame $\{\overline{\omega }^{\prime ,i}\}_{i=1}^{n}$. For a subset of indices
$I\subset \{1,\ldots , n-1\}$ with $|I|=j$, we set
$\complement I:=\{1,\cdots ,n-1\}\setminus I$, and clearly we have
$|\complement I|=n-j-1\geq q$.

Now we write
\begin{equation}
s=\sum _{|I|=j,\, 1\leq \ell \leq \rank (F)} s^{\prime ,\ell}_{I}
\overline \omega ^{\prime ,I}\otimes f_{\ell }%
\label{eq4.15}
\end{equation}
on the small domain $U$, where
$s^{\prime ,\ell}_{I}=s^{\prime ,\ell}_{I}(x)\in \mathscr{C}^{\infty}(U)$
and $\{f_{\ell}: 1 \leq \ell \leq \rank (F)\}$ is a local orthonormal frame
of $(F,h^{F})$ on $U$. For $\varepsilon >0$,
\begin{equation}
\begin{split} \langle [\sqrt{-1}\mathscr{L}_{\varrho _{c}},\Lambda ]^{bM}s,s
\rangle _{h_{0}}+ \varepsilon \tr _{bM}[\sqrt{-1}\mathscr{L}_{
\varrho _{c}}]|s|_{h_{0}}^{2} =\sum _{I,\, \ell}\left (\sum _{k\in I}
\lambda '_{k}+ (\varepsilon -1)\sum _{k=1}^{n-1} \lambda '_{k}\right )|s^{
\prime ,\ell}_{I}|^{2}.
\end{split}
\label{eq:4.12dec24}
\end{equation}
Note that
\begin{equation}
\begin{split} \sum _{k\in I} \lambda '_{k}+(\varepsilon -1)\sum _{k=1}^{n-1}
\lambda '_{k}&= \sum _{k\in {\complement I}}(-\lambda '_{k})+
\varepsilon \sum _{k=1}^{n-1} \lambda '_{k}
\\
&\geq (-\lambda '_{n-q})+(|\complement I|-1)(-\lambda '_{n-1})+
\varepsilon (n-1) \lambda '_{1}.
\end{split}
\label{eq:4.13dec24}
\end{equation}

Let $C>0$ be sufficiently large such that $C-\frac{n-1}{C}>C_{1}$. Note that
we have the Hermitian metric on $H_{b}X$ induced from $\Theta _{0}$ (see {\eqref{eq:4.7dec24}})
\begin{equation}
\Theta ^{b}_{0}=\sqrt{-1}\sum _{k=1}^{n-1} \omega ^{\prime ,k}
\wedge \overline{\omega }^{\prime ,k}.
\label{eq4.18}
\end{equation}
Then we modify the metric $g^{HX}_{0}$ to be a new metric $g^{HX}$ on
$U$; if $\Theta ^{b}$ is the Hermitian metric on $H_{b}X$ induced by
$g^{HX}$, then
\begin{equation}
\Theta ^{b}=\sum _{k=1}^{n-1} \kappa _{k} \sqrt{-1} \omega ^{\prime ,k}
\wedge \overline{\omega }^{\prime ,k},
\label{eq:4.14dec24}
\end{equation}
where the frame $\{\omega ^{\prime ,k}\}_{k=1}^{n-1}$ is orthonormal with
respect to the original $\Theta _{0}^{b}$ and satisfies {\eqref{eq:4.10dec24}}, and each $\kappa _{k}>0$ is some constant to be fixed.

Now we define the new metric $g^{TX}$ associated to the Hermitian metric
\begin{equation}
\Theta =\sum _{k=1}^{n-1} \kappa _{k} \sqrt{-1} \omega ^{\prime ,k}
\wedge \overline{\omega }^{\prime ,k} + 2 \sqrt{-1} \partial \varrho
\wedge \overline{\partial }\varrho .
\label{eq:4.15dec24}
\end{equation}
With respect to the metric $g^{TX}$, the Levi form
$\mathscr{L}_{\varrho _{c}}$ on $bM$ has the eigenvalues
$\{\lambda _{k}\}_{k=1}^{n-1}$ such that
\begin{equation}
\lambda _{k}=\frac{\lambda '_{k}}{\kappa _{k}}.
\label{eq4.21}
\end{equation}
Subsequently, we select each $\kappa _{k} > 0$ appropriately to obtain
\begin{equation}
2C >|\lambda _{1}|\geq \cdots \geq |\lambda _{n-q}|>C, \qquad |
\lambda _{k}|<1/C,\; \text{ for } n-q+1\leq k\leq n-1.
\label{eq:4.17dec24}
\end{equation}
Next we choose $\varepsilon _{1}>0$ such that
$\varepsilon _{1}<\min \{1,\frac{1}{2(n-1)C^{2}}\}$. By {\eqref{eq:4.12dec24}} and {\eqref{eq:4.13dec24}}, with respect to
$g^{TX}$ and for $\varepsilon \in (0,\varepsilon _{1})$, we obtain at any
point $x\in U\cap bM$,
\begin{equation}
\begin{split} &\langle [\sqrt{-1}\mathscr{L}_{\varrho _{c}},\Lambda ]^{bM}s,s
\rangle _{h}(x) +\varepsilon \tr _{bM}[\sqrt{-1}\mathscr{L}_{\varrho _{c}}]|s(x)|_{h}^{2}
\\
&\quad \geq \left ((-\lambda _{n-q})+(n-j-2)(-\lambda _{n-1})+
\varepsilon (n-1) \lambda _{1} \right ) |s(x)|_{h}^{2}
\\
&\quad \geq ( C-\frac{n-j-1}{C})|s(x)|_{h}^{2}\geq C_{1}|s(x)|_{h}^{2}.
\end{split}
\label{eq:4.19dec24}
\end{equation}
Note that $e_{\mathfrak{n}}$ is still the dual vector of
$-d\varrho $ with respect to this new metric $g^{TX}$.

Since $\mathscr{L}_{\varrho _{c}}$ depends smoothly on $c\in [u,v]$, the ordered eigenvalues of $\mathscr{L}_{\varrho _{c}}$ vary continuously with
$c$. Then we can always obtain a local Hermitian metric
$\Theta $ by modifying $\Theta _{0}$ locally on $U$ such that {\eqref{eq:4.19dec24}} holds for any $c\in (u,v)$ with $bM=bX_{c}$. Patching
the local metrics as constructed above with the help of a partition of
unity, we conclude our lemma.
\end{proof}

From now on, we consider $(X,\Theta )$ as a Hermitian manifold, where
$\Theta $ will be taken as constructed in {Lemma~\ref{lem_metric}}. Let
$dv_{X}=\Theta ^{n}/n!$ be the volume form on $X$ with the induced metric
$dv_{bX_{c}}$ on the boundary $bX_{c}$.

\begin{lemma}%
\label{lem-question2}
Let $L$ be a holomorphic line bundle on $X$. For any given constant
$C_{0}>0$, there exists a Hermitian metric $\Theta $ (as well as the associated
Riemannian metric $g^{TX}$), a Hermitian metric $h_{\chi}^{L}$ on
$L$, and $\varepsilon _{1}\in (0,1)$ such that $\Theta $ and
$\varepsilon _{1}$ satisfy all the conditions of {Lemma~\ref{lem_metric}} and, for any
$\varepsilon \in (0,\varepsilon _{1})$, $\tau \in [0,1]$,
$s\in \Omega _{0}^{0,j}(X_{v}\setminus \overline X_{u},L^{k})$ with
$j\leq n-q-1$ and any $k\in \mathbb{N}_{\geq 1}$, we have, at any point
of $X_{v}\setminus \overline X_{u}$,
\begin{equation}
\label{eq-q2}
\,\left \langle [\sqrt{-1}R^{L}_{\chi},\Lambda ]_{\varepsilon} s,s
\right \rangle _{h} +2(1-\varepsilon )\tau \left \langle R^{L}_{\chi}(e^{(1,0)}_{
\mathfrak{n}}, e^{(0,1)}_{\mathfrak{n}})s, s\right \rangle _{h}\geq C_{0}|s|_{h}^{2}
\end{equation}
holds, where $e_{\mathfrak{n}}$ is the metric dual of $-d\varrho $ that has norm
$1$ on $X_{v}\setminus \overline X_{u}$, and
$R^{L}_{\chi}:=R^{(L,h^{L}_{\chi})}$ is the Chern curvature form of
$(L,h^{L}_{\chi})$.
\end{lemma}
\begin{proof}
We take the geometric data $(X,J,\Theta , g^{TX}, \varepsilon _{1})$ as
constructed in {Lemma~\ref{lem_metric}} with respect to the given constant
$C_{1}>0$ satisfying {\eqref{eq:4.5dec24}}. In particular, we have
$|d\varrho |_{g^{T^{\ast }X}}\equiv 1$ on a neighborhood of
$X_{v}\setminus \overline X_{u}$. Recall that $e_{\mathfrak{n}}$ is the
dual vector field to $-d\varrho $ with respect to $g^{TX}$.

Regarding the decomposition
$T^{(1,0)}X=H^{(1,0)}X\oplus \mathbb{C}\omega _{n}$, following {\eqref{eq:4.14dec24}} and {\eqref{eq:4.15dec24}} after appropriately adjusting
the constants $\kappa_k$, we obtain
\begin{equation}
\Theta =\Theta ^{b} + \sqrt{-1} \omega ^{n}\wedge \overline{\omega }^{n},
\label{eq4.25}
\end{equation}
where $\omega ^{n}=-\sqrt{2}\partial \varrho $. Let
$\{\omega_{j}\}_{j=1}^{n-1}$ be a local orthonormal frame of $H^{(1,0)}X$ with
respect to the Hermitian metric $\Theta ^{b}$ such that we can write
\begin{equation}
\begin{split} &\Theta ^{b}=\sqrt{-1}\sum _{j=1}^{n-1}\omega ^{j}
\wedge \overline{\omega }^{j},
\\
&\mathscr{L}_{\varrho}|_{H^{(1,0)}X\times H^{(0,1)}X}=\sum _{j=1}^{n-1}
\lambda _{j}\omega ^{j}\wedge \overline{\omega }^{j},
\end{split}
\label{eq:4.21dec24}
\end{equation}
where
$\lambda _{1}\leq \ldots \leq \lambda _{n-q} \leq \ldots \leq
\lambda _{n-1}$ satisfy {\eqref{eq:4.17dec24}} and $\lambda _{n-q}<0$.

Let $h^{L}$ be an arbitrary smooth Hermitian metric on $L$. For
$T>0$, we introduce a new family of Hermitian metrics on
$X_{v}\setminus \overline X_{u}$
\begin{equation}
\Theta _{T}= \frac{1}{T^{2}}\Theta ^{b} + \sqrt{-1} \omega ^{n}
\wedge \overline{\omega }^{n}.
\label{eq:4.20dec24}
\end{equation}
Let $g^{TX}_{T}$ be the Riemannian metric associated with
$\Theta _{T}$, and let $\langle \cdot ,\cdot \rangle _{h,T}$,
$|\cdot |_{h,T}$ denote the (pointwise) Hermitian inner product or norm
with respect to $\Theta _{T}$ and $h^{L}$. Note that for any $T>0$, we
always have $|e_{\mathfrak{n}}|_{g^{TX}_{T}}\equiv 1$ on
$X_{v}\setminus \overline X_{u}$. For simplicity, we set
$\mathscr{L}_{\varrho}:=\partial \overline{\partial }\varrho |_{
\Lambda ^{2}\mathbb{C}H_{b}X}$.

Then we have, by our proof of {Lemma~\ref{lem_metric}}, that if
$s\in \Omega _{0}^{0,j}(X_{v}\setminus \overline X_{u},L^k)$ with
\begin{equation}
\iota _{e_{\mathfrak{n}}}s\equiv 0,
\label{eq:4.23dec}
\end{equation}
then for $T>0$, $c\in (u,v)$, $0<\varepsilon <\varepsilon _{1}$, and at
each point $x\in bX_{c}$,
\begin{equation}
\langle [\sqrt{-1}\mathscr{L}_{\varrho},\Lambda _{T}]^{bX_{c}}_{
\varepsilon }s,s \rangle _{h,T}(x)\geq T^{2} C_{1}|s(x)|_{h,T}^{2},
\label{eq:4.24dec24}
\end{equation}
where $\Lambda _{T}=\iota _{\Theta _{T}}$ is the adjoint of the Lefschetz
operator $\mathrm{Lef}_{\Theta _{T}}$ with respect to
$\langle \cdot ,\cdot \rangle _{h,T}$.

First, we need to prove the following claim: For any given
$C_{2}>0$, there exists a constant $T_{0}>0$ such that for
$s\in \Omega _{0}^{0,j}(X_{v}\setminus \overline X_{u},L^{k})$ (with
$0\leq j \leq n-q-1$ and $k\in \mathbb{N}$),
$0<\varepsilon <\varepsilon _{1}$, $\tau \in [0,1]$, and at
$x\in X_{v}\setminus \overline X_{u}$, we have
\begin{equation}
\left \langle [\sqrt{-1}\partial \overline{\partial }\varrho ,
\Lambda _{T_{0}}]_{\varepsilon} s,s\right \rangle _{h,T_{0}}+2(1-
\varepsilon )\tau \left \langle \partial \overline{\partial }\varrho (e^{(1,0)}_{
\mathfrak{n}}, e^{(0,1)}_{\mathfrak{n}})s,s\right \rangle _{h,T_{0}}
\geq C_{2}|s|_{h,T_{0}}^{2}.
\label{eq-25dec}
\end{equation}
Note that the power $k$ for the line bundle $L^{k}$ does not affect the
estimate since all the involved operators are zeroth order linear operators
acting only on the differential forms part.

Indeed, after comparing {\eqref{eq:4.24dec24}} and {\eqref{eq-25dec}}, we have
to deal with the $\overline{\omega }^{n}$-component of the section
$s$. For the above orthonormal frame
$\{\omega _{j}\}_{j=1}^{n}$, set 
$\omega _{j,T}=T\omega _{j}$ for $j=1,\ldots , n-1$ and
$\omega _{n,T}=\omega _{n}$. Then $\{\omega _{j,T}\}_{j=1}^{n}$ is a local orthonormal frame of
$(T^{(1,0)}X, h^{T^{(1,0)}X}_{T})$ (where the Hermitian inner product is
induced by $\Theta _{T}$); denote its dual frame by $\{\omega ^{j}_{T}\}_{j=1}^{n}$. By {\eqref{eq:4.21dec24}}, we
have
\begin{equation}
\begin{split} \partial \overline \partial \varrho = &T^{2}\sum _{j=1}^{n-1}
\lambda _{j}\omega _{T}^{j}\wedge \overline {\omega }_{T}^{j}+\nu _{n}
\omega ^{n}\wedge \overline {\omega }^{n}
\\
&+T\sum _{j=1}^{n-1}\nu _{j} \omega ^{j}_{T}\wedge
\overline {\omega }^{n}+ T\sum _{j=1}^{n-1}\overline \nu _{j} \omega ^{n}
\wedge \overline {\omega }^{j}_{T},
\end{split}
\label{eq4.31}
\end{equation}
where
$\nu _{j}=\partial \overline \partial \varrho (\omega _{j},
\overline{\omega }_{n})$ for $1\leq j\leq n$.

For any
$s\in \Omega _{0}^{0,j}(X_{v}\setminus \overline X_{u},L^{k})$ ($0
\leq j \leq n-q-1$), we have an orthogonal decomposition
\begin{equation}
s=s_{0}+\overline{\omega }^{n}\wedge s_{1}
\label{eq4.32}
\end{equation}
such that
\begin{equation}
\begin{split} & s_{0}\in \Omega _{0}^{0,j}(X_{v}\setminus \overline X_{u},L^{k}),
\; \iota _{e_{\mathfrak{n}}} s_{0}\equiv 0,
\\
& s_{1}\in \Omega _{0}^{0,j-1}(X_{v}\setminus \overline X_{u},L^{k}),
\; \iota _{e_{\mathfrak{n}}} s_{1}\equiv 0.
\end{split}
\label{eq:4.28Jan25}
\end{equation}
Then an elementary computation shows that%
\begin{align}
 &[\sqrt{-1}\partial \overline{\partial }\varrho ,
\Lambda _{T}]_{\varepsilon} s +(1-\varepsilon )\partial
\overline{\partial }\varrho (\omega _{n}, \overline{\omega }_{n})s
\notag \\
&=[\sqrt{-1}\mathscr{L}_{\varrho},\Lambda _{T}]^{bM}_{\varepsilon} s_{0}+
T\sum _{k=1}^{n-1} \overline{\nu }_{k} \overline{\omega }^{k}_{T}
\wedge s_{1}
\notag \\
&\quad + \overline{\omega }^{n}\wedge \left ([\sqrt{-1}\mathscr{L}_{
\varrho},\Lambda _{T}]^{bM}_{\varepsilon} s_{1}+T\sum _{k=1}^{n-1}
\nu _{k} \iota _{\overline{\omega }_{k,T}} s_{0} +\partial
\overline{\partial }\varrho (\omega _{n},\overline{\omega }_{n})s_{1}
\right ).
\label{eq:4.29dec24}
\end{align}

Note that {\eqref{eq:4.29dec24}} holds for any $s$ as above and for any
$T>0$ and $0<\varepsilon <\varepsilon _{1}$. Then by {\eqref{eq:4.24dec24}}, {\eqref{eq:4.28Jan25}} and the Cauchy-Schwarz inequality,
we obtain that there exists $C_{3}>0$ such that for
$s\in \Omega _{0}^{0,j}(X_{v}\setminus \overline X_{u},L^{k})$ ($0
\leq j\leq n-q-1$ and $k\in \mathbb{N}$), $T>0$,
$0<\varepsilon <\varepsilon _{1}$, $\tau \in [0,1]$, and at
$x\in X_{v}\setminus \overline X_{u}$
\begin{align}
& \left \langle [\sqrt{-1}\partial \overline{\partial }\varrho ,
\Lambda _{T}]_{\varepsilon} s,s\right \rangle _{h,T}+2(1-\varepsilon )
\tau \left \langle \partial \overline{\partial }\varrho (e^{(1,0)}_{
\mathfrak{n}}, e^{(0,1)}_{\mathfrak{n}})s,s\right \rangle _{h,T}\notag\\
&\quad \geq (T^{2}
C_{1} - TC_{3} -C_{3})|s|_{h,T}^{2}.
\label{eq-30Jan25}
\end{align}
Therefore, after fixing a sufficiently large $T_{0}>0$ such that
$T_{0}^{2} C_{1} - T_{0}C_{3} -C_{3}\geq C_{2}$, we conclude {\eqref{eq-25dec}}.

The next step is to find a new Hermitian metric $h^{L}_{\chi}$ as required
in our lemma with the help of {\eqref{eq-25dec}}. From now on, we fix
$C_{2}>0$ and $T_{0}\gg 1$ as in {\eqref{eq-25dec}}.

Let $\chi =\chi (t)$ be a real smooth function on $\mathbb{R}$, to be specified
later. Let $h^{L}_{\chi}:=h^{L}e^{-\chi (\varrho )}$ be the modified Hermitian
metric on $L$. The Chern curvature of $(L,h^{L}_{\chi})$ is given by
\begin{equation}
R^{L}_{\chi}=R^{L}+\chi '(\rho )\partial \overline \partial \varrho +
\chi ''(\varrho )\partial \varrho \wedge \overline \partial \varrho .
\label{eq4.36}
\end{equation}
Then the term
$\langle [\sqrt{-1}R^{L}_{\chi},\Lambda _{T_{0}}]_{\varepsilon }s,s
\rangle _{h,T_{0}}$ equals
\begin{equation}
\begin{split} \langle [\sqrt{-1}R^{L},\Lambda _{T_{0}}]_{\varepsilon }s,s
\rangle _{h,T_{0}}+\chi '(\varrho )\langle [\sqrt{-1}\partial
\overline \partial \varrho ,\Lambda _{T_{0}}]_{\varepsilon }s,s
\rangle _{h,T_{0}}&
\\
+\chi ''(\varrho )\langle [\sqrt{-1}\partial \varrho \wedge
\overline \partial \varrho ,\Lambda _{T_{0}}]_{\varepsilon }s,s
\rangle _{h,T_{0}}&.
\end{split}
\label{eq:4.32Jan}
\end{equation}
Similarly, note that
$\partial \varrho \wedge \overline{\partial }\varrho (e^{(1,0)}_{
\mathfrak{n}}, e^{(0,1)}_{\mathfrak{n}})=\frac{1}{4}$, we have
\begin{equation}
\begin{split} &\left \langle R^{L}_{\chi}(e^{(1,0)}_{\mathfrak{n}}, e^{(0,1)}_{
\mathfrak{n}})s,s\right \rangle _{h,T_{0}}
\\
=&\left \langle R^{L}(e^{(1,0)}_{\mathfrak{n}}, e^{(0,1)}_{
\mathfrak{n}})s,s\right \rangle _{h,T_{0}}+\chi '(\varrho )\left
\langle \partial \overline{\partial }\varrho (e^{(1,0)}_{\mathfrak{n}},
e^{(0,1)}_{\mathfrak{n}})s,s\right \rangle _{h,T_{0}}+\frac{1}{4}
\chi ^{\prime }(\varrho )|s|^{2}_{h,T_{0}}.
\end{split}
\label{eq:4.33Jan}
\end{equation}

There exists a constant $C_{4}>0$ (depending on $T_{0}$,
$\Theta $ and $h^{L}$) such that for any
$s\in \Omega _{0}^{0,j}(X_{v}\setminus \overline X_{u},L^{k})$ (with any
$0\leq j\leq n-q-1$ and $k\in \mathbb{N}$),
$0<\varepsilon <\varepsilon _{1}$, $\tau \in [0,1]$, and at
$x\in X_{v}\setminus \overline X_{u}$, we have
\begin{equation}
\begin{split} &\left |\left \langle [\sqrt{-1}R^{L},\Lambda _{T_{0}}]_{
\varepsilon} s,s\right \rangle _{h,T_{0}}+2(1-\varepsilon )\tau
\left \langle R^{L}(e^{(1,0)}_{\mathfrak{n}}, e^{(0,1)}_{\mathfrak{n}})s,s
\right \rangle _{h,T_{0}}\right |\leq C_{4}|s|_{h,T_{0}}^{2},
\\
&\left |\left \langle [\sqrt{-1}\partial \varrho \wedge
\overline{\partial }\varrho ,\Lambda _{T_{0}}]_{\varepsilon} s,s
\right \rangle _{h,T_{0}}+\frac{(1-\varepsilon )\tau}{2}|s|^{2}_{h,T_{0}}
\right |\leq C_{4}|s|_{h,T_{0}}^{2}.
\end{split}
\label{eq:4.34Jan}
\end{equation}

Combining {\eqref{eq-25dec}}, {\eqref{eq:4.32Jan}}--
{\eqref{eq:4.34Jan}}, we obtain
\begin{equation}
\begin{split} &\left \langle [\sqrt{-1}R^{L}_{\chi},\Lambda ]_{
\varepsilon} s,s\right \rangle _{h,T_{0}}+2(1-\varepsilon )\tau
\left \langle R^{L}_{\chi}(e^{(1,0)}_{\mathfrak{n}}, e^{(0,1)}_{
\mathfrak{n}})s,s\right \rangle _{h,T_{0}}
\\
&\quad \geq \left ( C_{2} \chi '(\varrho )-C_{4} \left (1+|\chi ^{
\prime }(\varrho )|\right )\right ) |s|^{2}_{h,T_{0}}.
\end{split}
\label{eq:4.35Jan}
\end{equation}

Then we choose the function $\chi $ such that for $t\in [u,v]$, we have
\begin{equation}
C_{2} \chi '(t)-C_{4} (1+|\chi ^{\prime }(t)|)\geq C_{0}.
\label{eq:4.36Jan}
\end{equation}
Such a function always exists. Finally, since
$h^{L}_{\chi}:=h^{L}e^{-\chi (\varrho )}$ and
$h^{L^{k}}_{\chi}:=h^{L^{k}}e^{-k\chi (\varrho )}$, we obtain from {\eqref{eq:4.35Jan}} and {\eqref{eq:4.36Jan}} that at each point
$x\in X_{v}\setminus \overline X_{u} $
\begin{equation}
\left \langle [\sqrt{-1}R^{L}_{\chi},\Lambda ]_{\varepsilon} s,s
\right \rangle _{h^{L}_{\chi},T_{0}}+2(1-\varepsilon )\tau \left
\langle R^{L}_{\chi}(e^{(1,0)}_{\mathfrak{n}}, e^{(0,1)}_{
\mathfrak{n}})s,s\right \rangle _{h^{L}_{\chi},T_{0}} \geq C_{0} |s|^{2}_{h^{L}_{
\chi},T_{0}}.
\label{eq:4.37Jan}
\end{equation}
Consequently, the Hermitian metrics $\Theta _{T_{0}}$ and
$h^{L}_{\chi}$ satisfy the requirements of the lemma, thus completing the
proof.
\end{proof}
\begin{proof}[Proof of {Proposition~\ref{prop:4.3Jan25}}]
As explained in the text before {Definition~\ref{def:2.6domain}}, since the
domain $M$ is assumed to be Levi $q$-concave, we can take a defining function
$\varrho \in \mathscr{C}^{\infty}(X,\mathbb{R})$ for $M$ with
$M=\{x\in X\;:\; \varrho (x)<0\}$, and there exists a sufficiently small
$\delta >0$ such that for all $c\in [-\delta ,\delta ]$, $c$ is a regular
value of $\varrho $, and $\mathscr{L}_{bX_{c}}$ has at least $n-q$ negative
eigenvalues on $bX_{c}$. Moreover,
$\sqrt{-1}\partial \overline{\partial }\varrho $ has at least
$n-q+1$ negative eigenvalues at each point
$x\in X_{\delta}\setminus X_{-\delta}$.

Apply {Lemmas~\ref{lem_metric} and \ref{lem-question2}} to the function
$\varrho $ with $[u,v]=[-\delta ,\delta ]$. Fix the metrics
$\widetilde{\Theta}$, $h^{L}_{\chi}$ and constants
$\varepsilon _{1}$, $C_{1}$, $C_{0}$ so that {\eqref{eq:4.8april2025}}, {\eqref{eq:4.5dec24}} and {\eqref{eq-q2}} hold.

For any $c\in [-\delta ,\delta ]$, apply the Nakano-Griffiths formula
from {Theorem~\ref{thm:NG-inequality}} to $X_{c}$ with the defining function
$\varrho _{c}:=\varrho - c$. For this purpose, let
$\tau \in \mathscr{C}^{\infty}(X,[0,1])$ be such that
$\tau \equiv 1$ near $X_{\delta}\setminus \overline{X}_{-\delta}$ and
$\supp \tau $ is included in a small open neighborhood of
$X_{\delta}\setminus \overline{X}_{-\delta}$ where $d\varrho $ does not
vanish.

By the proof of {\eqref{eq:3.44Nov}}, when applying it to $\varrho _{c}$ and $X_c$, all the involved constants can be chosen uniformly in $c\in [-\delta ,\delta ]$. Therefore, we conclude that there exists a constant
$C=C(X,\widetilde{\Theta},\delta,\varrho ,\tau )>0$ such that for $c\in [-\delta ,\delta ]$, and for any
$\varepsilon \in (0,\varepsilon _{1}]$,
$s\in \Omega _{0}^{0,j}(X_{\delta}\setminus \overline X_{-\delta},L^{k}
\otimes E)$ with $\iota _{e_{\mathfrak{n}}}s=0$ on $bX_{c}$ for
$0\leq j\leq n-q-1$, we have
\begin{equation}
\begin{split} &(1+\varepsilon )\left (\|\overline \partial ^{E}_{k} s
\|^{2}_{L^{2}(X_{c})}+\|\overline \partial ^{E*}_{k} s\|^{2}_{L^{2}(X_{c})}
\right )+C(1+\frac{1}{\varepsilon})\|s\|^{2}_{L^{2}(X_{c})}
\\
&\quad \geq \, k \left (\left \langle [\sqrt{-1}R^{L}_{\chi},\Lambda ]_{
\varepsilon} s,s\right \rangle _{L^{2}(X_{c})}+2(1-\varepsilon )
\left \langle \tau R^{L}_{\chi}(e^{(1,0)}_{\mathfrak{n}}, e^{(0,1)}_{
\mathfrak{n}})s,s\right \rangle _{L^{2}(X_{c})}\right )
\\
&\qquad \qquad +\left \langle [\sqrt{-1}R^{E},\Lambda ]_{\varepsilon} s,s
\right \rangle _{L^{2}(X_{c})}+2(1-\varepsilon )\left \langle \tau R^{E}(e^{(1,0)}_{
\mathfrak{n}}, e^{(0,1)}_{\mathfrak{n}})s,s\right \rangle _{L^{2}(X_{c})}
\\
&\qquad \qquad +\int _{bX_{c}}\langle [\sqrt{-1}\mathscr{L}_{\varrho},
\Lambda ]^{bX_{c}}_{\varepsilon}s,s\rangle _{h} dv_{bX_{c}}+(1-
\varepsilon )\int _{bX_{c}} e_{\mathfrak{n}}(|d\varrho |)|s|^{2}_{h} dv_{bX_{c}}
\\
&\qquad \qquad +(1-\varepsilon )\left (\|(\widetilde{\nabla}^{TX})^{1,0}s
\|^{2}_{L^{2}(X_{c})}-2\|\sqrt{\tau}\widetilde{\nabla}^{TX}_{e^{(1,0)}_{
\mathfrak{n}}}s\|^{2}_{L^{2}(X_{c})}\right ).
\end{split}
\label{eq:4.42Jan25}
\end{equation}
Since
$\supp \, s \subset X_{\delta}\setminus \overline{X}_{-\delta}$,
$\sqrt{\tau}\in [0,1]$, and $\sqrt{2}e^{(1,0)}_{\mathfrak{n}}$ has norm
$1$ on the support of $\tau $,
\begin{equation}
\|(\widetilde{\nabla}^{TX})^{1,0}s\|^{2}_{L^{2}(X_{c})}-2 \|\sqrt{
\tau}\widetilde{\nabla}^{TX}_{e^{(1,0)}_{\mathfrak{n}}}s\|^{2}_{L^{2}(X_{c})}
\geq 0
\label{eq:4.43Jan2025}
\end{equation}
follows. Then, by {Lemmas~\ref{lem_metric} and \ref{lem-question2}}, we have
\begin{equation}
\begin{split} &e_{\mathfrak{n}}(|d\varrho |)\equiv 0 \,\text{ on }\, X_{
\delta}\setminus \overline{X}_{-\delta}\;,
\\
&\langle [\sqrt{-1}\mathscr{L}_{\varrho},\Lambda ]^{bX_{c}}_{
\varepsilon}s,s\rangle _{h} \geq C_{1}|s|^{2}_{h} \,\text{ on }\, bX_{c}
\,,
\\
&\left \langle [\sqrt{-1}R^{L}_{\chi},\Lambda ]_{\varepsilon} s,s
\right \rangle _{h}+ 2(1-\varepsilon )\tau \left \langle R^{L}_{\chi}(e^{(1,0)}_{
\mathfrak{n}}, e^{(0,1)}_{\mathfrak{n}})s,s\right \rangle _{h}\geq C_{0}|s|_{h}^{2}
\,\text{ on }\, X_{\delta}\setminus \overline{X}_{-\delta}\,.
\end{split}
\label{eq:4.44Jan2025}
\end{equation}
Combining {\eqref{eq:4.42Jan25}}--{\eqref{eq:4.44Jan2025}}, there exists
$C'>0$ such that for any $\varepsilon \in (0,\varepsilon _{1}]$, and any
$s\in \Omega _{0}^{0,j}(X_{\delta}\setminus \overline X_{-\delta},L^{k}
\otimes E)$ with $\iota _{e_{\mathfrak{n}}}s=0$ on $bX_{c}$, we have
\begin{equation}
\begin{split} (1+\varepsilon )\left (\|\overline \partial ^{E}_{k} s\|^{2}_{L^{2}(X_{c})}+
\|\overline \partial ^{E*}_{k} s\|^{2}_{L^{2}(X_{c})}\right )+ C(1+
\frac{1}{\varepsilon})\|s\|^{2}_{L^{2}(X_{c})}&
\\
\geq \left (C_{0} k - C' \right )\|s\|^{2}_{L^{2}(X_{c})} + C_{1} \|s
\|^{2}_{L^{2}(bX_{c})}&.
\end{split}
\label{eq:4.45Jan2025}
\end{equation}
The conclusion now follows from the inequality {\eqref{eq:4.45Jan2025}}.  
\end{proof}

\begin{rem}
\label{rem4.7}
Note that the inequality {\eqref{eq:4.45Jan2025}} is stronger than the one
in {Proposition~\ref{prop:4.3Jan25}} because of the boundary term
$C_{1} \|s\|^{2}_{L^{2}(bX_{c})}$. Actually, the inequality {\eqref{eq:4.45Jan2025}} implies the basic estimate in the sense of
\cite[Chapter II]{FK:72} for sufficiently large $k$.
\end{rem}

\subsection{Holomorphic Morse inequalities; Proof of {Theorems~\ref{thm:4.10important} and \ref{TheoremB}}}
\label{ss:4.3Jan25}

Now we prove {Theorem~\ref{thm:4.10important}} based on {Theorem~\ref{TheoremA}} and {Lemmas~\ref{lem_metric} and \ref{lem-question2}}.

\begin{proof}[Proof of {Theorem~\ref{thm:4.10important}}]
By the assumption that the domain $M$ is $q$-concave, as in the first part
of the proof of {Theorem~\ref{TheoremA}}, we get a defining function
$\varrho \in \mathscr{C}^{\infty}(X,\mathbb{R})$ for $M$ with the properties
mentioned there. In particular, set $K:=\overline{X}_{-\delta}$, so that
$\sqrt{-1}\partial \overline{\partial }\varrho $ has at least
$n-q+1$ negative eigenvalues at each point $x\in M\setminus K$. Let
$K_{L}\subset M$ be the compact subset as given in {Theorem~\ref{thm:4.10important}} such that $\sqrt{-1}R^{L}\leq 0$ on
$M\setminus K_{L}$. Then we can fix $u_{0}\in [-\delta ,0)$ such that
\begin{equation}
K\cup K_{L} \subset X_{u_{0}}.
\label{eq4.47}
\end{equation}
Then we fix a small interval $[u',v]$ of regular values of
$\varrho $ with $u_{0}< u'<0< v\leq \delta $. Fix another number $u$ such
that $u'< u<0 <v $. Now, we give a proof of {Theorem~\ref{thm:4.10important}} by refining the construction of
$h^{L}_{\chi}$ in the proof of {Lemma~\ref{lem-question2}}.

Note that we can apply {Lemmas~\ref{lem_metric} and \ref{lem-question2}} to $\varrho $ and the interval $[u',v]$. First, we
fix a modified Hermitian metric $\widetilde{\Theta}$ from {Lemma~\ref{lem_metric}} such that the results hold for all $c\in (u',v)$. In particular,
the unit vector field $e_{\mathfrak{n}}$ on $X_{v}\setminus X_{u'}$ is
defined as the metric dual of $-d\varrho $. Recall
that $M=X_{0}$.

Let $h^{L}$ be the smooth Hermitian metric given in the assumption of {Theorem~\ref{thm:4.10important}}. Then $\sqrt{-1}R^{L}\leq 0$ on
$M\setminus K_{L}$. As in the proof of {Lemma~\ref{lem-question2}}, we consider
a smooth real function $\chi $ and set a new Hermitian metric
$h^{L}_{\chi}:=h^{L}e^{-\chi (\varrho )}$. Then on
$M\setminus X_{u_{0}}$, we have the comparison of real Hermitian forms
\begin{equation}
\sqrt{-1}R^{L}_{\chi}\leq \sqrt{-1}\chi '(\rho )\partial
\overline \partial \varrho +\sqrt{-1}\chi ''(\varrho )\partial
\varrho \wedge \overline \partial \varrho .
\label{eq:4.47plus}
\end{equation}

The first condition for the function $\chi $ is that for
$t\in [u,v]$, we have
\begin{equation}
C_{2} \chi '(t)-C_{4} (1+|\chi ^{\prime }(t)|)\geq C_{0},
\label{eq:4.47Jan2025}
\end{equation}
where the constants $C_{0}$, $C_{2}$, $C_{4}$ are given in the proof of
{Lemma~\ref{lem-question2}}.

In addition to {\eqref{eq:4.47Jan2025}}, we also require the function
$\chi $ to satisfy the following conditions,
\begin{equation}
\chi (t)=0\, \text{ for } \, t \leq u',\quad \chi '(t)>0 \,
\text{ for }\, t > u'.
\label{eq:4.48Jan2025}
\end{equation}
Note that a smooth function
$\chi : \mathbb{R}\rightarrow \mathbb{R}$ satisfying {\eqref{eq:4.47Jan2025}} and {\eqref{eq:4.48Jan2025}} always exists, for example,
near $u'$, $\chi (t)$ could have a formula like $e^{-(t-u')^{-2}}$ for
$t>u'$.

As in the proof of {Lemma~\ref{lem_metric}}, on
$X_{v}\setminus X_{u'}$, we have the orthogonal splitting of $(1,0)$ tangent bundle
\begin{equation}
T^{(1,0)}X=H^{(1,0)}X\oplus \mathbb{C}\omega _{n},
\label{eq:4.51Jan-2}
\end{equation}
where $\omega _{n}=\sqrt{2}e^{(1,0)}_{\mathfrak{n}}$. Since
$\mathscr{L}_{\varrho}$ on $H^{(1,0)}X$ has at least $n-q$ negative eigenvalues,
we can find a local orthonormal family
$\{\omega _{\ell}\}_{\ell =1}^{n-q}$ in $H^{(1,0)}X$ such that, the restriction of $\sqrt{-1}\partial \overline{\partial }\varrho $ to the subspace
$V_{n-q}:=\mathrm{span}_{\mathbb{C}}\{\omega _{\ell},\; \ell =1,
\ldots , n-q\}$ is strictly
negative definite.

By our choice of $V_{n-q}\subset H^{(1,0)}X$, for $v\in V_{n-q}$, we always
have
\begin{equation}
\label{eq-inters}
(\partial \varrho \wedge \overline \partial \varrho )(v,\overline v)=0.
\end{equation}
Then by {\eqref{eq:4.47plus}}, we have at all points in
$M\setminus \overline{X}_{u'}$, for any nonzero $v\in V_{n-q}$,
\begin{equation}
\label{eq-preserve-curv}
R^{L}_{\chi}(v,\overline{v})=R^{L}(v,\overline{v})+\chi '(\varrho )
\partial \overline{\partial }\varrho (v,\overline{v})\leq \chi '(
\varrho )\partial \overline{\partial }\varrho (v,\overline{v}) <0,
\end{equation}
which is equivalent to $\sqrt{-1}R^{L}_{\chi}<0$ therein. As a consequence (and by Sylvester's law of inertia for Hermitian matrices),
for $j=0, \ldots , n-q-1$, we have
\begin{equation}
M(j,h^{L}_{\chi})=\overline{X}_{u'}(j, h^{L}_{\chi})=\overline{X}_{u'}(j,
h^{L}).
\label{eq4.54}
\end{equation}
Since $\sqrt{-1}R^{L}\leq 0$ on $\overline{X}_{u'}\setminus K_{L}$, we have
for $j=0, \ldots , n-q-1$,
\begin{equation}
M(j,h^{L}_{\chi})= \overline{X}_{u'}(j, h^{L})=K_{L}(j,h^{L}).
\label{eq:4.53Jan25}
\end{equation}
Moreover, on a neighborhood of $K_{L}$, we have
\begin{equation}
c_{1}(L, h^{L}_{\chi})=c_{1}(L,h^{L}).
\label{eq:4.54Jan25}
\end{equation}

Finally, combining the inequalities in {Theorem~\ref{TheoremA}} and {\eqref{eq:4.53Jan25}}, {\eqref{eq:4.54Jan25}}, we conclude exactly {\eqref{eq:4.58Jan25-1}} and {\eqref{eq:4.58Jan25-2}} for
$H^{0,\ell}_{(2)}(M,L^{k}\otimes E)$ defined with respect to the original
metrics $\Theta $, $h^{E}$, $h^{L}$, and the integrals on the right-hand
side are exactly on $K_{L}(\ell ,h^{L})$ or
$K_{L}(\leq \ell \, ,h^{L})$. Moreover, since $M$ satisfies the $Z(j)$ condition
for $0\leq j\leq n-q-1$, we apply {Theorem~\ref{zqzq1}} to identify
$H^{0,\ell}_{(2)}(M,L^{k}\otimes E)$ with
$H^{0,\ell}(\overline{M},L^{k}\otimes E)$ and with
$\mathscr{H}^{0,\ell}(M, L^{k}\otimes E)$. This completes the proof.
\end{proof}

\begin{rem}%
\label{rmk:4.9new}
Note that the equation {\eqref{eq-inters}} plays a crucial role in achieving {\eqref{eq-preserve-curv}}, which is a consequence of constructing
$\Theta $ adapted to the splitting {\eqref{eq:4.51Jan-2}}. Generally, it
is difficult to control $R^{L}_{\chi}$ in {\eqref{eq:4.47plus}}, since
$\sqrt{-1}\partial \varrho \wedge \overline \partial \varrho \geq 0$, while
the term $\sqrt{-1}\partial \overline \partial \varrho $ could be negative
along certain directions in $T^{(1,0)}X$.

Another important observation from {\eqref{eq-preserve-curv}} is that, instead of $\sqrt{-1}R^{L}\leq 0$ we
only need $\sqrt{-1}R^{L}|_{H^{(1,0)}X\otimes H^{(0,1)}X} \leq 0$ on
$X_{v}\setminus \overline{X}_{u'}$. 
\end{rem}

When $K_L=\varnothing $ in {Theorem~\ref{thm:4.10important}}, we obtain:
\begin{cor}%
\label{cor_qconcave}
If $c_{1}(L,h^{L})\leq 0$ on a smooth domain $M\Subset X$ of a complex
manifold $X$ of dimension $n \geq 2$ such that its Levi form
$\mathscr{L}_{bM}$ has at least $n-q$ negative eigenvalues on $bM$ with
$1\leq q\leq n-1$, then for $0\leq \ell \leq n-q-1$, as
$k\rightarrow \infty $,
\begin{equation}
\nonumber
\dim \mathscr{H}^{0,\ell}(M, L^{k}\otimes E)=o(k^{n}).
\end{equation}
Note that the spaces $\mathscr{H}^{0,\ell}(M,L^{k}\otimes E)$ above can
be replaced by the cohomology groups
$H^{0,\ell}(\overline{M},L^{k}\otimes E)$.
\end{cor}

Now we turn to {Theorem~\ref{TheoremB}}. Let $X$ be a $q$-concave manifold
equipped with a sublevel exhaustion function
$\varrho \in \mathscr{C}^{\infty}(X,\mathbb{R})$ and the associated exceptional
subset $K\Subset X$, then the $(1,1)$-form
$\partial \overline \partial \varrho $ has at least $n-q+1$ negative eigenvalues
on $X\setminus K$. Now we fix $u_{0} \in \mathbb{R}$ such that
$K \Subset X_{u_{0}}$. Then for a regular value $c>u_{0}$, by {Proposition~\ref{prop:2.3dec24}}, $Z(n-q-1)$ holds on $X_{c}$ and then we have for
$j\leq n-q-2$ (see {Theorem~\ref{zqzq1}} and {Theorem~\ref{thm:2.9AG}}),
\begin{equation}
H^{j}(X,L^{k}\otimes E)\simeq H^{0,j}(\overline X_{c},L^{k}\otimes E)
\simeq \mathscr{H}^{0,j}(X_{c},L^{k}\otimes E),
\label{eq:4.5Jan}
\end{equation}
where $\mathscr{H}^{0,j}(X_{c},\cdot )$ is the space of the harmonic forms
smooth up to the boundary of $X_{c}$ (for fixed smooth metrics
$\Theta $, $h^{L}$ and $h^{E}$). By {Theorem~\ref{thm:2.12H65}}, we have
\begin{equation}
\dim H^{n-q-1}(X,L^{k}\otimes E)\leq \dim H^{0,n-q-1}(\overline X_{c},L^{k}
\otimes E).
\label{eq:4.5Jan-new}
\end{equation}

\begin{proof}[Proof of {Theorem~\ref{TheoremB}}]
By the assumption, we take a sublevel exhaustion function
$\varrho \in \mathscr{C}^{\infty}(X,\mathbb{R})$ and the associated exceptional
subset $K\Subset X$. Let $u_{0}\in \mathbb{R}$ be such that
$K \cup K_{L} \subset X_{u_{0}}\Subset X$. Then we fix a regular value
$c\geq u_{0}$ of $\varrho $ such that $M:=X_{c}\Subset X$. Then $M$ is
a Levi $q$-concave domain with a smooth boundary.

By our assumption, $\sqrt{-1}R^{L}\leq 0$ on $M\setminus K_{L}$, we apply
{Theorem~\ref{thm:4.10important}} to this domain $M$ and line bundle
$(L,h^{L})$, combining with {\eqref{eq:4.5Jan}} and {\eqref{eq:4.5Jan-new}}, we obtain exactly {\eqref{eq:1.4Jan25}},
{\eqref{eq:1.3Jan25}} and {\eqref{eq:1.5Jan25}}.
\end{proof}

Combining {Theorem~\ref{TheoremB}} with \cite[\S 2.1]{CMW23}, we can also
obtain a singular version of {Theorem~\ref{TheoremB}}. See
\cite[\S 2.1]{CMW23} for the detailed definition of singular metrics with
algebraic singularities and the Nadel multiplier ideal sheaf.

\begin{cor}[Singular holomorphic Morse inequalities for $q$-concave manifolds]%
\label{cor:4.12}
Let $X$ be a connected $q$-concave manifold of dimension $n$. Let
$L,E$ be holomorphic vector bundles on $X$ with $\rank (L)=1$. Let
$h^{L}$ be a Hermitian metric on $L$ with algebraic singularities such
that $S(h^{L})\subset K$ and $c_{1}(L,h^{L})\leq 0$ on
$X\setminus K$ for some compact subset $K$. Then for
$\ell \leq n-q-2$,
\begin{equation}
\begin{split} \limsup _{k\rightarrow \infty}k^{-n} \dim H^{\ell}(X,L^{k}
\otimes E\otimes \mathscr{I}(h^{L^{k}}))& \leq \frac{\rank (E)}{n!}
\int _{K(\ell , h^{L})}(-1)^{\ell}c_{1}(L,h^{L})^{n},
\\
\limsup _{k\rightarrow \infty}\sum _{j=0}^{\ell }
\frac{(-1)^{\ell -j}}{k^{n}}\dim H^{j}(X,L^{k}\otimes E\otimes
\mathscr{I}(h^{L^{k}}))&\leq \frac{\rank (E)}{n!}\int _{K(\leq \ell
\, , h^{L})}(-1)^{\ell}c_{1}(L,h^{L})^{n},
\end{split}
\label{eq:4.61singular}
\end{equation}
where $\mathscr{I}(h^{L^{k}})$ is the Nadel multiplier ideal sheaf of
$h^{L^{k}}$, and
$K(\ell , h^{L}):=K \cap (X\setminus S(h^{L}))(\ell , h^{L})$. Moreover,
we have
\begin{align}
&\limsup _{k\rightarrow \infty}k^{-n} \dim H^{n-q-1}(X,L^{k}\otimes E
\otimes \mathscr{I}(h^{L^{k}}))\notag\\
&\quad \leq \frac{\rank (E)}{n!}\int _{K(n-q-1,
h^{L})}(-1)^{n-q-1}c_{1}(L,h^{L})^{n}.
\label{eq:4.62singular}
\end{align}
\end{cor}
The inequalities {\eqref{eq:4.61singular}} for $\ell \leq n-q-2$ already
appeared in \cite[Theorem 3.12]{CMW23} but the weak Morse inequality {\eqref{eq:4.62singular}} for degree $n-q-1$ is new.

\subsection{Volume of semi-positive line bundles;
Proofs of {Theorems~\ref{thm:4.13new} and \ref{thm:stability}}}
\label{ss:4.4semipositive}

In this subsection, we first present a proof of {Theorem~\ref{thm:4.13new}}. Our approach is based on the established proofs of {Theorems~\ref{thm:4.10important} and \ref{TheoremA}} as well as the subsequent corollary
derived from {Theorem~\ref{thm:4.10important}}, where we identify
$H^{0,0}(\overline{M},L^{k}\otimes E)$ with
$H^{0}(M,L^{k}\otimes E)$ by {Theorem~\ref{zqzq1}}.

\begin{cor}%
\label{cor:4.14new}
Let $X$ be a connected complex manifold of dimension $n\geq 3$. Let
$M$ be a relatively compact domain in $X$ with a smooth boundary
$bM$, and assume that $\varrho $ is a defining function for $bM$ with
$M:=\{x\in X\;:\; \varrho (x)<0\}$ (and $d\varrho $ does not vanish near
$bM$). Assume that $\mathscr{L}_{\varrho}$ is strictly negative on
$bM$. Let $E$ be a holomorphic vector bundle on $X$. Let $(L,h^{L})$ be
a Hermitian holomorphic line bundle on $X$ such that
$c_{1}(L,h^{L})\leq 0$ on $M\setminus K$ for some compact subset
$K\Subset M$. Then we have
\begin{equation}
\begin{split} &\limsup _{k\rightarrow \infty}\frac{n!}{k^{n}} \dim H^{0}(M,L^{k}
\otimes E)\leq \rank (E)\int _{K(0,h^{L})}c_{1}(L,h^{L})^{n},
\\
&\liminf _{k\rightarrow \infty}\frac{n!}{k^{n}} \dim H^{0}(M,L^{k}
\otimes E)\geq \rank (E)\int _{K(\leq 1\,,h^{L})} c_{1}(L,h^{L})^{n}.
\end{split}
\label{eq:4.58Jan25-4}
\end{equation}
\end{cor}

\begin{proof}[Proof of {\eqref{eq:4.61part1}}]
For any sufficiently small $\delta >0$, we first construct a Hermitian
metric $h^{L}_{\delta}$ such that $c_{1}(L,h^{L}_{\delta})\leq 0$ outside
$X_{-\delta}$; we then apply {\eqref{eq:4.58Jan25-4}} and let
$\delta \rightarrow 0$.

Let $\delta _{0}>0$ be small enough such that the assumptions for {\eqref{eq:4.61part1}} hold on
$\overline{M}_{\delta _{0}}\setminus M_{-\delta _{0}}\subset U$, where
$M_{c}:=\{x\in X\;:\; \varrho (x)<c\}$. We also assume that
$\phi |_{\overline{U}}>0$, and $c_{1}(L,h^{L})\geq 0$ on
$\overline{M}\cap \overline{U}$.

For $\delta \in (0,\delta _{0})$, let
$\chi _{\delta}:\mathbb{R}\rightarrow [0,1]$ be an increasing smooth function
such that $\chi _{\delta}(t)=0$ for $t\leq -\delta $, and
$\chi _{\delta}(t)=1$ for $t\geq -\delta /2$. We may also assume that
$\chi _{\delta}'(t)>0$ for $-\delta /2>t>-\delta $. Set a new Hermitian
metric on $L$ by
\begin{equation}
h^{L}_{\delta}:=h^{L}\exp (\chi _{\delta}(\varrho )\phi ).
\label{eq4.62}
\end{equation}
Then the corresponding Chern curvature form
\begin{equation}
R^{L}_{\delta}=R^{L}- \chi _{\delta}(\varrho )\partial
\overline{\partial }\phi -\chi '_{\delta}(\varrho ) (\partial \phi
\wedge \overline{\partial }\varrho +\partial \varrho \wedge
\overline{\partial }\phi )-\chi _{\delta}^{\prime }(\varrho )
\phi \partial \varrho \wedge \overline{\partial }\varrho -\chi '_{
\delta}(\varrho )\phi \partial \overline{\partial }\varrho .
\label{eq4.63}
\end{equation}
Then it is clear that
\begin{equation}
R^{L}_{\delta }=R^{L}\;\text{ on }\, M_{-\delta},\quad \text{ and }\, R^{L}_{
\delta}\equiv 0 \;\text{ on }\, M\setminus M_{-\delta /2}.
\label{eq4.64}
\end{equation}
Hence {Corollary~\ref{cor:4.14new}} applies.

For any $c\in (-\delta ,-\frac{\delta}{2})$ and restricting to
$T^{(1,0)}bM_{c}\otimes T^{(0,1)}bM_{c}$, we have
\begin{equation}
R^{L}_{\delta}|_{bM_{c}}=(1-\chi _{\delta}(c))R^{L}|_{bM_{c}}-\chi '_\delta(c)
\phi \mathscr{L}_{\varrho}|_{bM_{c}}>0 .
\label{eq4.65}
\end{equation}
As a consequence, on $M_{-\delta /2}\setminus M_{-\delta}$, the Hermitian
form $\sqrt{-1}R^{L}_{\delta}$ has at most one negative eigenvalue. Then
\begin{equation}
\begin{split} \int _{(M_{-\delta /2}\setminus M_{-\delta})(\,\leq 1,h^{L}_{
\delta})} c_{1}(L,h^{L}_{\delta})^{n}= \int _{M_{-\delta /2}
\setminus M_{-\delta}} c_{1}(L,h^{L}_{\delta})^{n}.
\end{split}
\label{eq4.66}
\end{equation}
Write
\begin{equation*}
(R^{L}_{\delta})^{n}= \left (\partial \overline \partial (\phi -\chi _{
\delta}(\varrho )\phi )\right )^{n}=-d\left \{\left (\partial
\overline \partial (\phi -\chi _{\delta}(\varrho )\phi )\right )^{n-1}
\partial (\phi -\chi _{\delta}(\varrho )\phi )\right \}
\end{equation*}
By Stokes' theorem and our condition on $\chi _{\delta}$, we get
\begin{equation}
\begin{split} \int _{M_{-\delta /2}\setminus M_{-\delta}} c_{1}(L,h^{L}_{
\delta})^{n} &= \frac{\sqrt{-1}}{2\pi}\int _{bM_{-\delta}} c_{1}(L,h^{L})^{n-1}\wedge
\partial \phi .
\end{split}
\label{eq4.67}
\end{equation}
Then by the second line of {\eqref{eq:4.58Jan25-4}}, we get
\begin{equation}
\begin{split} &\liminf _{k\rightarrow \infty}\frac{n!}{k^{n}} \dim H^{0}(M,L^{k}
\otimes E)
\\
&\quad \geq \rank (E)\int _{M_{-\delta}(\leq 1, h^{L})}c_{1}(L,h^{L})^{n}+
\rank (E)\frac{\sqrt{-1}}{2\pi}\int _{bM_{-\delta}} c_{1}(L,h^{L})^{n-1}
\wedge \partial \phi .
\end{split}
\label{eq:4.69Jan-25}
\end{equation}
Finally {\eqref{eq:4.61part1}} follows from making
$\delta \rightarrow 0$ in {\eqref{eq:4.69Jan-25}}.
\end{proof}

\begin{rem}
\label{rem4.12}
From the proof of {\eqref{eq:4.61part1}}, we also obtain that
\begin{equation}
\begin{split} &\limsup _{k\rightarrow \infty}\frac{n!}{k^{n}} \dim H^{0}(M,L^{k}
\otimes E)
\\
&\quad \leq \rank (E)\int _{M_{-\delta /2}(0,h^{L}_{\delta})}c_{1}(L,h^{L}_{
\delta})^{n}
\\
&\quad =\rank (E)\int _{M_{-\delta}(0,h^{L})}c_{1}(L,h^{L})^{n} +
\rank (E)\int _{(M_{-\delta /2}\setminus M_{-\delta})(0,h^{L}_{\delta})}c_{1}(L,h^{L}_{
\delta})^{n}.
\end{split}
\label{eq:4.71Jan-25}
\end{equation}
The integral over
$(M_{-\delta /2}\setminus M_{-\delta})(0,h^{L}_{\delta})$ appears to be
difficult to compute explicitly, making it challenging to determine the
limit as $\delta \rightarrow 0$.
\end{rem}

The proof of {\eqref{eq:4.62part2}} is more intricate compared to the preceding
one, requiring a thorough review of the arguments presented in the proofs
of {Theorems~\ref{TheoremB} and \ref{TheoremA}}.
\begin{proof}[Proof of {\eqref{eq:4.62part2}}]
We use the same notations as in the proof of {\eqref{eq:4.61part1}}. Let
$\delta '>0$ be a small constant such that on
$M_{\delta '}\setminus M_{-\delta '}$, we have
\begin{equation*}
R^{L}=-\partial \overline \partial \varrho \,.
\end{equation*}
Note that this implies that
$\sqrt{-1}\partial \overline \partial \varrho \leq 0$.

Fix any $\delta \in (0,\delta ')$, and take a smooth real function
$\eta $ on $\mathbb{R}$ such that
\begin{itemize}
\item $\eta (t)=0$ when $t\leq -\delta $;
\item $\eta '(-\delta /2)= 1$;
\item $\eta '(t)>0$ and $\eta ^{\prime }(t)>0$, for
$t>-\delta $.
\end{itemize}
We then consider the following Hermitian metric
$h^{L}_{\eta}:=h^{L}e^{-\eta (\varrho )}$ on $L$. Its curvature is
\begin{equation}
R^{L}_{\eta}=R^{L}+\eta '(\varrho )\partial \overline \partial
\varrho +\eta ^{\prime }(\varrho )\partial \varrho \wedge
\overline \partial \varrho .
\label{eq:4.73new}
\end{equation}
We have the following properties:
\begin{itemize}
\item $R^{L}_{\eta}=R^{L}$ on $M_{-\delta}$.
\item $R^{L}_{\eta }\geq 0$ on $M_{-\delta /2}$.
\item $R^{L}_{\eta}$ has at least $2$ negative eigenvalues on
$M_{\delta}\setminus \overline{M}_{-\delta /2}$ (since $n\geq 3$).
\end{itemize}

Following {Remark~\ref{rmk:4.9new}} and repeating the arguments in the
proofs of {Proposition~\ref{prop:4.3Jan25}} and {Theorem~\ref{thm:4.10important}}, we first obtain the optimal fundamental estimates
for $(0,0)$ and $(0,1)$-forms on $M_{0}=M$, with the core subset
$K=\overline{M}_{-\delta /2}$ and with respect to a modified Hermitian
metric on $L$ that coincides with the above $h^{L}_{\eta}$ on
$\overline{M}_{-\delta /2}$; then we obtain for any
$\delta \in (0,\delta ')$ that
\begin{equation}
\begin{split} &\limsup _{k\rightarrow \infty}\frac{n!}{k^{n}} \dim H^{0}(M,L^{k}
\otimes E)\leq \rank (E)\int _{M_{-\delta /2}(0,h^{L}_{\eta})}c_{1}(L,h^{L}_{
\eta})^{n},
\\
&\liminf _{k\rightarrow \infty}\frac{n!}{k^{n}} \dim H^{0}(M,L^{k}
\otimes E)\geq \rank (E)\int _{M_{-\delta /2}(\leq 1\,,h^{L}_{\eta})} c_{1}(L,h^{L}_{
\eta})^{n}.
\end{split}
\label{eq:4.74Jan25}
\end{equation}

By our construction, we have
$M_{-\delta /2}(0,h^{L}_{\eta})=M_{-\delta /2}(\leq 1\,,h^{L}_{\eta})$,
so that
\begin{equation}
\begin{split} \int _{M_{-\delta /2}(0,h^{L}_{\eta})}c_{1}(L,h^{L}_{
\eta})^{n} & =\int _{M_{-\delta /2}} c_{1}(L,h^{L}_{\eta})^{n}
\\
&=\int _{M_{-\delta}} c_{1}(L,h^{L})^{n} +\int _{M_{-\delta /2}
\setminus M_{-\delta}} c_{1}(L, h^{L}_{\eta})^{n}.
\end{split}
\label{eq:4.75Jan25}
\end{equation}
By {\eqref{eq:4.73new}}, we have
\begin{equation}
\begin{split} &\int _{M_{-\delta /2}\setminus M_{-\delta}} c_{1}(L, h^{L}_{
\eta})^{n}
\\
=&\int _{M_{-\delta /2}\setminus M_{-\delta}} (1-\eta '(\varrho ))^{n}c_{1}(L,h^{L})^{n}
\\
&+ \frac{\sqrt{-1}}{2\pi}\int _{M_{-\delta /2}\setminus M_{-\delta}} n(1-
\eta '(\varrho ))^{n-1} \eta ''(\varrho ) \partial \varrho \wedge
\overline \partial \varrho \wedge c_{1}(L,h^{L})^{n-1}.
\end{split}
\label{eq:4.76Jan25}
\end{equation}
It is clear that there exists a constant $C$ such that for all
$\delta \in (0,\delta ')$, we have
\begin{equation}
\left | \int _{M_{-\delta /2}\setminus M_{-\delta}} (1-\eta '(
\varrho ))^{n}c_{1}(L,h^{L})^{n} \right |\leq C\delta .
\label{eq:4.77Jan25}
\end{equation}

Another observation, by calculus in local coordinates, is that there exists
a constant $C'>0$ independent of $\delta $ and $\eta $ such that
\begin{equation}
\begin{split} \Big|&\frac{\sqrt{-1}}{2\pi}\int _{bM_{-\delta}} c_{1}(L,h^{L})^{n-1}
\wedge \partial \varrho
\\
& \quad + \frac{\sqrt{-1}}{2\pi}\int _{M_{-\delta /2}\setminus M_{-
\delta}} n(1-\eta '(\varrho ))^{n-1} \eta ''(\varrho ) \partial
\varrho \wedge \overline \partial \varrho \wedge c_{1}(L,h^{L})^{n-1}
\Big |\leq C' \delta .
\end{split}
\label{eq:4.78Jan25}
\end{equation}
In fact, since all points in $[-\delta ,0]$ are regular values of
$\varrho $, after using the gradient flow of $\varrho $ starting from
$bM_{-\delta}$, we can identify
\begin{equation*}
M_{-\delta /2}\setminus M_{-\delta} \simeq bM_{-\delta} \times [-
\delta ,-\delta /2)
\end{equation*}
so that for $t\in [-\delta ,-\delta /2)$,
\begin{equation*}
M_{t}\simeq bM_{-\delta}\times \{t\}.
\end{equation*}
Equivalently, for
$(y,t)\in bM_{-\delta} \times [-\delta ,-\delta /2)\simeq M_{-\delta /2}
\setminus M_{-\delta}$, we have $\varrho (y,t)=t$. We fix a smooth volume
form $d\sigma $ on $bM_{-\delta}$, and we write in the coordinate
$(y,t)$
\begin{equation*}
\frac{\sqrt{-1}}{2\pi}\partial \varrho \wedge \overline \partial
\varrho \wedge c_{1}(L,h^{L})^{n-1}=a(y,t)dt\wedge d\sigma (y).
\end{equation*}
Set
\begin{equation*}
A(t):=\int _{bM_{-\delta}}a(y,t)d\sigma (y).
\end{equation*}
In particular,
\begin{equation}
\frac{\sqrt{-1}}{2\pi}\int _{bM_{-\delta}} c_{1}(L,h^{L})^{n-1}
\wedge \partial \varrho =-A(-\delta ).
\label{eq4.76}
\end{equation}
Therefore, using the fact
\begin{equation}
\int _{-\delta}^{-\delta /2}n(1-\eta '(t))^{n-1}\eta ''(t)dt=1,
\label{eq:2026-6-3}
\end{equation}
we deduce
\begin{align}
 &\frac{\sqrt{-1}}{2\pi}\int _{bM_{-\delta}} c_{1}(L,h^{L})^{n-1}
\wedge \partial \varrho
\notag \\
& \quad + \frac{\sqrt{-1}}{2\pi}\int _{M_{-\delta /2}\setminus M_{-
\delta}} n(1-\eta '(\varrho ))^{n-1} \eta ''(\varrho ) \partial
\varrho \wedge \overline \partial \varrho \wedge c_{1}(L,h^{L})^{n-1}
\notag \\
&=\int _{-\delta}^{-\delta /2}n(1-\eta '(t))^{n-1}\eta ''(t)\left (A(t)-A(-
\delta )\right )dt
\label{eq:4.78Jan25-11}
\end{align}
Then {\eqref{eq:4.78Jan25}} follows from {\eqref{eq:2026-6-3}},
{\eqref{eq:4.78Jan25-11}} and the facts that there exists a constant
$C'>0$ independent of the function $\eta $ and sufficiently small
$\delta >0$ such that for $t\in [-\delta ,-\delta /2]$,
\begin{equation*}
(1-\eta '(t))^{n-1}\eta ''(t)\geq 0,\,\quad |A(t)-A(-\delta )|\leq C'
\delta .
\end{equation*}

Finally, combining {\eqref{eq:4.74Jan25}}--{\eqref{eq:4.78Jan25}} and letting
$\delta \rightarrow 0$, we conclude exactly {\eqref{eq:4.62part2}}.
\end{proof}

Now we provide the proof of {Theorem~\ref{thm:stability}}.
\begin{proof}[Proof of {Theorem~\ref{thm:stability}}]
The proof follows directly from \cite[\S 5]{M:2016-q}, relying on the fact
that, after a small deformation of the complex structures on $L$ and
$M$, the Chern curvature $c_{1}(L,h^{L})$, the defining function of
$M$, and the local potential $\phi $ remain close in the
$\mathscr{C}^{\infty}$-topology. Since $c_{1}(L,h^{L})$ is strictly positive
and $bM$ strictly pseudoconcave, these properties are preserved under such
deformations.

Moreover, since $ c_{1}(L,h^{L})>0$ on $X_{\mathrm{reg}}$, the integral
in {\eqref{eq:1.10new2025}} is strictly positive with respect to the original
complex structures. Thus, after sufficiently small deformation of complex
structures, the strict inequality {\eqref{eq:4.62partI}} remains valid. Then
we can apply {Corollary~\ref{cor:1.7new}} to $(M,J')$.
\end{proof}

\subsection{Weak holomorphic Morse inequalities for 
\texorpdfstring{$(p,q)$}{(p,q)}-convex-concave domains}
\label{ss:4.5pqcorona}

By Definition~\ref{def:Zq}, when $M$ satisfies the $Z(q)$ condition, the
boundary points look either $(n-q-1)$-concave 
or $q$-convex. It could
be an interesting question to combine the estimate in {Proposition~\ref{prop:4.3Jan25}} for $(n-q-1)$-concave points and the estimate in
\cite[(3.5.19)]{MM} for $q$-convex points to obtain an optimal fundamental
estimate for $(0,q)$-forms on $M$. In this subsection, we aim to apply
our methods from the previous subsections to study the domains with mixed
convexity/concavity, which shall be simpler than the case of the general
$Z(q)$ condition.

Following \cite{Andreotti-Siu} and the {Definition~\ref{defmfd}} of
$q$-convex/concave functions, we introduce the following notions:
\begin{defn}
\label{def:pqcorona}
Let $X$ be a connected complex manifold of dimension $n\geq 1$ and let
$p,q\in \{1,\ldots ,n\}$.

(i) The complex space $X$ is called \textit{$(p,q)$-convex-concave} if there
exists a smooth real function $\varphi : X\rightarrow (a,b)$ with
$a\in \mathbb{R}\cup \{-\infty \}$ and
$a<b\in \mathbb{R}\cup \{+\infty \}$ such that
\begin{itemize}
\item For any $c,d$ with $a<c<d<b$, we have
$X^{d}_{c}:=\{x\in X\;:\; c<\varphi (x)<d\}\Subset X$.
\item There exist $a_{0},b_{0}$ with $a<a_{0}<b_{0}<b$ such that
$\sqrt{-1}\partial \overline \partial \varphi $ has at least $n-p+1$ positive
eigenvalues on $\{\varphi > b_{0}\}$, and has at least $n-q+1$ positive
eigenvalues on $\{\varphi < a_{0}\}$.
\end{itemize}
In this case, we call $\varphi $ an \textit{exhaustion function} for
$X$ with exceptional interval $[a_{0}, b_{0}]$.

(ii) If $M\Subset X$ is a relatively compact domain with a smooth boundary,
and we write $bM=b_{+}M\cup b_{-}M$ where $b_{+}M$ and $b_{-}M$ are disjoint
from each other and consist of connected components of $bM$. We call
$M$ a \textit{$(p,q)$-corona} with respect to the splitting
$bM=b_{+}M\cup b_{-}M$ if, near $b_{+}M$, there is a local smooth definition
function $r_{+}$ for $M$ such that $r_{+}$ is $p$-convex; and near
$b_{-}M$, there is a local smooth definition function $r_{-}$ for
$M$ such that $r_{-}$ is $q$-concave.
\end{defn}

Note that when $p+q\leq n-1$, then a $(p,q)$-corona $M$ always satisfies
the $Z(j)$ condition for $p\leq j\leq n-q-1$. When
$b_{+}M=\varnothing $ or $b_{-}M=\varnothing $, we get $M$ to be $q$-concave
or $p$-convex, respectively. If $p=1$, then $b_{+}M$ (if nonempty) is a
strongly pseudoconvex boundary component; if $q=1$, $b_{-}M$ (if nonempty) is a strongly
pseudoconcave boundary component. The following result is clear by definition and
by {Proposition~\ref{prop:2.3dec24}}.
\begin{lemma}%
\label{lem:4.17AS}
If $X$ is a $(p,q)$-convex-concave complex manifold with an exhaustion
function $\varphi $ and an exceptional interval
$[a_{0},b_{0}]\subset (a,b)$, then for any regular values $a'$ and
$b'$ of $\varphi $ with $a<a'<a_{0}$ and $b_{0}<b'<b$, the domain
$X^{b'}_{a'}=\{x\in X\;:\; a' < \varphi (x)<b'\}$ is a $(p,q)$-corona in
$X$ with $b_{+}X^{b'}_{a'}=\varphi ^{-1}(b')$ and
$b_{-}X^{b'}_{a'}=\varphi ^{-1}(a')$.
\end{lemma}

Analogous to {Theorems~\ref{zqzq1}, \ref{thm:2.9AG}, and \ref{thm:2.12H65}}, we have the following result (see also
\cite[Propositions 21, 22 and Th\'{e}or\`{e}me 15]{AG:62},
\cite[Theorem 3.4.9 and the subsequent remark]{Hor:65}, and
\cite[Proposition 1.2]{Andreotti-Siu}).
\begin{prop}
\label{prop:4.18pqcoh}
With the same notation and hypotheses as in {Lemma~\ref{lem:4.17AS}}, let
$F$ be a holomorphic vector bundle on $X$. If $p+q\leq n-1$, then we have
\begin{itemize}
\item For $p\leq j\leq n-q-1$, the isomorphism among cohomology groups
\begin{equation}
H^{j}(X,F)\simeq H^{j}(X^{b'}_{a'},F)\simeq H^{0,j}(X^{b'}_{a'},F).
\label{eq4.79}
\end{equation}
\item For $p\leq j\leq n-q-2$, the isomorphism among cohomology groups
\begin{equation}
H^{j}(X,F)\simeq H^{0,j}_{(2)}(X^{b'}_{a'},F)\simeq H^{0,j}(
\overline{X}^{b'}_{a'},F).
\label{eq4.80}
\end{equation}
\item For $j= n-q-1$,
\begin{equation}
\dim H^{n-q-1}(X,F)\leq \dim H^{0,n-q-1}_{(2)}(X^{b'}_{a'},F)=\dim H^{0,n-q-1}(
\overline{X}^{b'}_{a'},F).
\label{eq4.81}
\end{equation}
\end{itemize}
\end{prop}

We can now extend the weak holomorphic Morse inequalities presented in
{Theorem~\ref{thm:4.10important}} to the setting of the $(p,q)$-corona.

\begin{thm}%
\label{thm:4.19pqcorona}
Let $(X,\Theta )$ be a Hermitian complex manifold of dimension
$n\geq 2$. Take $p,q\in \{1,\ldots ,n-1\}$ such that $p+q\leq n-1$. Let
$M\Subset X$ be a relatively compact $(p,q)$-corona with a smooth boundary
in $X$ with respect to the splitting $bM=b_{+}M\cup b_{-}M$. Let
$(E, h^{E})$, $(L,h^{L})$ be two Hermitian holomorphic vector bundles on
$X$ with $\rank (L)=1$.
\begin{itemize}
\item There exists a Hermitian metric $h^{L}_{\chi}$ on $L$ that coincides
with $h^{L}$ away from $bM$, such that for $p\leq j\leq n-q-1$,
\begin{equation}
\limsup _{k\rightarrow \infty}k^{-n} \dim H^{0,j}(\overline{M},L^{k}
\otimes E)\leq \frac{\rank (E)}{n!}\int _{M(j,h^{L}_{\chi})}(-1)^{j}c_{1}(L,h^{L}_{
\chi})^{n}.
\label{eq:4.80Jan25-1}
\end{equation}
\item If $c_{1}(L,h^{L})\leq 0$ on an open neighborhood of $b_{-}M$ in
$\overline{M}$, and $c_{1}(L,h^{L})\geq 0$ on an open neighborhood of
$b_{+}M$ in $\overline{M}$, then for $p\leq j\leq n-q-1$,
\begin{equation}
\limsup _{k\rightarrow \infty}k^{-n} \dim H^{0,j}(\overline{M},L^{k}
\otimes E)\leq \frac{\rank (E)}{n!}\int _{M(j,h^{L})}(-1)^{j}c_{1}(L,h^{L})^{n}.
\label{eq:4.82Jan25-1}
\end{equation}

Moreover, if $p=1$, we also have
\begin{equation}
\liminf _{k\rightarrow \infty}k^{-n} \dim H^{0,0}(\overline{M},L^{k}
\otimes E)\geq \frac{\rank (E)}{n!}\int _{M(\leq 1\,,h^{L})} c_{1}(L,h^{L})^{n}.
\label{eq:4.83Jan25-2}
\end{equation}
\end{itemize}
Note that the spaces $H^{0,j}(\overline{M},L^{k}\otimes E)$ above can be
replaced by the spaces of $\overline \partial $-Neumann harmonic forms
$\mathscr{H}^{0,j}(M,L^{k}\otimes E)$ induced by $\Theta $, $h^{E}$,
$h^{L}$.
\end{thm}
\begin{proof}
Since $b_{+}M$ and $b_{-}M$ are disjoint from each other, we could take
any small open neighborhoods $U_{+}$ and $U_{-}$ of them in $X$ respectively
such that $\overline{U}_{+}\cap \overline{U}_{-}=\varnothing $. Then on
$U_{+}$, by $p$-convex condition, we proceed as in
\cite[Lemmas 3.5.3, 3.5.4 and Theorem 3.5.8]{MM} to modify $\Theta $ and
$h^{L}$ locally on $U_{+}$, so that if
$s\in B^{0,j}(M,L^{k}\otimes E)$ for $j\geq p$ and
$\mathrm{supp}\, (s)\subset U_{+}$, we have
\begin{equation}
\begin{split} \|\overline \partial ^{E}_{k} s\|^{2}_{L^{2}(M)}+\|
\overline \partial ^{E*}_{k} s\|^{2}_{L^{2}(M)} \geq \left (C_{0} k - C'
\right )\|s\|^{2}_{L^{2}(M)}.
\end{split}
\label{eq:4.95Jan2025}
\end{equation}

Then on $U_{-}$, by $q$-concave condition, we proceed as in the proof of
{Proposition~\ref{prop:4.3Jan25}} to modify $\Theta $ and $h^{L}$ locally
on $U_{-}$, so that if $s\in B^{0,j}(M,L^{k}\otimes E)$ for
$j\leq n-q-1$ and $\mathrm{supp}\, (s)\subset U_{-}$, we have
\begin{equation}
\begin{split} \|\overline \partial ^{E}_{k} s\|^{2}_{L^{2}(M)}+\|
\overline \partial ^{E*}_{k} s\|^{2}_{L^{2}(M)} \geq \left (C_{0} k - C'
\right )\|s\|^{2}_{L^{2}(M)}.
\end{split}
\label{eq:4.95Jan2025-2}
\end{equation}
Where $C_{0}>0$ and $C'>0$ are constants independent of $k$. Let $V\Subset M$ contains
$M\setminus (\overline{U}_{+}\cup \overline{U}_{-})$. Combining {\eqref{eq:4.95Jan2025}} and {\eqref{eq:4.95Jan2025-2}}, we obtain that there
exists a constant $\widetilde{C}>0$ such that for
$p\leq j\leq n-q-1$, $s\in B^{0,j}(M,L^{k}\otimes E)$, $k\gg 0$,
\begin{equation}
\label{eq-opfe2}
(1-\frac{\widetilde{C}}{k})\|s\|_{L^{2}(M)}^{2}\leq
\frac{\widetilde{C}}{k}\left (\|\overline \partial ^{E}_{k} s\|_{L^{2}(M)}^{2}+
\|\overline \partial ^{E*}_{k} s\|_{L^{2}(M)}^{2} \right )+\int _{V} |s|^{2}
dv_{X},
\end{equation}
where we use the modified metrics $\widetilde{\Theta}$ and
$h^{L}_{\chi}$. Then {\eqref{eq:4.80Jan25-1}} follows from {Corollary~\ref{cor:2.13new}} and {\eqref{eq-opfe2}}. Moreover, if $p=1$, we also obtain
by {Theorem~\ref{thm-l2hmi-2}} that
\begin{equation}
\liminf _{k\rightarrow \infty}k^{-n} \dim H^{0,0}(\overline{M},L^{k}
\otimes E) \geq \frac{\rank (E)}{n!}\int _{V(\leq 1\,,h^{L}_{\chi})} c_{1}(L,h^{L}_{
\chi})^{n}.
\label{eq:4.83Jan25-4}
\end{equation}

Using the same method as in the proof of {Theorem~\ref{thm:4.10important}} and \cite[Theorem 3.5.8]{MM}, in particular, we
construct the function $\chi $ (according to the choice of
$V\Subset M$) so that results analogous to {\eqref{eq-preserve-curv}}--{\eqref{eq:4.54Jan25}} and
\cite[(3.5.20)]{MM} hold. More precisely, for $p\leq j \leq n-q-1$, we
have $c_{1}(L,h^{L}_{\chi})=c_{1}(L,h^{L})$ on
\begin{equation}
V(j\,,h^{L}_{\chi})=V(j\,,h^{L}).
\label{eq:VjMay25}
\end{equation}
Therefore, we obtain {\eqref{eq:4.82Jan25-1}} from {\eqref{eq:4.80Jan25-1}}, {\eqref{eq:VjMay25}} and the fact that
$(-1)^{j}c_{1}(L,h^{L})^{n}$ is non-negative on $M(j,h^{L})$.

Finally, we need to prove {\eqref{eq:4.83Jan25-2}}. By {\eqref{eq:4.83Jan25-4}} and {\eqref{eq:VjMay25}}, we obtain for a given
$V\Subset M$ as chosen in the beginning
\begin{equation}
\liminf _{k\rightarrow \infty}k^{-n} \dim H^{0,0}(\overline{M},L^{k}
\otimes E)\geq \frac{\rank (E)}{n!}\int _{V(\leq 1\,,h^{L})} c_{1}(L,h^{L})^{n}.
\label{eq:4.83Jan25-5}
\end{equation}
Since we can choose the open subset $U_{\pm}$ to be sufficiently small
(due to the smoothness of the boundary), then by letting $V$ exhaust
$M$ the integral in the right-hand side of {\eqref{eq:4.83Jan25-5}} converges
to the integral in {\eqref{eq:4.83Jan25-2}} and this finishes the proof.
\end{proof}

As a direct consequence of {Proposition~\ref{prop:4.18pqcoh}} and {Theorem~\ref{thm:4.19pqcorona}}, we obtain an extension of the weak holomorphic
Morse inequalities in {Theorem~\ref{TheoremB}} as follows.

\begin{thm}%
\label{TheoremB-pqcase}
Consider integers $n\geq 2$, $p, q\geq 1$ with
$1\leq p+q \leq n-1$, and let $X$ be a $n$-dimensional $(p,q)$-convex-concave
connected complex manifold. Using the same notation as in {Definition~\ref{def:pqcorona}}, let $\varphi : X\rightarrow (a,b)$ be an exhaustion
function with the exceptional interval $[a_{0}, b_{0}]$. Let $E,L$ be holomorphic
vector bundles on $X$ with $\rank (L)=1$. Let $h^{L}$ be a smooth Hermitian
metric on $L$ such that there are $a'<b'$ satisfying
$a<a'\leq a_{0}$, $b_{0}\leq b'<b$ such that $c_{1}(L,h^{L})\leq 0$ on
$\{\varphi < a'\}$ and $c_{1}(L,h^{L})\geq 0$ on $\{\varphi > b'\}$. Then
for $p\leq j\leq n-q-1$, we have the weak Morse inequalities:
\begin{equation}
\begin{split} \limsup _{k\rightarrow \infty}k^{-n} \dim H^{j}(X,L^{k}
\otimes E)\leq \frac{\rank (E)}{n!}\int _{X_{a'}^{b'}(j,h^{L})}(-1)^{j}c_{1}(L,h^{L})^{n}&
\\
\leq \frac{\rank (E)}{n!}\int _{X(j,h^{L})}(-1)^{j}c_{1}(L,h^{L})^{n}.&
\end{split}
\label{eq:4.85Jan25}
\end{equation}
Moreover, if $p=1$, we also have
\begin{equation}
\liminf _{k\rightarrow \infty}k^{-n} \dim H^{0}(X,L^{k}\otimes E)
\geq \frac{\rank (E)}{n!}\int _{X(\leq 1\,,h^{L})} c_{1}(L,h^{L})^{n},
\label{eq:4.86Jan25}
\end{equation}
where the right or left side might be $+\infty $.
\end{thm}

\begin{rem}
\label{rem4.19}
Note that analogous versions of {Theorems~\ref{thm:4.13new},
\ref{thm_vanishing}} and {\ref{vanish}} can also be established for
$(p,q)$-coronas or $(p,q)$-convex-concave manifolds. Moreover, by combining
the above strategy with the techniques in \cite{Oh:82} and
\cite{M:96}, the results can be extended to compact complex spaces with
smooth boundary and isolated singularities, by constructing a complete
Hermitian metric in a punctured neighborhood of each singularity.
\end{rem}

\section{Proof of {Theorem~\ref{thmextension1}}}
\label{proofddbar1}

Let $(X,\Theta )$ be an $n$-dimensional Hermitian complex manifold ($n
\geq 2$). Let $(E,h^{E})$ and $(L,h^{L})$ be holomorphic Hermitian vector
bundles over $X$ with $\rank L=1$. Let $M$ be a relatively compact domain
in $X$ with smooth boundary $bM$. We denote by $L^{-1}=L^{*}$ and
$E^{*}$ the dual bundle of $L$ and $E$, respectively. We also recall our assumptions
in {Theorem~\ref{thmextension1}} as follows:
\begin{itemize}
\item[(A)] $L$ is semi-positive on $X$ and positive at one point;
\item[(B)] The Levi form $\mathscr{L}_{\varrho}$ of a defining function
for $M\Subset X$ has at least $n-q$ negative eigenvalues on $bM$ ($1
\leq q\leq n-1$).
\end{itemize}
Note that (B) implies that $M$ satisfies $Z(j)$ for
$0\leq j\leq n-q-1$.

We need the following two lemmas to prove {Theorem~\ref{thmextension1}}.

\begin{lemma}
\label{propprojection}
Assume $X$ is compact. Under the assumptions (A) and (B), for the harmonic
projection with $q \leq \ell \leq n-1$,
\begin{equation}
H: H^{0}(X,L^{k})\times \mathscr{H}^{0,n-\ell -1}({M},L^{-k}\otimes E^{*})
\rightarrow \mathscr{H}^{0,n-\ell -1}({M,E^{*}}),
\label{eq5.1}
\end{equation}
there exists a non-zero holomorphic section
$s\in H^{0}(X,L^{k_{0}})$ for some $k_{0}\in \mathbb{N}$ such that
$H(s\theta )=0$ for every
$\theta \in \mathscr{H}^{0,n-\ell -1}({M},L^{-k_{0}}\otimes E^{*})$.
\end{lemma}
\begin{proof}
To simplifying notation, we set $V:=H^{0}(X,L^{k})$,
$U:=\mathscr{H}^{0,n-\ell -1}({M},E^{*})$ and
$W:=\mathscr{H}^{0,n-\ell -1}({M},L^{-k}\otimes E^{*})$. Moreover, we define
a bilinear map $F(s,\alpha ):=H(s\alpha )$ for $s\in V$ and
$\alpha \in W$. Then we obtain a linear map
$G: V\rightarrow W^{*} \otimes U$ by $G(s):=F(s,\cdot )$. Suppose the assertion
were false, that is, for any $k \in \mathbb{N}$ and any non-zero
$s\in V$, there exists $\alpha _{0}\in W$ such that
$F(s,\alpha _{0})\neq 0$. If $G(s_{0})=0$ for some $s_{0}\in V$, then
$F(s_{0},\alpha ) =0$ for any $\alpha \in W$. Thus $s_{0}$ is zero, that
is, $G$ is injective. And it follows that
$\dim V\leq \dim (W^{*} \otimes U)=\dim W \times \dim U$, that is, for
any $k\in \mathbb{N}$,
\begin{equation}
\dim H^{0}(X,L^{k}) \leq \dim \mathscr{H}^{0,n-\ell -1}({M},L^{-k}
\otimes E^{*}) \times \dim \mathscr{H}^{0,n-\ell -1}({M},E^{*}).
\label{eq5.2}
\end{equation}

From {Corollary~\ref{cor_qconcave}}, the assumptions (B), $L\geq 0$ on
$X$ and $q \leq \ell $ imply that as $k\rightarrow \infty $,
\begin{equation}
\dim \mathscr{H}^{0,n-\ell -1}({M},L^{-k}\otimes E^{*})=o(k^{n}).
\label{eq5.3}
\end{equation}
By the bigness of $L$, we have
\begin{equation}
\dim H^{0}(X,L^{k})\geq Ck^{n},\quad k\rightarrow \infty .
\label{eq5.4}
\end{equation}
Under the assumption (B), it follows that $M$ satisfies $Z(j)$ for all
$0\leq j\leq n-q-1$, thus {Theorem~\ref{neumannop}}(d) implies that, for
$q\leq \ell $,
\begin{equation}
\dim \mathscr{H}^{0,n-\ell -1}({M,E^{*}})< \infty .
\label{eq5.5}
\end{equation}
We obtain $Ck^{n}\leq o (k^{n})$ for $k$ large which cannot hold. Thus,
the assertion is true.
\end{proof}

\begin{lemma}%
\label{lemdbar}
Assume $X$ is compact and the assumptions (A) and (B) hold. For
$q \leq \ell \leq n-1$, let $s\in H^{0}(X,L^{k_{0}})$ be the non-zero holomorphic
section in {Lemma~\ref{propprojection}} and
$\sigma \in \Omega ^{n,\ell}(bM, E)$ be $\overline \partial _{b}$-closed.
Then, there exists a $\overline \partial $-closed extension $S$ of
$s\otimes \sigma \in \Omega ^{n,\ell}(bM, L^{k_{0}}\otimes E)$, i.e.,
$S\in \Omega ^{n,\ell}(\overline{M}, L^{k_{0}}\otimes E)$ such that
$\overline \partial S=0$ on $M$ and
$\mu (S|_{bM})=\mu (s\otimes \sigma )$.
\end{lemma}

\begin{proof}
Let $\sigma '\in \Omega ^{n,\ell}(\overline{M}, E)$ such that
$\sigma '|_{bM}=\sigma $. Let
$\theta \in \mathscr{H}^{0,n-\ell -1}({M},L^{-k_{0}}\otimes E^{*})$. Thus
\begin{equation}
\theta \wedge (s\sigma ')\in \Omega ^{n,n-1}(\overline{M})
\label{eq5.6}
\end{equation}
and $s\theta \in \Omega ^{0,n-\ell -1}(\overline{M},E^{*})$ is
$\overline \partial $-closed. By {Theorem~\ref{neumannop}} (a)(b) and
$H(s\theta )=0$, there exists
$\xi \in \Omega ^{0,n-\ell -2}(\overline{M},E^{*})$ such that
\begin{equation}
s\theta =\overline \partial \xi + H(s\theta )=\overline \partial \xi .
\label{eq5.7}
\end{equation}
Then, by the Stokes formula,
\begin{equation}
\label{eqstokes1}
\int _{bM}\theta \wedge (s\sigma ) = \int _{bM}(s\theta )\wedge
\sigma =\int _{bM}(\overline \partial \xi ) \wedge \sigma =(-1)^{n-
\ell -1}\int _{bM}\xi \wedge (\overline \partial \sigma ').
\end{equation}

By {(\ref{eqdecomp})}, {(\ref{eqdbarb})} and
$\overline \partial _{b} \sigma =0$, we have
\begin{equation}
(\overline \partial \sigma ')|_{bM}=\overline \partial _{b} \sigma +(
\psi \wedge \overline \partial \varrho )|_{bM} =(\psi \wedge
\overline \partial \varrho )|_{bM}
\label{eq5.9}
\end{equation}
where $\psi \in \Omega ^{n,\ell}(\overline{M},E)$, and $\varrho $ is the
defining function of $M$. By the Stokes formula and
$bM=\{x\in X: \varrho (x)=0\}$, we have
\begin{equation}
\label{eqstokes2}
\int _{bM}\xi \wedge (\overline \partial \sigma ') = \int _{bM}\xi
\wedge \psi \wedge \overline \partial \varrho =-\int _{bM}
\overline \partial (\xi \wedge \psi )\wedge \varrho =0.
\end{equation}

Then, by {(\ref{eqstokes1})} and {(\ref{eqstokes2})}, it follows that
\begin{equation}
\int _{bM}\theta \wedge (s\sigma )=0
\label{eq5.11}
\end{equation}
for the $\overline \partial _{b}$-closed form
$s\sigma \in \Omega ^{n,\ell}(bM, L^{k_{0}}\otimes E)$ and any
$\theta \in \mathscr{H}^{0,n-\ell -1}(M,L^{-k_{0}}\otimes E^{*})$. Finally,
there exists a $\overline \partial $-closed extension $S$ of
$s\sigma $ by {Theorem~\ref{extensioncr}} for $Z(n-\ell -1)$.
\end{proof}

\begin{proof}[Proof of {Theorem~\ref{thmextension1}}]
The result follows by replacing $E$ by $E\otimes \Lambda ^{n} T^{1,0}X$ in {Lemma \ref{propprojection}} and {\ref{lemdbar}}.
\end{proof}

{Theorem~\ref{thmextension1}} can be slightly generalized to the following
case without strict positivity.

\begin{thm}%
\label{thmextension2}
Let $X$ be a Moishezon manifold of $\dim X=n \geq 2$. Let
$1\leq q\leq n-1$ and let $M\Subset X$ be a domain with smooth boundary $bM$ on which
the Levi form restricted to the analytic tangent space has at least
$n-q$ negative eigenvalues. Let $L$ be a big line bundle on $X$ with a
singular Hermitian metric $h^{L}$. Assume that $h^{L}$ is smooth on
$M$ and $c_{1}(L,h^{L})\geq 0$ on $M$. Then the conclusion of {Theorem~\ref{thmextension1}} holds.
\end{thm}
\begin{proof}
Since $X$ is compact and $L$ is big, we have
$\dim H^{0}(X,L^{k})\geq C k^{n}$. Since the assumption (B) holds, {Theorem~\ref{neumannop}} (d) entails
$\dim \mathscr{H}^{0,n-\ell -1}(M,F)<\infty $. Meanwhile, {Corollary~\ref{cor_qconcave}} implies
$\dim \mathscr{H}^{0,n-\ell -1}(M,L^{-k}\otimes E^{*})=o(k^{n})$. As in
{Lemma~\ref{propprojection}}, for $q \leq \ell \leq n-1$ and the harmonic
projection
$H: H^{0}(X,L^{k})\times \mathscr{H}^{0,n-\ell -1}({M},L^{-k}\otimes E^{*})
\rightarrow \mathscr{H}^{0,n-\ell -1}({M,E^{*}})$, there exist
$k_{0}\in \mathbb{N}$ and $s\in H^{0}(X,L^{k_{0}})\setminus \{0\}$ such
that $H(s\theta )=0$ for every
$\theta \in \mathscr{H}^{0,n-\ell -1}({M},L^{-k_{0}}\otimes E^{*})$. Finally,
as in {Lemma~\ref{lemdbar}}, we obtain $\overline \partial $-closed extension
$S\in \Omega ^{n,\ell}(\overline{M}, L^{k_{0}}\otimes E)$, and then replace
$E$ by $E\otimes \Lambda ^{n} T^{1,0}X$.
\end{proof}

\section{Vanishing theorems and applications}
\label{sec-vanish}

In this section, we present some vanishing theorems as consequences of
the Bochner-Kodaira-Nakano formula {\eqref{eq:3.42Nov}} with boundary term.
These include the classical Andreo\-tti-Tomassini vanishing theorem and a
specific case for negative line bundles on domains with Levi-flat boundaries.

\subsection{Andreotti-Tomassini vanishing theorem;
Proof of {Theorem~\ref{thm_vanishing}}}
\label{ss:6.1AT}

In this subsection we revisit the Andreotti-Tomassini vanishing theorem
\cite{AT:69} (especially for the $q$-concave case), and give a proof based
on {Theorem~\ref{thm:NG-inequality}} and inspired by Griffiths
\cite[Theorem 7.1]{Griffiths1966}.

\begin{proof}[Proof of {Theorem~\ref{thm_vanishing}}]
The proof of {Theorem~\ref{thm_vanishing}} follows from a modification of
the proof of {Proposition~\ref{prop:4.3Jan25}}. We solely prove (a) via {Theorem~\ref{thm:NG-inequality}}, as an analogous proof of (b) can be achieved by
applying {\eqref{eq:3.13Nov}} (refer to \cite[Theorem 3.5.9]{MM}).

Let $\varrho $ be a defining function for $M$ such that
$M:=\{x\in X\;:\; \varrho (x)<0\}$ and $d\varrho $ doesn't vanish near
$bM$. First, by {Lemma~\ref{lem_metric}}, since
$\mathscr{L}_{\varrho}$ has at least $n-q$ negative eigenvalues along
$bM$, there exists a smooth Hermitian metric $\Theta $ in a neighborhood
of $\overline{M}$ such that $d\varrho $ has norm $1$ in a neighborhood
of $bM$ and for a sufficiently small $\varepsilon >0$
\begin{equation}
[\sqrt{-1}\mathscr{L}_{\varrho}, \Lambda ]^{bM}_{\varepsilon }\geq 1
\label{eq:6.3Jan25}
\end{equation}
as a positive operator acting on
$\oplus _{j\leq n-q-1} \Lambda ^{j}T^{(0,1)\ast}bM$.

Since $c_{1}(L,h^{L})<0$ on $\overline{M}$, we need not modify further
the metric $h^{L}$ for there exists $C'>0$ such that for all
$s\in \Omega ^{0,j}(\overline{M},L)$ ($j\leq n-1$) and
$x\in \overline{M}$,
\begin{equation}
\langle [\sqrt{-1}R^{L},\Lambda ]s,s\rangle _{h}(x)\geq C'|s|^{2}_{h}(x).
\label{eq:6.4Jan25}
\end{equation}

Then following the first part of the proof of {Lemma~\ref{lem-question2}}, by replacing
$\sqrt{-1}\partial \overline \partial \varrho $ with
$\sqrt{-1}R^{L}$, we can conclude that, after modifying $\Theta $ again
with the parameter $T\geq 1$ as in {\eqref{eq:4.20dec24}}, we may always
assume that for the Hermitian metric $\Theta $, {\eqref{eq:6.3Jan25}} still
holds as in {\eqref{eq:4.24dec24}} and moreover, there exists a constant
$C_{0}>0$ such that for any sufficiently small $\varepsilon >0$, and
for any $\tau \in [0,1]$,
$s\in \Omega _{0}^{0,j}(\overline{M},L^{k})$ with $\supp \,(s)$ included
in a small neighborhood of $bM$, with $j\leq n-q-1$ and any
$k\in \mathbb{N}_{\geq 1}$, we have
\begin{equation}
\label{eq-q2-new}
\,\left \langle [\sqrt{-1}R^{L},\Lambda ]_{\varepsilon} s,s\right
\rangle _{h}+ 2(1-\varepsilon )\tau \left \langle R^{L}(e^{(1,0)}_{
\mathfrak{n}}, e^{(0,1)}_{\mathfrak{n}})s, s\right \rangle _{h}\geq C_{0}|s|_{h}^{2}.
\end{equation}
When we work for the points in $M$ that are away from $bM$,
{\eqref{eq:6.4Jan25}} still holds no matter how we modify the metric
$\Theta $ (with possibly a different constant $C'>0$). So if
$\varepsilon >0$ is sufficiently small, we also have
\begin{equation}
\langle [\sqrt{-1}R^{L},\Lambda ]_{\varepsilon }s,s\rangle _{h}\geq C'|s|^{2}_{h}.
\label{eq:6.4Jan25-new}
\end{equation}

By {Theorem~\ref{thm:NG-inequality}}, combining {\eqref{eq:6.3Jan25}},
{\eqref{eq:6.4Jan25}}, {\eqref{eq-q2-new}}, and {\eqref{eq:6.4Jan25-new}}, we
conclude that there exist $C>0$, $C_{0}>0$, $C^{\prime }>0$ and a
sufficiently small $\varepsilon \in (0,1)$ such that for all
$k\in \mathbb{N}_{\geq 1}$,
$s\in \Omega ^{0,j}(\overline{M},L^{k}\otimes E)$ with $j\leq n-q-1$ and
$\iota _{e_{\mathfrak{n}}}s=0$ on $bM$, we have
\begin{equation}
\begin{split} (1+\varepsilon )\left (\|\overline \partial ^{E}_{k} s\|^{2}_{L^{2}(M)}+
\|\overline \partial ^{E*}_{k} s\|^{2}_{L^{2}(M)}\right )+ C(1+
\frac{1}{\varepsilon})\|s\|^{2}_{L^{2}(M)} \geq \left (C_{0} k - C^{
\prime } \right )\|s\|^{2}_{L^{2}(M)}.
\end{split}
\label{eq:6.5Jan2025}
\end{equation}

The above inequality still holds for
$s\in \Dom (\overline \partial ^{E}_{k})\cap \Dom (
\overline \partial ^{E*}_{k})\cap L^{2}_{0,j}(M,L^{k}\otimes E)$ for
$j\leq n-q-1$, by the density of $B^{0,j}(M,L^{k}\otimes E)$ with respect
to the graph norm of
$\overline \partial ^{E}_{k}+\overline \partial ^{E*}_{k}$. Therefore,
we conclude that if
$k> \frac{1}{C_{0}}(C^{\prime }+C+C/\varepsilon )$, for all
$j\leq n-q-1$
\begin{equation}
\label{eq_vanxc}
\mathscr{H}^{0,j}(M,L^{k}\otimes E)=0.
\end{equation}
Using {Theorem~\ref{zqzq1}} to identify
$\mathscr{H}^{0,j}(M,L^{k}\otimes E)$ with
$H^{0,j}(\overline{M},L^{k}\otimes E)$ for $j\leq n-q-1$ completes the
proof.
\end{proof}

As a consequence of the above theorem, we recover the Andreotti-Tomassini
vanishing theorem \cite{AT:69}.
\begin{thm}[{\cite[Theorems 1 and 2]{AT:69}}]%
\label{vanish}
Let $n\geq 2$ and $1\leq q\leq n-1$.

(a) Let $X$ be a connected $q$-concave complex manifold of dimension
$n$, and let $\varrho $ be an exhaustion function with an exceptional set
$K\Subset X$. Let $E$ and $L$ be holomorphic vector bundles on $X$ and
$\rank (L)=1$. Suppose $L$ is negative in the neighborhood of
$\overline{X}_{c}$ for some regular value $c$ of $\varrho $ such that
$K\subset X_{c}$. Then, there exists an integer $k_{0}=k_{0}(L,E)$ such
that
\begin{equation}
H^{\ell}(X,L^{k}\otimes E)=0 \quad \,\text{ for }\, k\geq k_{0}\,
\text{ and }\,\ell \leq n-q-1.
\label{eq:6.7AT}
\end{equation}

(b) Let $X$ be a connected $q$-convex complex manifold of dimension
$n$, and let $\varrho $ be an exhaustion function with an exceptional set
$K\Subset X$. Let $E$ and $L$ be holomorphic vector bundles on $X$ and
$\rank (L)=1$. Suppose $L$ is positive in the neighborhood of
$\overline{X}_{c}$ for some regular value $c$ of $\varrho $ such that
$K\subset X_{c}$. Then there exists an integer $k_{0}=k_{0}(L,E)$ such
that
\begin{equation}
H^{\ell}(X,L^{k}\otimes E)=0 \quad \,\text{ for }\,k\geq k_{0}\,
\text{ and }\,\ell \geq q.
\label{eq6.8}
\end{equation}
\end{thm}
\begin{proof}
For (a): by {Theorem~\ref{thm_vanishing}}(a), for all sufficiently large
$k$, we have for $j\leq n-q-1$,
\begin{equation}
H^{0,j}(\overline{X}_{c}, L^{k}\otimes E)=0.
\label{eq6.9}
\end{equation}
By {Theorems~\ref{zqzq1}, \ref{thm:2.9AG} and \ref{thm:2.12H65}}, for
$j\leq n-q-2$, we have
\begin{equation}
H^{j}(X,L^{k}\otimes E)\simeq H^{0,j}(\overline X_{c},L^{k}\otimes E),
\label{eq6.10}
\end{equation}
and 
\begin{equation}
\dim H^{n-q-1}(X,L^{k}\otimes E)\leq \dim H^{0,n-q-1}(\overline X_{c},L^{k}
\otimes E).
\label{eq6.11}
\end{equation}
So we obtain {\eqref{eq:6.7AT}}.

Part (b) follows from {Theorem~\ref{thm_vanishing}}(b), {Theorem~\ref{zqzq1}}, and from theorems of H\"{o}rmander and Andreotti-Grauert
\cite[Theorem 3.5.6 and 3.5.7]{MM}.
\end{proof}

\subsection{The \texorpdfstring{$\overline \partial _{b}$}{dbar	extsubscript{b}} extension via vanishing theorems}
\label{sec_related}

If we replace the assumption of {Theorem~\ref{thmextension1}} that $L$ is
semi-positive on $X$ and positive at least at one point with the assumption
$L>0$ on $M$, the $\overline \partial _{b}$ extension result follows via
our proof of {Theorem~\ref{thm_vanishing}}(a). In this case, the compactness
assumption on $X$ can be dropped.

\begin{thm}%
\label{extpos}
Let $X$ be a complex manifold of dimension $n$. Let $(E,h^{E})$ and
$(L,h^{L})$ be holomorphic Hermitian vector bundles over $X$ and
$\rank (L)=1$. Let $1\leq q\leq n-1$ and assume that $M$ is a relatively
compact domain in $X$ with a smooth boundary $bM$ on which the Levi form
, restricted to the analytic tangent space, has at least $n-q$ negative eigenvalues.
Assume that $c_{1}(L,h^{L})>0$ on $\overline{M}$.

Then there exists $k_{0}\in \mathbb{N}$ such that for each
$k\geq k_{0}$, every $q\leq \ell \leq n-1$, and every
$\overline \partial _{b}$-closed form
$\sigma \in \Omega ^{0,\ell}(bM,L^{k}\otimes E)$, there exists a
$\overline \partial $-closed extension $\Sigma $ of $\sigma $, i.e.,
\begin{equation}
\Sigma \in \Omega ^{0,\ell}(\overline{M}, L^{k}\otimes E)
\label{eq6.12}
\end{equation}
such that $\overline \partial \Sigma =0$ on $M$ and
$\mu (\Sigma |_{bM})=\mu (\sigma )$ on $bM$.
\end{thm}
\begin{proof}
Firstly, $M$ satisfies $Z(i)$ for all $0\leq i\leq n-q-1$. By {\eqref{eq_vanxc}} for $i\leq n-q-1$ and $k$ large,
$\mathscr{H}^{n,i}(M,L^{-k}\otimes E^{*})=0$ after identifying
$(n,i)$-forms with $(0,i)$-forms valued in $K_X$. Finally, we apply {Theorem~\ref{extensioncr}} for $Z(n-\ell -1)$ with $\ell \geq q$.
\end{proof}

Similarly, if $L<0$ on $M$, the $\overline \partial _{b}$ extension result
follows via {Theorem~\ref{thm_vanishing}}(b).

\begin{thm}%
\label{extneg}
Let $X$ be a complex manifold of dimension $n$. Let $(E,h^{E})$ and
$(L,h^{L})$ be holomorphic Hermitian bundles over $X$ and
$\rank (L)=1$. Let $1\leq q\leq n-1$ and assume that $M$ is a relatively
compact domain in $X$ with a smooth boundary $bM$ on which the Levi form
restricted to the analytic tangent space has at least $q$ positive eigenvalues.
Assume that $c_{1}(L,h^{L})<0$ on $\overline M$.

Then there exists $k_{0}\in \mathbb{N}$ such that for each
$k\geq k_{0}$, every $0\leq \ell \leq q-1$, and every
$\overline \partial _{b}$-closed form
$\sigma \in \Omega ^{0,\ell}(bM,L^{k}\otimes E)$, there exists a $\overline \partial $-closed extension $\Sigma $ of $\sigma $, i.e.,
\begin{equation}
\Sigma \in \Omega ^{0,\ell}(\overline{M}, L^{k}\otimes E)
\label{eq6.13}
\end{equation}
such that $\overline \partial \Sigma =0$ on $M$ and
$\mu (\Sigma |_{bM})=\mu (\sigma )$ on $bM$.
\end{thm}

\begin{proof}
Firstly, $M$ satisfies $Z(j)$ for $n-q\leq j\leq n-1$. Since
$M=\{\varrho<0\}$ and $M_{\varepsilon}=\{ \varrho<\varepsilon \}$ for sufficiently
small $\varepsilon >0$ are $j$-convex manifolds for $j\geq n-q$, then  {Theorem~\ref{thm_vanishing}}(b) with the replacement of the degree
$q$ by $n-q$ gives
$\mathscr{H}^{n,j}(M,L^{-k}\otimes E^{*})=0$ for $k$ large after identifying
$(n,j)$-forms with $(0,j)$-forms valued in $K_X$. Finally, we apply {Theorem~\ref{extensioncr}} for $Z(n-\ell -1)$ with $0\leq \ell \leq q-1$.
\end{proof}

{Theorem~\ref{extneg}} with $q=1$ means that, under the assumption
$L<0$ on $\overline M$, we can generalize the classical result
\cite[Theorem 7.5]{KR:65} to forms with values in $L^{k}\otimes E$.

As a final remark, we proceed in the same vein by using {Theorem~\ref{extensioncr}} for $q=p=0$ and \cite[Proposition 4.17]{Wh:19} that
$\mathscr{H}^{n,n-1}(M,L^{*})=H^{0,n-1}_{(2)}(M,L^{*}\otimes K_{X})=0$.
It follows that: Let $(X,\omega )$ be a K\"{a}hler manifold and
$M\Subset X$ be a domain with smooth boundary $bM$ such that the Levi form on $bM$ have
one positive eigenvalue and semi-positive everywhere. If
$(L,h^{L})$ is negative on $bM$ and semi-negative on $M$, then every
$\overline \partial _{b}$-closed smooth section of $L$ on $bM$ has a holomorphic
extension to all of $M$.

\subsection{Negative line bundles on domains with Levi-flat boundary}
\label{ss:6.3Leviflat}

A domain $M=\{x\in X: \varrho (x)<0\}\Subset X$ with defining function
$\varrho \in \mathscr{C}^{\infty}(X,\mathbb{R})$ is said to have a Levi-flat
boundary if the Levi form $\mathscr{L}_{\varrho}$ vanishes identically
along $bM$. In the Levi-flat case, the boundary term on the right-hand
side of {\eqref{eq:3.44Nov}} vanishes, allowing us to obtain the following
result for the negative line bundles.

\begin{prop}%
\label{prop:6.4may}
Let $X$ be a complex manifold of dimension $n\geq 2$, and let
$M\Subset X$ be a domain with a smooth Levi-flat boundary $bM$. Let
$\varrho $ be a defining function of $M$, and let $\Theta _{0}$ be a smooth
Hermitian metric on $X$ such that $|d\varrho |=1$ near $bM$. Let
$(L,h^{L})$ be a negative line bundle in a neighborhood of
$\overline{M}$ and let $(E,h^{E})$ be any holomorphic Hermitian vector bundle on
$X$. Then there exists a Hermitian metric $\Theta $ on
$\overline{M}$, which is given by modifying $\Theta _{0}$ locally near
$bM$, and an integer $k_{0}>0$ such that for $q=0,\ldots , n-1$ and
$k\geq k_{0}$, we have
\begin{equation}
\mathscr{H}^{0,q}_{(2)}(M, L^{k}\otimes E; \Theta ):= \ker \square ^{L^{k}
\otimes E}_{k}|_{L^{2}_{0,q}(M, L^{k}\otimes E)}=0,
\label{eq6.14}
\end{equation}
where the $L^{2}$-cohomology and Kodaira Laplacians are given with respect
to the metric $\Theta $.
\end{prop}
\begin{proof}
Using the fact that $L$ is negative near $\overline{M}$, we can proceed
as in the proof of {Lemma~\ref{lem-question2}} to modify $\Theta _{0}$ near
$bM$ such that $d\varrho $ still has norm $1$ near $bM$. Moreover, there
exists $C_{0}>0$ such that for a small $0<\varepsilon \ll 1$, any
$\tau \in [0,1]$, and for any $s\in \Omega ^{0,j}(\overline M, L^{k}\otimes E)$ (with
$j\leq n-1$ and $k\in \mathbb{N}$), we have near $bM$,
\begin{equation}
\label{eq:Leviflat}
\,\left \langle [\sqrt{-1}R^{L},\Lambda ]_{\varepsilon} s,s\right
\rangle _{h}+2(1-\varepsilon )\tau \left \langle R^{L}(e^{(1,0)}_{
\mathfrak{n}}, e^{(0,1)}_{\mathfrak{n}})s,s\right \rangle _{h}\geq C_{0}|s|_{h}^{2},
\end{equation}
where the metric is taken with respect to $\Theta $ and
$(h^{L})^{\otimes k}$. Since $bM$ is Levi-flat and $d\varrho $ has the
norm $1$ near $bM$, the boundary terms vanish identically in {\eqref{eq:3.44Nov}}. Finally, we repeat the proof of {Proposition~\ref{prop:4.3Jan25}}. In particular, by applying {Theorem~\ref{thm:NG-inequality}}, we obtain, with some constants $C>0$ and
$C'>0$ independent of $k$,
\begin{equation}
(1+\varepsilon )\left (\|\overline \partial ^{E}_{k} s\|^{2}_{L^{2}(M)}+
\|\overline \partial ^{E*}_{k} s\|^{2}_{L^{2}(M)}\right )+C\left (1+
\frac{1}{\varepsilon}\right )\|s\|^{2}_{L^{2}(M)}\geq \left (C_{0} k -
C' \right )\|s\|^{2}_{L^{2}(M)}.
\label{eq:6.16Levi}
\end{equation}
The proposition follows.
\end{proof}

As a consequence, for any domain $M$ with a smooth Levi-flat boundary ($
\dim _{\mathbb{C}}M\geq 2$), there exists no strictly plurisubharmonic
bounded smooth function in any open neighborhood of $\overline{M}$. This
also follows from the maximum principle directly applied on the boundary
of $M$. When we specialize to projective spaces, we conclude from {Proposition~\ref{prop:6.4may}} that if $M\Subset \mathbb{CP}^{n}$ is a domain with a
smooth Levi-flat boundary, then each hyperplane in $\mathbb{CP}^{n}$ has
a nonempty intersection with $\overline{M}$. This also follows from the
solution to the Levi problem over $\mathbb{CP}^{n}$ by Fujita
\cite{MR159034} and Takeuchi \cite{MR173789}. Namely, any proper pseudoconvex
domain in $\mathbb{CP}^{n}$ is Stein and neither $M$ nor the complement of
$\overline{M}$ can contain a divisor when the boundary of $M$ is Levi-flat.
Furthermore, for $n \geq 3$, this is evidently due to the nonexistence
of smooth Levi-flat real hypersurfaces in $\mathbb{CP}^{n}$
\cite{Siu2000}.

The conjecture of the nonexistence of smooth real Levi-flat hypersurfaces
in projective spaces first appeared in \cite{MR1001454} and
\cite{MR1275208} for the study of minimal sets of holomorphic foliations
on $\mathbb{CP}^{n}$. Then it was proven for $n\geq 3$ by Lins Neto
\cite{MR1703092} in the real analytic case and by Siu \cite{Siu2000} in
the smooth case. Despite many attempts, the nonexistence of smooth Levi-flat
real hypersurfaces in $\mathbb{CP}^{2}$ is still considered an open problem;
see \cite{CS07,IM2008,MR3587467,Shafikov:2025aa} and references therein. In
\cite{AdachiBri2015}, Adachi and Brinkschulte deduced a curvature restriction
of Levi-flat real hypersurfaces in $\mathbb{CP}^{2}$. Their method was
inspired by the integral formula of Griffiths
\cite[VII.\S 5]{Griffiths1966}; see also {Lemma~\ref{lm:7.3}} and {Corollary~\ref{cor:3.5new}}.


\begin{thebibliography}{{CMW}23}
\bibliographystyle{apalike}	

\bibitem[AB15]{AdachiBri2015}
M.~Adachi and J.~Brinkschulte.
\newblock Curvature restrictions for {L}evi-flat real hypersurfaces in complex
  projective planes.
\newblock {\em Ann. Inst. Fourier (Grenoble)}, 65(6):2547--2569, 2015.



	\bibitem[And73]{And:73}
	A. Andreotti.
	\newblock Nine lectures on complex analysis.
	\newblock In {\em Complex Analysis}, volume~62 of {\em C.I.M.E. Summer Schools
		(1973)}, pages 1--166. Springer-Verlag Berlin Heidelberg, 2011.

	\bibitem[AG62]{AG:62}
	A. Andreotti and H. Grauert.
	\newblock Th\'eor\`eme de finitude pour la cohomologie des espaces complexes.
	\newblock {\em Bull. Soc. Math. France}, 90:193--259, 1962.
	
	\bibitem[AH72]{AH:72}
	A. Andreotti and C.~Denson Hill.
	\newblock E. {E}. {L}evi convexity and the {H}ans {L}ewy problem. {I}.
	{R}eduction to vanishing theorems.
	\newblock {\em Ann. Scuola Norm. Sup. Pisa (3)}, 26:325--363, 1972.

    \bibitem[AS70]{Andreotti-Siu}
A.~Andreotti and Y.-T. Siu.
\newblock Projective embedding of pseudoconcave spaces.
\newblock {\em Ann. Scuola Norm. Sup. Pisa Cl. Sci. (3)}, 24:231--278, 1970.


	
	\bibitem[AT69]{AT:69} 
	A.~Andreotti and G.~Tomassini.
	\newblock A remark on the vanishing of certain cohomology groups.
	\newblock {\em Compositio Math.}, 21:417--430, 1969.

\bibitem[AV65]{AV65}
A. Andreotti and E. Vesentini, Carleman estimates for the Laplace-Beltrami equation
on complex manifolds, Inst. Hautes Etudes Sci. Publ. Math., 25 (1965), 81--130.	
	
	
	\bibitem[Ber05]{RB:05}
	R. Berman.
	\newblock Holomorphic {M}orse inequalities on manifolds with boundary.
	\newblock {\em Ann. Inst. Fourier (Grenoble)}, 55(4):1055--1103, 2005.

    \bibitem[Ber08]{RB:08}
R.~Berman.
\newblock Corrigendum to: ``{H}olomorphic {M}orse inequalities on manifolds
  with boundary'' [{A}nn. {I}nst. {F}ourier ({G}renoble) {\bf 55} (2005), no.
  4, 1055--1103; mr2157164].
\newblock {\em Ann. Inst. Fourier (Grenoble)}, 58(1):377--381, 2008.

\bibitem[Ber10]{RB:10}
R. Berman.
\newblock Bergman kernels and equilibrium measures for polarized pseudo-concave
  domains.
\newblock {\em Internat. J. Math.}, 21(1):77--115, 2010.



	\bibitem[Boc43]{Bo:43}
	S.~Bochner.
	\newblock Analytic and meromorphic continuation by means of {G}reen's formula.
	\newblock {\em Ann. of Math. (2)}, 44:652--673, 1943.

    \bibitem[Bou89]{Bouche1989}
T.~Bouche.
\newblock In\'{e}galit\'{e}s de {M}orse pour la {$d''$}-cohomologie sur une
  vari\'{e}t\'{e} holomorphe non compacte.
\newblock {\em Ann. Sci. \'{E}cole Norm. Sup. (4)}, 22(4):501--513, 1989.


\bibitem[CLS89]{MR1001454}
C.~Camacho, A.~Lins~Neto, and P.~Sad.
\newblock Minimal sets of foliations on complex projective spaces.
\newblock {\em Inst. Hautes \'{E}tudes Sci. Publ. Math.}, 68:187--203 (1989),1988.

\bibitem[CS07]{CS07}
J.~Cao and M.-C.~Shaw. 
\newblock{The $\overline\partial$-Cauchy problem and nonexistence
of Lipschitz Levi-flat hypersurfaces in
$\C\mathbb{P}^n$ with $n\geq3$.} 
\newblock {\em Math. Z.}, 256:175-192, 2007 .


\bibitem[Ce93]{MR1275208}
D.~Cerveau.
\newblock Minimaux des feuilletages alg\'{e}briques de {${\bf C}{\rm P}(n)$}.
\newblock {\em Ann. Inst. Fourier (Grenoble)}, 43(5):1535--1543, 1993.


	
	\bibitem[CMW23]{CMW23}
	D. {Coman}, G. {Marinescu}, and H. Wang.
	\newblock Singular holomorphic Morse inequalities on non-compact manifolds. \newblock {\em Rev. Roum. Math. Pures Appl.} No. 1-2 (2023).

    \bibitem[D85]{Demailly1985}
J.-P. {Demailly}.
\newblock Champs magn\'{e}tiques et in\'{e}galit\'{e}s de {M}orse pour la
  {$d''$}-cohomologie.
\newblock {\em Ann. Inst. Fourier (Grenoble)}, 35(4):189--229, 1985.



\bibitem[D86]{Dem1986}
J.-P. Demailly.
\newblock Sur l'identit\'{e} de {B}ochner-{K}odaira-{N}akano en
  g\'{e}om\'{e}trie hermitienne.
\newblock In {\em S\'{e}minaire d\/'analyse {P}. {L}elong-{P}. {D}olbeault-{H}.
  {S}koda, ann\'{e}es 1983/1984}, volume 1198 of {\em Lecture Notes in Math.},
  pages 88--97. Springer, Berlin, 1986.

  \bibitem[Eh56]{Ehr:56}
L.~Ehrenpreis.
\newblock Sheaves and differential equations.
\newblock {\em Proc. Amer. Math. Soc.}, 7:1131--1138, 1956.




	\bibitem[FK72]{FK:72}
	G.~B. Folland and J.~J. Kohn.
	\newblock {\em The {N}eumann problem for the {C}auchy-{R}iemann complex}.
	\newblock Princeton University Press, Princeton, N.J.; University of Tokyo
	Press, Tokyo, 1972.
	\newblock Annals of Mathematics Studies, No. 75.

    \bibitem[Fu63]{MR159034}
R.~Fujita.
\newblock Domaines sans point critique int\'{e}rieur sur l'espace projectif
  complexe.
\newblock {\em J. Math. Soc. Japan}, 15:443--473, 1963.


	
\bibitem[Gr66]{Griffiths1966}
P.~A. Griffiths.
\newblock The extension problem in complex analysis. {II}. {E}mbeddings with
  positive normal bundle.
\newblock {\em Amer. J. Math.}, 88:366--446, 1966.


\bibitem[HL88]{HL88}
G.~M. Henkin and J.~Leiterer.
\newblock {\em Andreotti-{G}rauert theory by integral formulas}, volume~74 of
  {\em Progress in Mathematics}.
\newblock Birkh\"{a}user Boston, Inc., Boston, MA, 1988.



\bibitem[H65]{Hor:65}
L. H\"ormander, $L^2$-estimates and existence theorem for the 
$\overline\partial$operator, Acta
Math. 113 (1965), 89--152.
	
	\bibitem[H90]{Hor:90}
	L. H\"ormander.
	\newblock {\em An introduction to complex analysis in several variables},
	volume~7 of {\em North-Holland Mathematical Library}.
	\newblock North-Holland Publishing Co., Amsterdam, third edition, 1990.

\bibitem[HLM22]{HLM22}
C.-Y.~Hsiao, X.~Li and G.~Marinescu.
\newblock{Semi-classical asymptotics of partial Bergman kernels on $\R$-symmetric complex manifolds with boundary}.
\newblock{arXiv:2208.12412v2.}

\bibitem[HMZ26]{HMZ26}
C.-Y.~Hsiao, G.~Marinescu and W.~Zhu.
\newblock{Semi-classical heat kernel asymptotics on complex manifolds with boundary}.
\newblock{arXiv: 2606.18689.}


    \bibitem[IM08]{IM2008}
A.~Iordan and F.~Matthey.
\newblock R\'{e}gularit\'{e} de l'op\'{e}rateur {$\overline\partial$} et
  th\'{e}or{\`e}me de {S}iu sur la non-existence d'hypersurfaces {L}evi-plates
  dans l'espace projectif.
\newblock {\em C. R. Math. Acad. Sci. Paris}, 346(7-8):395--400, 2008.

    
	\bibitem[KR65]{KR:65}
	J.~J. Kohn and H. Rossi.
	\newblock On the extension of holomorphic functions from the boundary of a
	complex manifold.
	\newblock {\em Ann. of Math. (2)}, 81:451--472, 1965.
	
	\bibitem[LTL89]{TL:89}
	C. Laurent-Thi\'ebaut and J. Leiterer.
	\newblock On the Hartogs-Bochner extension phenomenon for differential
	forms.
	\newblock {\em Math. Ann.}, 284(1):103--119, 1989.  
	
	\bibitem[LSW23]{LSW23}
	X. Li, G. Shao and H. Wang.
	\newblock 
	On Bergman kernel functions and weak holomorphic Morse inequalities.
	\newblock
	{\em Anal. Math. Phys.}13(2023), no.5, Paper No. 68, 22 pp.

    \bibitem[LN99]{MR1703092}
A.~Lins~Neto.
\newblock A note on projective {L}evi flats and minimal sets of algebraic
  foliations.
\newblock {\em Ann. Inst. Fourier (Grenoble)}, 49(4):1369--1385, 1999.

\bibitem[MM07]{MM}
	X. Ma and G. Marinescu.
	\newblock {\em Holomorphic {M}orse inequalities and {B}ergman kernels}, volume
	254 of {\em Progress in Mathematics}.
	\newblock Birkh\"auser Verlag, Basel, 2007.
	
\bibitem[M92]{Mar92b}
{G.~Marinescu}, {Morse inequalities for {$q$}-positive line bundles
  over weakly {$1$}-complete manifolds}, {\em C. R. Acad. Sci. Paris S\'er. I Math.},
  315 (1992), pp.~895--899.
	
	\bibitem[M96]{M:96}
	G. Marinescu.
	\newblock Asymptotic {M}orse inequalities for pseudoconcave manifolds.
	\newblock {\em Ann. Scuola Norm. Sup. Pisa Cl. Sci. (4)}, 23(1):27--55, 1996.

    	\bibitem[M97]{M97} 
	G. Marinescu.
	\newblock {Some applications of the Kohn-Rossi extension theorems}.
	\newblock {\em arXiv:math/9705203}.



    \bibitem[M16]{M:2016-q}
G.~Marinescu.
\newblock Existence of holomorphic sections and perturbation of positive line
  bundles over {$q$}-concave manifolds.
\newblock {\em Bull. Inst. Math. Acad. Sin. (N.S.)}, 11(3):579--602, 2016.


	
	

\bibitem[NT88]{NT88}
A. Nadel and H. Tsuji.
\newblock{Compactification of complete K\"ahler manifolds of negative Ricci curvature}, 
\newblock{\em{J. Differential Geom.}} 28 (1988), no. 3, 503-512.


\bibitem[Oh82]{Oh:82}T. Ohsawa, Isomorphism theorems for cohomology groups of weakly 1-
complete manifolds, Publ. Res. Inst. Math. Sci. 18 (1982), 191-232.

\bibitem[Oh15]{MR3587467}
T.~Ohsawa.
\newblock A survey on {L}evi flat hypersurfaces.
\newblock In {\em Complex geometry and dynamics}, volume~10 of {\em Abel
  Symp.}, pages 211--225. Springer, Cham, 2015.


\bibitem[Sh25]{Shafikov:2025aa}
R.~Shafikov.
\newblock Levi-flats in {$\mathbb{CP}^n$}: a survey for nonexperts.
\newblock {\em J. Geom. Anal.} 35 (2025) 181.



\bibitem[S84]{Siu84}
Y.-T. Siu.
\newblock
A vanishing theorem for semipositive line bundles over
  non--K\"ahler manifolds.
\newblock{\em J. Differential Geom.}, 20 (1984), pp.~431--452.




\bibitem[S00]{Siu2000}
Y.-T. Siu.
\newblock Nonexistence of smooth {L}evi-flat hypersurfaces in complex
  projective spaces of dimension {$\geq 3$}.
\newblock {\em Ann. of Math. (2)}, 151(3):1217--1243, 2000.

\bibitem[Ta64]{MR173789}
A.~Takeuchi.
\newblock Domaines pseudoconvexes infinis et la m\'{e}trique riemannienne dans
  un espace projectif.
\newblock {\em J. Math. Soc. Japan}, 16:159--181, 1964.



\bibitem[TCM01]{TCM:01}
\newblock{R.~Todor, I.~Chiose, and G.~Marinescu}, \newblock{Morse inequalities for
  covering manifolds}, \newblock{}{\em Nagoya Math. J.}, 163 (2001), 145--165.
    
	\bibitem[W21]{Wh:19}
	H. Wang.
	\newblock {Cohomology dimension growth for Nakano $q$-semipositive line
		bundles}.
	\newblock {\em J. of Geom. Anal.} 31 (2021), 4934--4965.


\end{thebibliography}
\end{document}